\documentclass[preprint]{imsart}

\RequirePackage[OT1]{fontenc}
\RequirePackage{amsthm,amsmath,amssymb}
\RequirePackage{natbib}
\RequirePackage[colorlinks,citecolor=blue,urlcolor=blue]{hyperref}
\usepackage{graphicx}
\usepackage{amsthm}
\usepackage{pgf,tikz}
\usepackage{url}
\usepackage{thmtools}
\usetikzlibrary{arrows}
\usepackage{pdfcomment}
\usepackage[final]{pdfpages}

\usepackage{diagbox}
\usepackage{multirow}
\usepackage{booktabs}
\usepackage{adjustbox}
\usepackage{hhline}

\numberwithin{equation}{section}
\theoremstyle{plain}
\declaretheorem[name=Theorem,numberwithin=section]{Th}
\declaretheorem[name=Lemma,numberwithin=section]{Lem}

\theoremstyle{definition}
\newtheorem{Def}{Definition}[section]

\theoremstyle{remark}

 \newtheorem{Remth}{Remark}[Th]

 \newtheorem{Rem}{Remark}[section]

\DeclareMathOperator{\E}{E}                       

\renewcommand{\P}{\mathit{P}}
\renewcommand{\i}{{\,\text{\it \@i}\,}}

\newcommand{\Tr}{{{\rm Trace}}}
\newcommand{\Normal}{\mathbf{N}}

\newcommand{\I}{{\rm I}}

\newcommand{\X}{{\mathbf{X}}}
\newcommand{\Y}{{\mathbf{Y}}}

\newcommand{\1}{{\mathbf{1}}}

\newcommand{\x}{{\textbf{\textit{x}}}}
\newcommand{\y}{{\textbf{\textit{y}}}}

\makeatother

\begin{document}

\begin{frontmatter}
\title{Statistical applications of random matrix theory:\\ comparison of two populations I\thanksref{T1}}
\runtitle{Comparison of two populations}
\thankstext{T1}{This paper is based on the PhD thesis of R\'emy Mari\'etan that will be divided into three parts.}

\begin{aug}
\author{\fnms{R\'emy} \snm{Mari\'etan}\thanksref{t1}\ead[label=e1]{remy.marietan@epfl.alumni.ch}}
\and
\author{\fnms{Stephan} \snm{Morgenthaler}\thanksref{t2}\ead[label=e2]{stephan.morgenthaler@epfl.ch}}

\thankstext{t1}{PhD Student at EPFL in mathematics department}
\thankstext{t2}{Professor at EPFL in mathematics department}
\runauthor{R\'emy Mari\'etan and Stephan Morgenthaler}

\affiliation{\'Ecole Polytechnique F\'ed\'eral de Lausanne, EPFL}

\address{Department of Mathematics\\
\'Ecole Polytechnique F\'ed\'eral de Lausanne\\1015 Lausanne\\
\printead{e1}, 
\printead*{e2}}
\end{aug}

\begin{abstract}
This paper investigates a statistical procedure for testing the equality of two independent estimated covariance matrices when the number of potentially dependent data vectors is large and proportional to the size of the vectors, that is, the number of variables.
Inspired by the spike models used in random matrix theory, we concentrate on the largest eigenvalues of the matrices in order to determine significance. To avoid false rejections we must guard against residual spikes and need a sufficiently precise description of the behaviour of the largest eigenvalues under the null hypothesis. 

In this paper, we lay a foundation by treating alternatives based on perturbations of order $1$, that is, a single large eigenvalue. Our statistic allows the user to test the equality of two populations. Future work will extend the result to perturbations of order $k$ and demonstrate conservativeness of the procedure for more general matrices.
\end{abstract}



\end{frontmatter}

\section{Introduction} 

In the last two decades, random matrix theory (RMT) has produced numerous results that offer a better understanding of large random matrices. These advances have enabled interesting applications in communication theory and even though it can potentially contribute to many other data-rich domains such as brain imaging or genetic research, it has rarely been applied.
The main barrier to the adoption of RMT may be the lack of concrete statistical results from the probability side. The straightforward adaptation of classical multivariate theory to high dimensions can sometimes be achieved, but such  procedures are only valid under strict assumptions about the data such as normality or independence. Even minor differences between the model assumptions and the actual data lead to catastrophic results and such procedures also often do not have enough power.

This paper proposes a statistical procedure for testing the equality of two covariance matrices when the number of potentially dependent data vectors $n$ and the number of variables $m$ are large. RMT denotes the investigation of estimates of covariance matrices $\hat{\Sigma}$ or more precisely their eigenvalues and eigenvectors when both $n$ and $m$ tend to infinity with $\lim \frac{m}{n}=c>0$.
When $m$ is finite and $n$ tends to infinity the behaviour of the random matrix is well known and presented in the books of \cite{multi3}, \cite{multi} and \cite{multi2} (or its original version \cite{multi22}).
In the RMT case, the behaviour is more complex, but many results of interest are known. \cite{Alice}, \cite{Tao} and more recently \cite{bookrecent} contain comprehensive introductions to RMT and \cite{Appliedbook} covers the case of empirical (estimated) covariance matrices.

Although the existing theory builds a good intuition of the behaviour of these matrices, it does not provide enough of a basis to construct a test with good power, which is robust with respect to the assumptions.
Inspired by the spike models, we define the residual spikes and provide a description of the behaviour of this statistic under a null hypothesis when the perturbation is of order $1$. These results enable the user to test the equality of two populations as well as other null hypotheses such as the independence of two sets of variables.
Later papers will extend the results to perturbations of order $k$ and demonstrate the robustness of our test's level for more general matrices.

The remainder of the paper is organized as follows. In the next section, we develop the test statistic and discuss the problems associated with high dimensions. Then we present the main theorem \ref{TH=Main}. Various results necessary for the proof are introduced in Section \ref{sec:Theorems}. The proofs themselves are technical and presented in the supplementary material \cite{Suppmaterial} included in the second part of this paper. The last section contains case studies and a comparison with alternative tests. 

\section{Test statistic}\label{section:teststat}

We compare the spectral properties of two covariance estimators $\hat{\Sigma}_X$ and $\hat{\Sigma}_Y$ of dimension $m\times m$ which can be represented as
\begin{eqnarray*}
\hat{\Sigma}_X=P_X^{1/2} W_X P_X^{1/2} \text{ and } \hat{\Sigma}_Y=P_Y^{1/2} W_Y P_Y^{1/2}.
\end{eqnarray*}
In this equation, $W_X$ and $W_Y$ are of the form
\begin{eqnarray*}
W_X=O_X \Lambda_X O_X \text{ and } W_Y=O_Y \Lambda_Y O_Y,
\end{eqnarray*}
with $O_X$ and $O_Y$ being independent unit orthonormal random matrices whose distributions are invariant under rotations, while  $\Lambda_X$ and $\Lambda_Y$ are independent positive random diagonal matrices, independent of  $O_X, O_Y$ with trace equal to m and a bound on the diagonal elements. Note that the usual RMT assumption, $\frac{m}{n}=c$ is replaced by this bound! The (multiplicative) spike model of order 1 determines the form of 
$P_X= \I_m + (\theta_{X}-1) u_X u_X^t$ and $P_Y= \I_m+ (\theta_{Y}-1) u_Y u_Y^t$. 

Our results will apply to any two centered data matrices $\X \in \mathbb{R}^{m\times n_X}$ and $\Y \in \mathbb{R}^{m\times n_Y}$ which are such that 
\begin{eqnarray*}
\hat{\Sigma}_X= \frac{1}{n_X}\X \X^t \text{ and } \hat{\Sigma}_Y= \frac{1}{n_Y}\Y \Y^t
\end{eqnarray*}
and can be decomposed in the manner indicated. This is the basic assumption concerning the covariance matrices. We will assume throughout the paper that $n_X\geq n_Y$. Note that because $O_X$ and $O_Y$ are independent and invariant by rotation we can assume without loss of generality that $u_X=e_1$ as in \cite{deformedRMT}. Under the null hypothesis, $P_X=P_Y$  and we use the simplified notation $P$ for both matrices where $\theta_{X}=\theta_{Y}=\theta$ and $u_X=u_Y(=e_1)$. 

To test $H_0:P=P_X=P_Y$ against $H_1:P_X\neq P_Y$ it is natural to consider the extreme eigenvalues of 
\begin{eqnarray}
\hat{\Sigma}_X^{-1/2} \hat{\Sigma}_Y \hat{\Sigma}_X^{-1/2}\,.\label{eq:base1}
\end{eqnarray}
We could also swap the subscripts, but it turns our to be preferable to use the inversion on the matrix with larger sample size. 

The distributional approximations we will refer to are based on RMT, that is, they are derived by embedding a given data problem into a sequence of random matrices for which both $n$ and $m$ tend to infinity such that $m/n$ tends to a positive constant $c$. The most celebrated results of RMT describe the almost sure weak convergence of the empirical distribution of the eigenvalues (spectral distribution) to a non-random compactly supported limit law. An extension of this theory to the "Spike Model" suggests that we should modify $\hat{\Sigma}$ because estimates of isolated eigenvalues derived from these estimates are asymptotically biased. The following corrections will be used. 
\begin{Def} \label{Def=unbiased} 
Suppose $\hat{\Sigma}$ is of the form described at the start of the section. 
The \textbf{unbiased estimator of }$\theta$ is defined as 
\begin{equation} \hat{\hat{\theta}}=1+\frac{1}{\frac{1}{m-1} \sum_{i=2}^{m} \frac{\hat{\lambda}_{\hat{\Sigma},i}}{\hat{\theta}-\hat{\lambda}_{\hat{\Sigma},i}}}\,,\label{eq:corrlambda}
\end{equation}
where $\hat{\lambda}_{\hat{\Sigma},i}$ is the $i^{\text{th}}$ eigenvalue of $\hat{\Sigma}$. When $\hat{\Sigma}=P^{1/2}WP^{1/2}$ as above, it is asymptotically equivalent to replace $\frac{1}{m-1}\sum_{i=2}^{m} \frac{\hat{\lambda}_{\hat{\Sigma},i}}{\hat{\theta}-\hat{\lambda}_{\hat{\Sigma},i}}$ by $\frac{1}{m}\sum_{i=1}^{m} \frac{\hat{\lambda}_{W,i}}{\hat{\theta}-\hat{\lambda}_W,i}$.\\
Suppose that $\hat{u}$ is the eigenvector corresponding to $\hat{\theta}$, then the \textbf{filtered estimated covariance matrix } is defined as 
\begin{equation}
\hat{\hat{\Sigma}}= \I_m+(\hat{\hat{\theta}}-1) \hat{u} \hat{u}^t\,.\label{eq:corrSigma}
\end{equation}
\end{Def}

The matrix (\ref{eq:base1}) which serves as the basis for the test then becomes either
\begin{equation}
\hat{\hat{\Sigma}}_X^{-1/2} \hat{\hat{\Sigma}}_Y \hat{\hat{\Sigma}}_X^{-1/2} 
\text{ or }\hat{\hat{\Sigma}}_X^{-1/2} \hat{\Sigma}_Y \hat{\hat{\Sigma}}_X^{-1/2} \,.\label{eq:base2}
\end{equation}

\begin{figure}[htbp] 
\centering
\begin{tabular}{cc}
 \includegraphics[width=0.26\paperwidth]{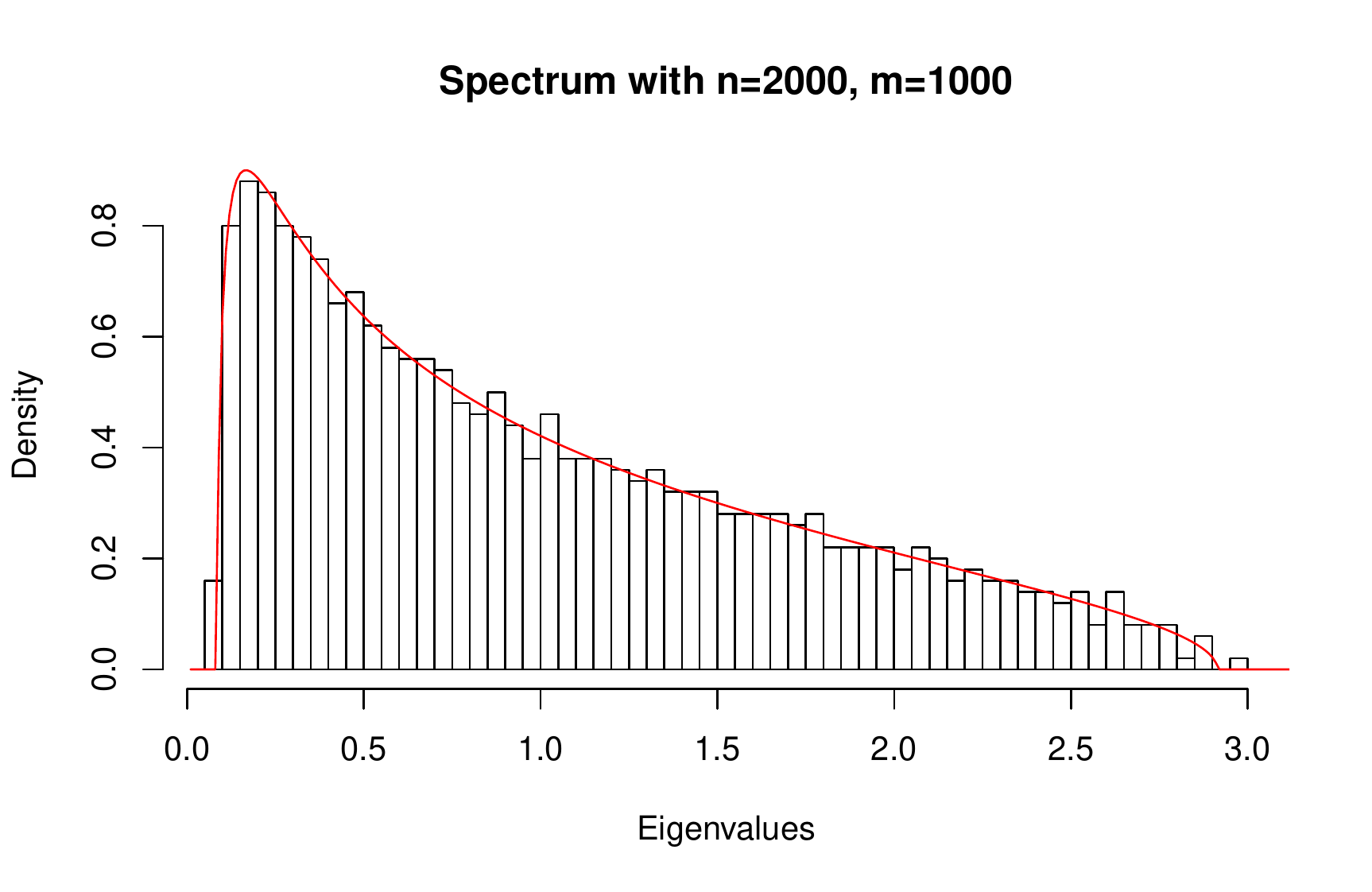} & 
 \includegraphics[width=0.26\paperwidth]{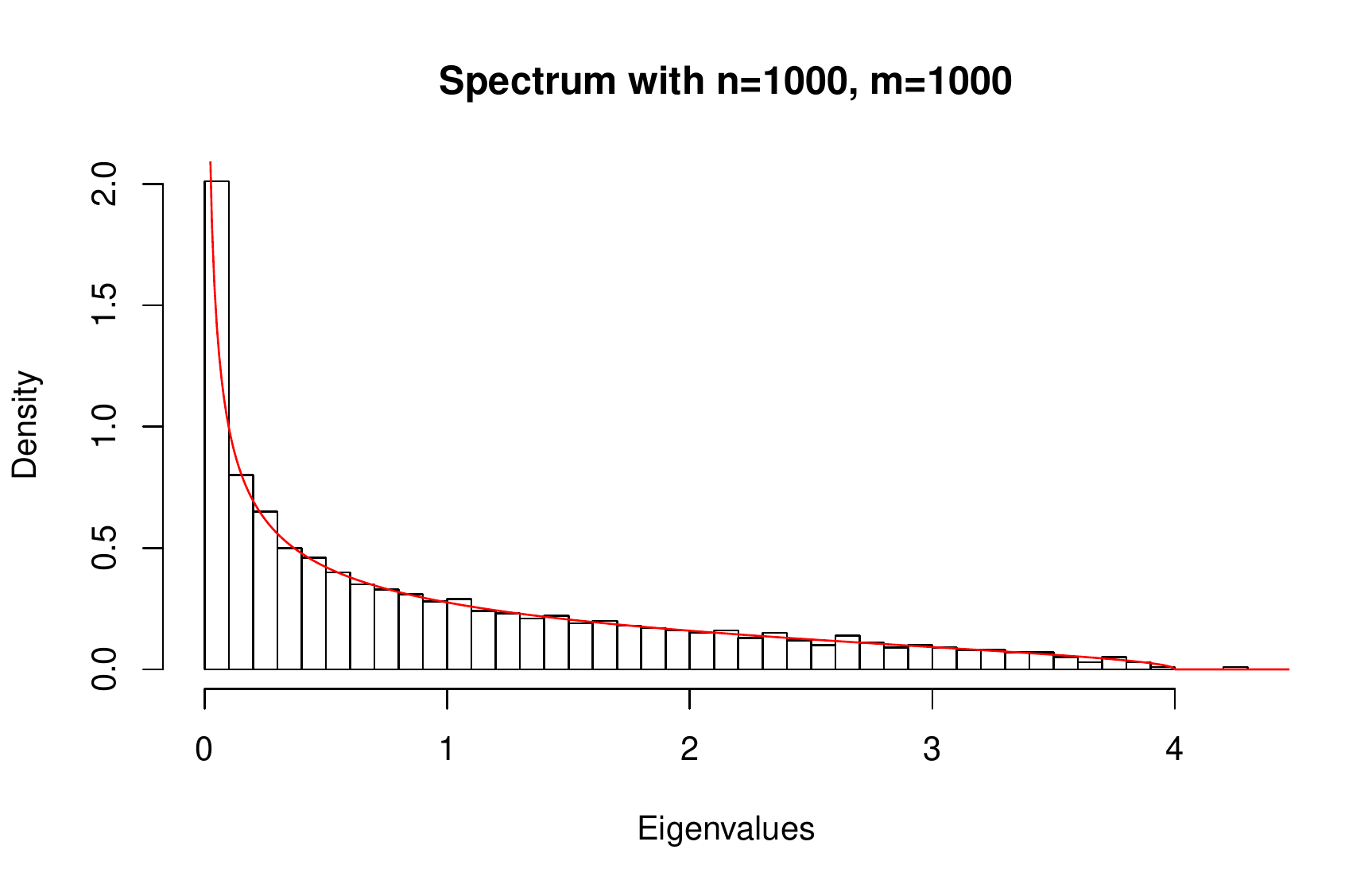} \\
 \includegraphics[width=0.26\paperwidth]{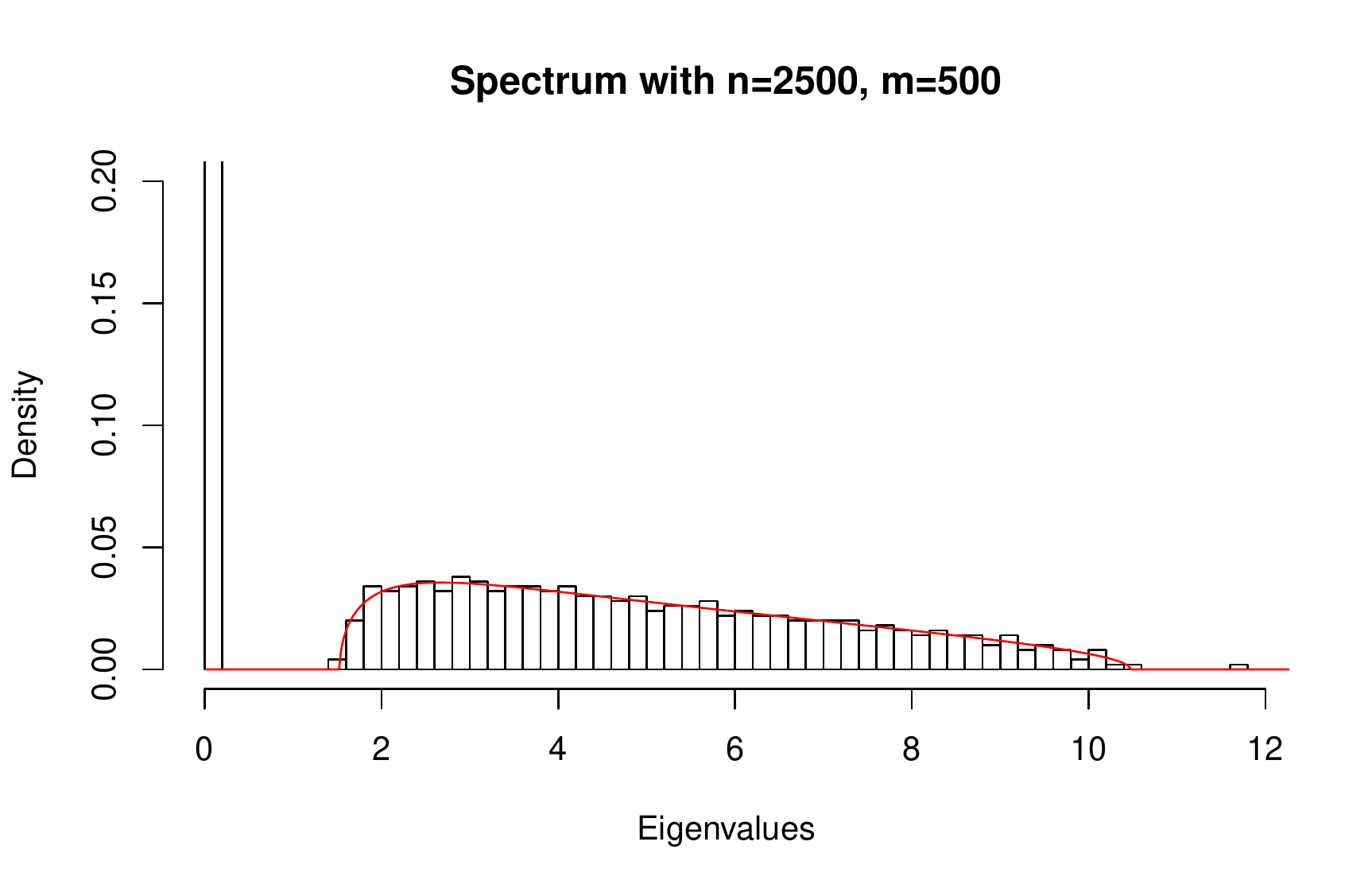}&
  \includegraphics[width=0.26\paperwidth]{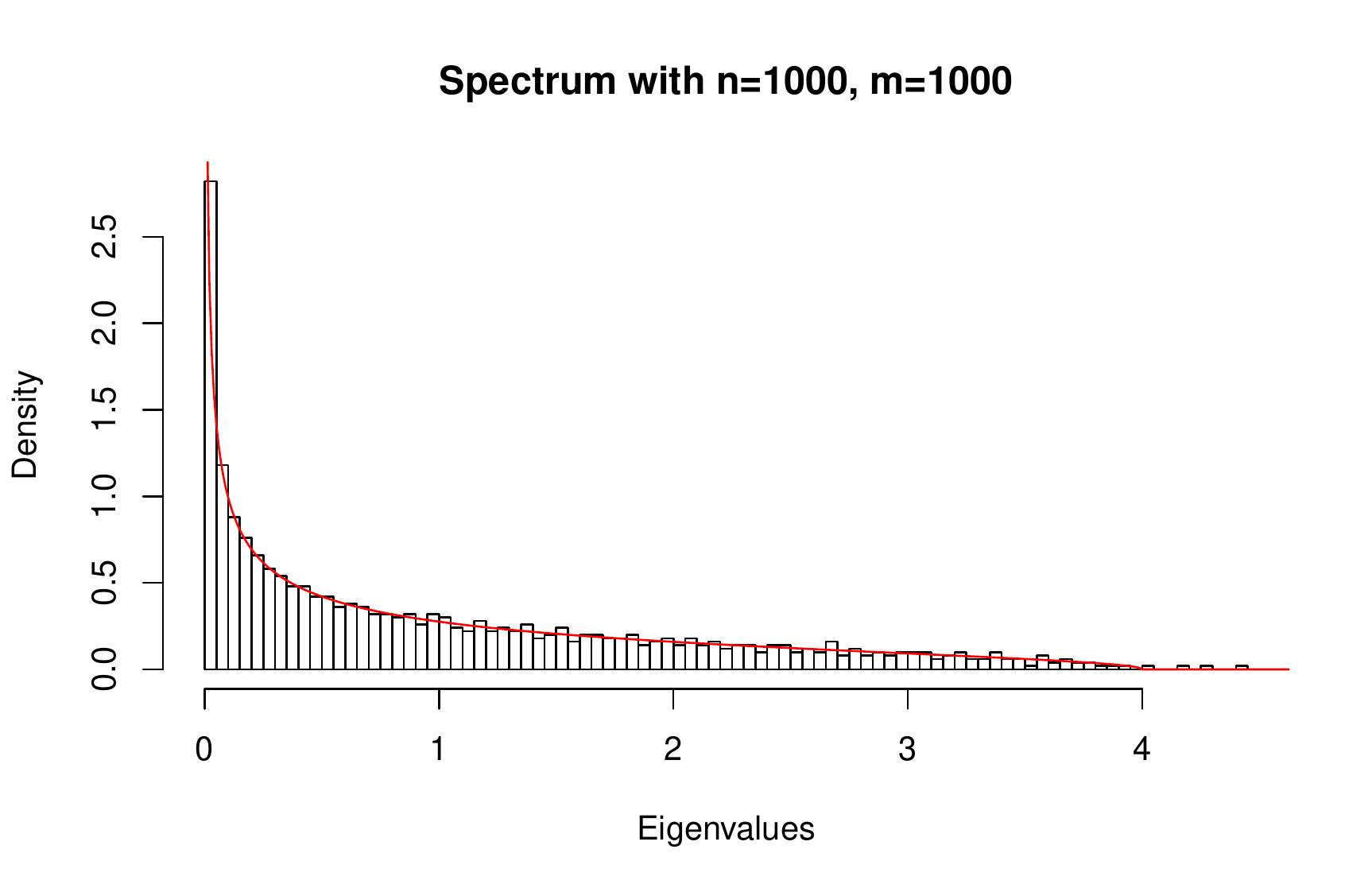}   
 \end{tabular}
\caption{Example of residual spikes of $\hat{\hat{\Sigma}}_X^{-1/2} \hat{\Sigma}_Y \hat{\hat{\Sigma}}_X^{-1/2}$ when $\theta=10$ for the first three figures and $\theta_{1,2,3,4}=10,15,20,25$ for the last figure.}\label{fig=spike}
\end{figure}

In the particular case where $\X$ and $\Y$ have independent jointly normal columns vector with constant variance $P=P_X=P_Y$, the distribution of the spectrum of the second of the above matrices is approximately Marcenko-Pastur distributed (see \cite{deformed}). This follows because $\hat{\hat{\Sigma}}_X$ is a finite perturbation. However, because of the non-consistency of the eigenvectors presented in \cite{deformedRMT}, we may observe \textbf{residual spikes} in the spectra, as shown in Figure \ref{fig=spike}. Thus, even if the two random matrices are based on the same perturbation, we see some spikes outside the bulk. This observation is worse in the last plot because four spikes fall outside the bulk even if there is actually no difference! 
This poses a fundamental problem for our test, because we must be able to distinguish the spikes indicative of a true difference from the residual spikes.  These remarks lead to the following definition.
\begin{Def}
The \textbf{residual spikes} are the isolated eigenvalues of 
$$\hat{\hat{\Sigma}}_X^{-1/2} \hat{\hat{\Sigma}}_Y \hat{\hat{\Sigma}}_X^{-1/2}\text{ or of }
\hat{\hat{\Sigma}}_X^{-1/2} \hat{\Sigma}_Y \hat{\hat{\Sigma}}_X^{-1/2}$$
when $P_X=P_Y$ (under the null hypothesis).
The \textbf{residual zone}  is the interval where a residual spike can fall asymptotically. 
\end{Def}

\begin{figure}[htbp] 
\centering
\scalebox{0.7}{
\begin{tabular}{c}
\definecolor{wwqqzz}{rgb}{0.4,0.,0.6}
\definecolor{uuuuuu}{rgb}{0.26666666666666666,0.26666666666666666,0.26666666666666666}
\definecolor{ttqqqq}{rgb}{0.2,0.,0.}
\definecolor{ffqqqq}{rgb}{1.,0.,0.}
\begin{tikzpicture}[line cap=round,line join=round,>=triangle 45,x=3.0cm,y=4.0cm]
\draw[->,color=black] (-0.1260631217517215,0.) -- (5.068082604362673,0.);
\foreach \x in {,0.2,0.4,0.6,0.8,1.,1.2,1.4,1.6,1.8,2.,2.2,2.4,2.6,2.8,3.,3.2,3.4,3.6,3.8,4.,4.2,4.4,4.6,4.8,5.}
\draw[shift={(\x,0)},color=black] (0pt,-2pt);
\draw[->,color=black] (0.,-0.3711350758839508) -- (0.,0.9716909818012385);
\foreach \y in {-0.3,-0.2,-0.1,0.2,0.3,0.4,0.5,0.6,0.7,0.8,0.9}
\draw[shift={(0,\y)},color=black] (0pt,-2pt) -- (0pt,-2pt);
\clip(-0.1260631217517215,-0.3711350758839508) rectangle (5.068082604362673,0.9716909818012385);
\draw[line width=1.2pt,color=ffqqqq,smooth,samples=100,domain=0.0857904999793061:2.914212816177216] plot(\x,{sqrt(((1.0+sqrt(0.5))^(2.0)-(\x))*((\x)-(1.0-sqrt(0.5))^(2.0)))/2.0/3.1415926535/0.5/(\x)});
\draw [color=ffqqqq](0.42792951474866847,0.866057465940425) node[anchor=north west] {Marcenko-Pastur};
\draw (1.484249715766261,-0.04943300485329163) node[anchor=north west] {$1+\sqrt{c}$};
\draw (3,0.5) node[anchor=north west] {$\hat{\hat{\Sigma}}_X^{-1/2} \hat{\Sigma}_Y \hat{\hat{\Sigma}}_X^{-1/2} $};
\draw (2.644287988169762,-0.015822340715760077) node[anchor=north west] {$\left(1+\sqrt{c}\right)^2$};
\draw (0.8159438668999814,-0.06063655956580215) node[anchor=north west] {$1$};
\draw (3.5176468364387627,-0.017422848531833006) node[anchor=north west] {$T^{-1}\left(\frac{1}{\lambda-1}\right)$};
\draw (3.52508297921058,-0.16306905979446976) node[anchor=north west] {$\lambda=\frac{1}{2}\left( 2+c+\sqrt{c^2+4c} \right)$};
\draw [line width=5.2pt,color=wwqqzz] (2.914213562373095,0.)-- (3.722220771172429,0.0026511965675838013);
\draw [color=wwqqzz](3.02686141349882,0.19224367537372097) node[anchor=north west] {Residual zone};
\begin{scriptsize}
\draw [color=ttqqqq] (2.914213562373095,0.)-- ++(-4.0pt,0 pt) -- ++(8.0pt,0 pt) ++(-4.0pt,-4.0pt) -- ++(0 pt,8.0pt);
\draw [color=black] (1.69813,0.)-- ++(-4.0pt,0 pt) -- ++(8.0pt,0 pt) ++(-4.0pt,-4.0pt) -- ++(0 pt,8.0pt);
\draw [color=uuuuuu] (0.89908,0.)-- ++(-4.0pt,0 pt) -- ++(8.0pt,0 pt) ++(-4.0pt,-4.0pt) -- ++(0 pt,8.0pt);
\draw [color=black] (3.722220771172429,0.0026511965675838013)-- ++(-4.0pt,0 pt) -- ++(8.0pt,0 pt) ++(-4.0pt,-4.0pt) -- ++(0 pt,8.0pt);
\end{scriptsize}
\end{tikzpicture} \\ \  \definecolor{wwqqzz}{rgb}{0.4,0.,0.6}
\definecolor{uuuuuu}{rgb}{0.26666666666666666,0.26666666666666666,0.26666666666666666}
\begin{tikzpicture}[line cap=round,line join=round,>=triangle 45,x=2.98cm,y=8.0cm]
\draw[->,color=black] (-0.1260631217517215,0.) -- (5.068082604362673,0.);
\foreach \x in {,0.2,0.4,0.6,0.8,1.,1.2,1.4,1.6,1.8,2.,2.2,2.4,2.6,2.8,3.,3.2,3.4,3.6,3.8,4.,4.2,4.4,4.6,4.8,5.}
\draw[shift={(\x,0)},color=black] (0pt,-2pt);
\draw[->,color=black] (0.,-0.1710715988748344) -- (0.,0.16690661773675);
\foreach \y in {}
\draw[shift={(0,\y)},color=black] (0pt,-2pt) -- (0pt,-2pt);
\clip(-0.1260631217517215,-0.1710715988748344) rectangle (5.068082604362673,0.16690661773675);
\draw (1.584249715766261,-0.04943300485329163) node[anchor=north west] {$1+\sqrt{c}$};
\draw (0.8629438668999814,-0.06063655956580215) node[anchor=north west] {$1$};
\draw [color=wwqqzz](2.052726710390752,0.13302488617902253) node[anchor=north west] {Residual zone};
\draw (2.5881289899615982,-0.04623198922114577) node[anchor=north west] {$\lambda=1+c+\sqrt{c^2+2c} $};
\draw (3,0.15) node[anchor=north west] {$\hat{\hat{\Sigma}}_X^{-1/2} \hat{\hat{\Sigma}}_Y \hat{\hat{\Sigma}}_X^{-1/2} $};
\draw [line width=5.2pt,color=wwqqzz] (1.69813,0.)-- (2.9859626282538247,1.827374449692361E-4);
\begin{scriptsize}
\draw [color=black] (1.69813,0.)-- ++(-4.0pt,0 pt) -- ++(8.0pt,0 pt) ++(-4.0pt,-4.0pt) -- ++(0 pt,8.0pt);
\draw [color=uuuuuu] (0.89908,0.)-- ++(-4.0pt,0 pt) -- ++(8.0pt,0 pt) ++(-4.0pt,-4.0pt) -- ++(0 pt,8.0pt);
\draw [color=black] (2.9859626282538247,1.827374449692361E-4)-- ++(-4.0pt,0 pt) -- ++(8.0pt,0 pt) ++(-4.0pt,-4.0pt) -- ++(0 pt,8.0pt);
\end{scriptsize}
\end{tikzpicture}  
\end{tabular} 
}
\caption{Residual zone of $\hat{\hat{\Sigma}}_X^{-1/2} \hat{\Sigma}_Y \hat{\hat{\Sigma}}_X^{-1/2}$ and $\hat{\hat{\Sigma}}_X^{-1/2} \hat{\hat{\Sigma}}_Y \hat{\hat{\Sigma}}_X^{-1/2}$.}\label{fig=residualzone0}
\end{figure}
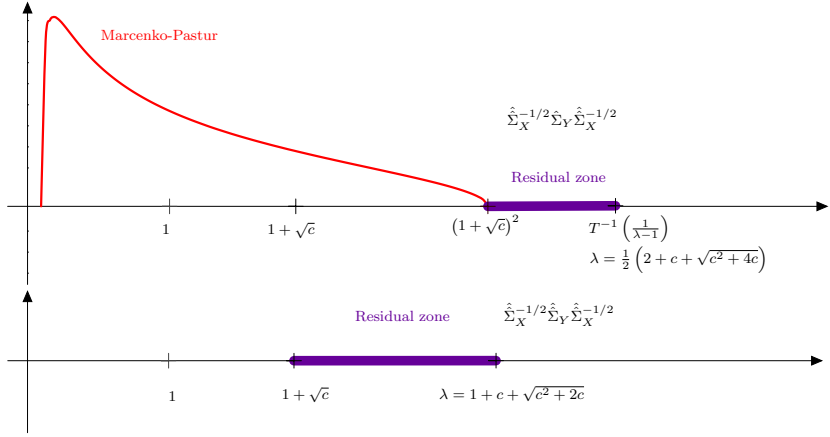

This paper studies these residual spikes by deriving the distribution of the extreme residual spikes under the null hypothesis. The philosophy is explained in Figure \ref{fig=residualzone0} with illustrations inspired by the i.i.d. normal case. All the eigenvalues lying in what we call the residual zone are potentially not indicative of real differences. However, when an eigenvalue is larger, we declare that this spike expresses a true difference. 

Most of our plots feature the seemingly more natural matrix
\begin{eqnarray*}
\hat{\hat{\Sigma}}_X^{-1/2} \hat{\Sigma}_Y \hat{\hat{\Sigma}}_X^{-1/2}\,.
\end{eqnarray*}
But, although this choice simplifies the study in terms of convergence in probability when the perturbation is of order $1$, this is no longer the case in more complex situations. In addition, the eigenvectors associated with the residual spikes are more accessible for the matrix in which all estimates are filtered.

Let $\hat\theta_X$ and $\hat\theta_Y$ be isolated eigenvalues and construct the asymptotic unbiased estimators as in Equation (\ref{eq:corrlambda})
\begin{eqnarray*}
\hat{\hat{\theta}}_{X}=1+\frac{1}{\frac{1}{m} \sum_{i=1}^m \frac{\hat{\lambda}_{W_X,i}}{\hat{\theta}_{X}-\hat{\lambda}_{W_X,i}}} \text{ and } \hat{\hat{\theta}}_{Y}=1+\frac{1}{\frac{1}{m} \sum_{i=1}^m \frac{\hat{\lambda}_{W_Y,i}}{\hat{\theta}_{Y}-\hat{\lambda}_{W_Y,i}}}\,.
\end{eqnarray*}
Here $\hat{\lambda}_{W_X,i}$ and $\hat{\lambda}_{W_Y,i}$ are the eigenvalues of $W_X$ and $W_Y$, respectively.
In practice we do of course not observe $W_X$ and $W_Y$, but a simple argument using Cauchy's interlacing law shows that we can replace the previous estimators by 
\begin{eqnarray*}
\hat{\hat{\theta}}_{X}=1+\frac{1}{\frac{1}{m-1} \sum_{i=2}^{m} \frac{\hat{\lambda}_{\hat{\Sigma}_{X},i}}{\hat{\theta}_{X}-\hat{\lambda}_{\hat{\Sigma}_{X},i}}} \text{ and } \hat{\hat{\theta}}_{Y}=1+\frac{1}{\frac{1}{m-1} \sum_{i=2}^m \frac{\hat{\lambda}_{\hat{\Sigma}_{Y},i}}{\hat{\theta}_{Y}-\hat{\lambda}_{\hat{\Sigma}_{Y},i}}},
\end{eqnarray*}
where $\hat{\lambda}_{\hat{\Sigma}_{X},i}$ and $\hat{\lambda}_{\hat{\Sigma}_{Y},i}$ are the i$^{th}$ ordered eigenvalue of $\hat{\Sigma}_{X}$ and $\hat{\Sigma}_{Y}$, respectively.
The test statistic is then 
\begin{eqnarray*}
 \lambda_{\min}\left(\hat{\hat{\Sigma}}_X^{-1/2} \hat{\hat{\Sigma}}_Y \hat{\hat{\Sigma}}_X^{-1/2} \right) \text{ and } \lambda_{\max}\left(\hat{\hat{\Sigma}}_X^{-1/2} \hat{\hat{\Sigma}}_Y \hat{\hat{\Sigma}}_X^{-1/2} \right)\,,
\end{eqnarray*}
where the filtered matrices are constructed as in (\ref{eq:corrSigma}).
These two statistics provide a basis for a powerful and robust test for the equality of (detectable) perturbations $P_X$ and $P_Y$.

\subsection{Null distribution} \label{sec=test}
Obviously under $H_0$, $\lambda_{\max}\left(\hat{\hat{\Sigma}}_X^{-1/2} \hat{\hat{\Sigma}}_Y \hat{\hat{\Sigma}}_X^{-1/2} \right) $ is a function of $\theta=\theta_X=\theta_Y$. The suspected \textbf{worst case} occurs in the limit as $\theta\to\infty$ and it is this limit which will determine the critical values of the test. This can be checked by Criterion \ref{CriteriaTsmall} which we discuss later. Let 
\begin{eqnarray*}
&&\lambda_{\max}\left(\hat{\hat{\Sigma}}_X^{-1/2} \hat{\hat{\Sigma}}_Y \hat{\hat{\Sigma}}_X^{-1/2} \right) \leqslant \lim_{\theta \rightarrow \infty} \lambda_{\max}\left(\hat{\hat{\Sigma}}_X^{-1/2} \hat{\hat{\Sigma}}_Y \hat{\hat{\Sigma}}_X^{-1/2} \right) =V_{\max} ,\\
&&\lambda_{\min}\left(\hat{\hat{\Sigma}}_X^{-1/2} \hat{\hat{\Sigma}}_Y \hat{\hat{\Sigma}}_X^{-1/2} \right) \geqslant \lim_{\theta \rightarrow \infty} \lambda_{\min}\left(\hat{\hat{\Sigma}}_X^{-1/2} \hat{\hat{\Sigma}}_Y \hat{\hat{\Sigma}}_X^{-1/2} \right) =V_{\min}.
\end{eqnarray*}
Because of our focus on the worst case scenario under $H_0$, we will investigate the asymptotics as $\frac{\theta}{\sqrt{m}} \rightarrow \infty$.
Our test rejects the null hypothesis of equal populations if 
either $\P\left( V_{\max} > \hat{\lambda}_{\max}   \right)$ or $\P\left( V_{\min} < \hat{\lambda}_{\min}   \right)$ is small, where $\hat{\lambda}_{\max}$ and $\hat{\lambda}_{\min}$ are the observed extreme residual spikes.

In the investigation of the extremal eigenvalues under multiplicative perturbations, we make use of the following random variables
\begin{eqnarray*}
M_{s_1,s_2,X}(\rho_X)&=&\frac{1}{m} \sum_{i=1}^m \frac{\hat{\lambda}_{W_X,i}^{s_1}}{\left( \rho_X -\hat{\lambda}_{W_X,i} \right)^2},\\
M_{s_1,s_2,Y}(\rho_Y)&=&\frac{1}{m} \sum_{i=1}^m \frac{\hat{\lambda}_{W_Y,i}^{s_1}}{\left( \rho_Y -\hat{\lambda}_{W_Y,i} \right)^2},\\
M_{s_1,s_2}(\rho_X,\rho_Y)&=& \frac{M_{s_1,s_2,X}(\rho_X)+M_{s_1,s_2,Y}(\rho_Y)}{2}.
\end{eqnarray*}
In particular, when $s_2=0$, we use $M_{s_1,X}=M_{s_1,0,X}$. When we only study one group, we use the simpler notation $M_{s_1,s_2}(\rho)$ when no confusion is possible. Note that $M_{1,1,X}(\rho_X)=\hat{T}_{W_X}(\rho_X)$ is the empirical T-transform. Moreover, for an applied perspective, all these variables can be estimated by
$$M_{s_1,s_2,X}(\rho_X)=\hat{M}_{s_1,s_2,X}(\rho_X)=\frac{1}{m-1} \sum_{i=2}^m \frac{\hat{\lambda}_{\hat{\Sigma}_X,i}^{s_1}}{\left( \rho_X -\hat{\lambda}_{\hat{\Sigma}_X,i} \right)^2}.$$

The following result describes the asymptotic behavior of the extreme eigenvalues and thus of $V_{\max}$ and $V_{\min}$. .
\begin{Th}
 \label{TH=Main} 
Suppose $W_X,W_Y \in \mathbb{R}^{m\times m}$ are as described at the start of Section 2 and 
$P= \I_m+(\theta-1)e_1 e_1^t \in  \mathbb{R}^{m\times m}$ with $\frac{\sqrt{m}}{\theta}=o(1)$ with regard to large $m$. Let  
\begin{eqnarray*}
\hat{\Sigma}_{X}=P^{1/2} W_X P^{1/2} \text{ and } \hat{\Sigma}_{Y}=P^{1/2} W_Y P^{1/2}.
\end{eqnarray*}
and $\hat{\hat{\Sigma}}_{X}$, $\hat{\hat{\Sigma}}_{Y}$ as described above (see, \ref{Def=unbiased}). 

Then, conditional on the spectra $S_{W_X}=\left\lbrace \hat{\lambda}_{W_X,1},\hat{\lambda}_{W_X,2},...,\hat{\lambda}_{W_X,m} \right\rbrace$ and $S_{W_Y}=\left\lbrace \hat{\lambda}_{W_Y,1},\hat{\lambda}_{W_Y,2},...,\hat{\lambda}_{W_Y,m} \right\rbrace$ of $W_X$ and $W_Y$,
\begin{eqnarray*}
\left. \sqrt{m} \frac{\left(\lambda_{\max}\left(\hat{\hat{\Sigma}}_{X}^{-1/2} \hat{\hat{\Sigma}}_{Y} \hat{\hat{\Sigma}}_{X}^{-1/2} \right) -\lambda^+ \right)}{\sigma^+} \right| S_{W_X}, S_{W_Y} \sim \Normal(0,1)+o_{p}(1), 
\end{eqnarray*}
where
\begin{eqnarray*}
&&\lambda^+= \sqrt{M_2^2-1}+M_2, \hspace{30cm}
\end{eqnarray*}
\scalebox{0.77}{
\begin{minipage}{1\textwidth}
\begin{eqnarray*}
&&{\sigma^+}^2=\frac{1}{\left(M_{2,X}+M_{2,Y}-2\right) \left(M_{2,X}+M_{2,Y}+2\right)} \hspace{30cm} \\
&&\hspace{2cm} \Bigg( 9 M_{2,X}^4 M_{2,Y}+4 M_{2,X}^3 M_{2,Y}^2+4 M_{2,X}^3 M_{2,Y}+2 M_{2,X}^3 M_{3,Y}-2 M_{2,X}^2 M_{2,Y}^3\\
&& \hspace{2cm}+4 M_{2,X}^2 M_{2,Y}^2-11 M_{2,X}^2 M_{2,Y}-8 M_{3,X} M_{2,X}^2 M_{2,Y}+2 M_{2,X}^2 M_{2,Y} M_{3,Y}\\
&& \hspace{2cm}-2 M_{2,X}^2 M_{3,Y}+M_{2,X}^2 M_{4,Y}+4 M_{2,X} M_{2,Y}^3+M_{2,X} M_{2,Y}^2+4 M_{2,X} M_{2,Y}\\
&& \hspace{2cm}-4 M_{3,X} M_{2,X} M_{2,Y}^2-4 M_{3,X} M_{2,X} M_{2,Y}-2 M_{2,X} M_{2,Y}^2 M_{3,Y}-4 M_{2,X} M_{2,Y} M_{3,Y}\\
&& \hspace{2cm}-6 M_{2,X} M_{3,Y}+2 M_{4,X} M_{2,X} M_{2,Y}+2 M_{2,X} M_{2,Y} M_{4,Y}-2 M_{3,X} M_{2,Y}^2\\
&& \hspace{2cm}+2 M_{3,X} M_{2,Y}+M_{4,X} M_{2,Y}^2+4 M_{2,X}^5+2 M_{2,X}^4-4 M_{3,X} M_{2,X}^3-13 M_{2,X}^3\\
&& \hspace{2cm}-2 M_{3,X} M_{2,X}^2+M_{4,X} M_{2,X}^2-2 M_{2,X}^2+10 M_{3,X} M_{2,X}+4 M_{2,X}+4 M_{3,X}\\
&& \hspace{2cm}-2 M_{4,X}+M_{2,Y}^5+2 M_{2,Y}^4-M_{2,Y}^3-2 M_{2,Y}^2+4 M_{2,Y}-2 M_{2,Y}^3 M_{3,Y}\\
&& \hspace{2cm}-2 M_{2,Y}^2 M_{3,Y}+2 M_{2,Y} M_{3,Y}+4 M_{3,Y}+M_{2,Y}^2 M_{4,Y}-2 M_{4,Y}-4 \Bigg)\\
&&\hspace{1.25cm} + \frac{1}{\sqrt{\left(M_{2,X}+M_{2,Y}-2\right) \left(M_{2,X}+M_{2,Y}+2\right)}}\\
&& \hspace{2cm} \Bigg( 5 M_{2,X}^3 M_{2,Y}-M_{2,X}^2 M_{2,Y}^2+2 M_{2,X}^2 M_{2,Y}+2 M_{2,X}^2 M_{3,Y}-M_{2,X} M_{2,Y}^3\\
&& \hspace{2cm}+2 M_{2,X} M_{2,Y}^2-4 M_{2,X} M_{2,Y}-4 M_{3,X} M_{2,X} M_{2,Y}-2 M_{2,X} M_{3,Y}+M_{2,X} M_{4,Y}\\
&& \hspace{2cm}-2 M_{3,X} M_{2,Y}+M_{4,X} M_{2,Y}+4 M_{2,X}^4+2 M_{2,X}^3-4 M_{3,X} M_{2,X}^2-5 M_{2,X}^2\\
&& \hspace{2cm}-2 M_{3,X} M_{2,X}+M_{4,X} M_{2,X}+2 M_{2,X}+2 M_{3,X}+M_{2,Y}^4+2 M_{2,Y}^3+M_{2,Y}^2\\
&& \hspace{2cm}+2 M_{2,Y}-2 M_{2,Y}^2 M_{3,Y}-2 M_{2,Y} M_{3,Y}-2 M_{3,Y}+M_{2,Y} M_{4,Y} \Bigg),
\end{eqnarray*}
\end{minipage}}
\begin{eqnarray*}
&&M_{s,X}= \frac{1}{m} \sum_{i=1}^m \hat{\lambda}_{W_X,i}^s,\hspace{30cm}\\
&&M_{s,Y}= \frac{1}{m} \sum_{i=1}^m \hat{\lambda}_{W_Y,i}^s,\\
&&M_s=\frac{M_{s,X}+M_{s,Y}}{2}.
\end{eqnarray*}
Moreover, 
\begin{eqnarray*}
\left.\sqrt{m} \frac{\left(\lambda_{\min}\left(\hat{\hat{\Sigma}}_{X}^{-1/2} \hat{\hat{\Sigma}}_{Y} \hat{\hat{\Sigma}}_{X}^{-1/2} \right) -\lambda^- \right)}{\sigma^-}\right| S_{W_X}, S_{W_Y} \sim \Normal(0,1)+o_{m}(1), 
\end{eqnarray*}
where 
\begin{eqnarray*}
&&\lambda^-= -\sqrt{M_2^2-1}+M_2, \hspace{30cm}\\
&&{\sigma^-}^2=\left(\lambda^-\right)^4 {\sigma^+}^2.
\end{eqnarray*}
The error $o_p(1)$ in the approximation is with regard to large values of $m$.
\end{Th}
\paragraph*{Special case}
If the spectra are Marcenko-Pastur distributed, then:

\scalebox{0.75}{
\begin{minipage}{1\textwidth}
\begin{eqnarray*}
&&c=\frac{c_X+c_Y}{2}, \text{ where } c_X=\frac{n_X}{m} \text{ and } c_Y=\frac{n_Y}{m},\hspace{30cm}\\
&&\lambda^+=  c+\sqrt{c (c+2)}+1,\\
&&{\sigma^+}^2= c_X^3+c_X^2 c_Y+3 c_X^2+4 c_X c_Y-c_X+c_Y^2+c_Y\\
&& \hspace{2cm}+\frac{ (8 c_X+2 c_X^2+\left(c_X^3+5 c_X^2+c_X^2 c_Y+4 c_X c_Y+5 c_X +3 c_Y+c_Y^2\right)\sqrt{c(c+2)}}{c+2}.
\end{eqnarray*}
\end{minipage}}\\
If $c_{X}$ tends to $0$, then \\
\scalebox{0.75}{
\begin{minipage}{1\textwidth}
\begin{eqnarray*}
{\sigma^+}^2&=&  \Bigg( M_{2,Y}^5+2 M_{2,Y}^4-2 M_{3,Y} M_{2,Y}^3+M_{2,Y}^3-4 M_{3,Y} M_{2,Y}^2+M_{4,Y} M_{2,Y}^2 +2 M_{2,Y}^2 \hspace{30cm}\\
&&\hspace{1cm} +2 M_{4,Y} M_{2,Y} +2 M_{2,Y}-2 M_{3,Y}-M_{4,Y}-2 \Bigg) \Bigg/ \Bigg( \left(M_{2,Y}-1\right) \left(M_{2,Y}+3\right)\Bigg)\\
&&+  \Bigg(M_{2,Y}^4+M_{2,Y}^3-2 M_{3,Y} M_{2,Y}^2+2 M_{2,Y}^2-2 M_{3,Y} M_{2,Y} +M_{4,Y} M_{2,Y}\\
&& \hspace{1cm} -2 M_{3,Y}+M_{4,Y} \Bigg)\Bigg/ \sqrt{\left(M_{2,Y}-1\right) \left(M_{2,Y}+3\right)}. 
\end{eqnarray*}
\end{minipage}}

(Proof in supplement material \cite{Suppmaterial}.)

\subsection{Discussion and simulation}\label{sec:smallsimulationmainth}

The above theorem gives the limiting distribution of $V_{\max}$ and $V_{\min}$. In this subsection, we first check the quality of the approximations in Theorem \ref{TH=Main}, then we investigate the worst case with regard to $\theta$.

\subsubsection{Some simulations}

Assume $\mathbf{X} \in \mathbb{R}^{m\times n_X}$ and $\mathbf{Y} \in \mathbb{R}^{m\times n_Y}$ with $\mathbf{X}=\left(X_1,X_2,...,X_{n_X}\right)$ and $\mathbf{Y}=\left(Y_1,Y_2,...,Y_{n_Y}\right)$. The components of the random vectors are independent and the covariance between the vectors is as follows:

\scalebox{0.65}{
\begin{minipage}{1\textwidth}
\begin{eqnarray*}
&&X_i \sim \Normal_m\left(\vec{0},\sigma^2 \I_m \right) \text{ with } X_1 = \epsilon_{X,1}\text{ and } 
X_{i+1}=\rho X_i+\sqrt{1-\rho^2} \ \epsilon_{X,i+1}, \text{ where }  \epsilon_{X,i} \overset{i.i.d}{\sim} \Normal_m\left(\vec{0},\sigma^2 \I_m\right),\\
&&Y_i \sim \Normal_m\left(\vec{0},\sigma^2 \I_m \right) \text{ with } Y_1 = \epsilon_{Y,1}\text{ and } Y_{i+1}=\rho Y_i+\sqrt{1-\rho^2} \ \epsilon_{Y,i+1}, \text{ where }  \epsilon_{Y,i} \overset{i.i.d}{\sim} \Normal_m\left(\vec{0},\sigma^2 \I_m\right)
\end{eqnarray*}
\end{minipage}}

Let $P_X= \I_m + (\theta_{X}-1) u_{X} u_{X}^t$ and $P_Y= \I_m+  (\theta_{Y}-1) u_{Y} u_{Y}^t$ be two perturbations in  $\mathbb{R}^{m\times m}$. Then,
\begin{eqnarray*}
\mathbf{X}_P=P_X^{1/2} \mathbf{X} \text{ and  } \mathbf{Y}_P=P_Y^{1/2} \mathbf{Y},\\
\hat{\Sigma}_X=\frac{\mathbf{X}_P^t \mathbf{X}_P}{n_X} \text{ and } \hat{\Sigma}_Y=\frac{\mathbf{Y}_P^t \mathbf{Y}_P}{n_Y}.
\end{eqnarray*}
We assume a common and large value for $\theta$ and $P_X=P_Y$.  
\paragraph*{Distribution}  \ \\
Figure \ref{tab:residual} presents empirical distributions of the extreme residual spikes in different scenarios together with the normal densities from our theorem.

\begin{table}[hbt]
 \begin{tabular}{ c c }
 \begin{minipage}{.45\textwidth} \centering
      \includegraphics[width=1\textwidth]{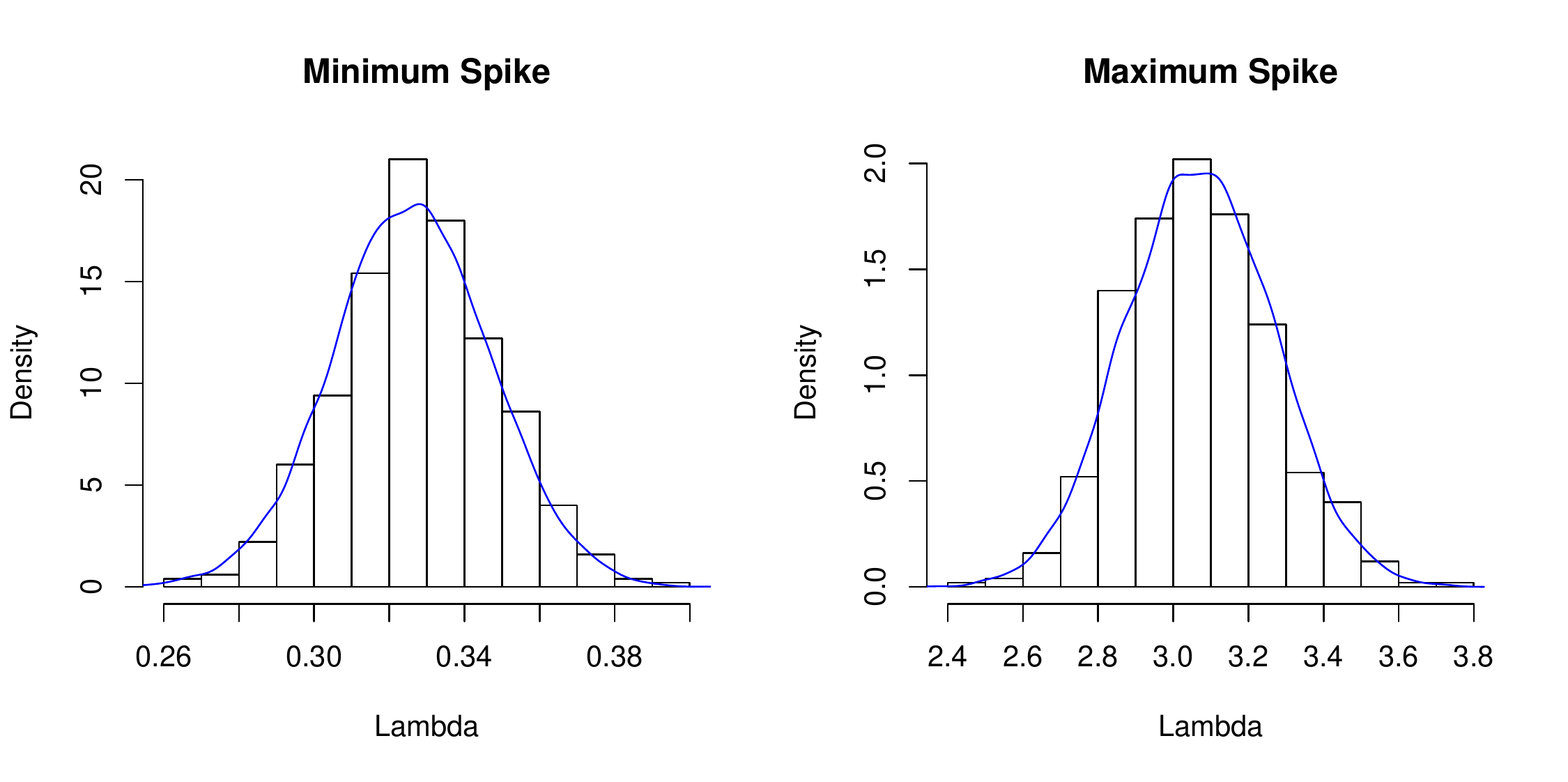}   
    \end{minipage} 
    &
    \begin{minipage}{.45\textwidth} \centering
      \includegraphics[width=1\textwidth]{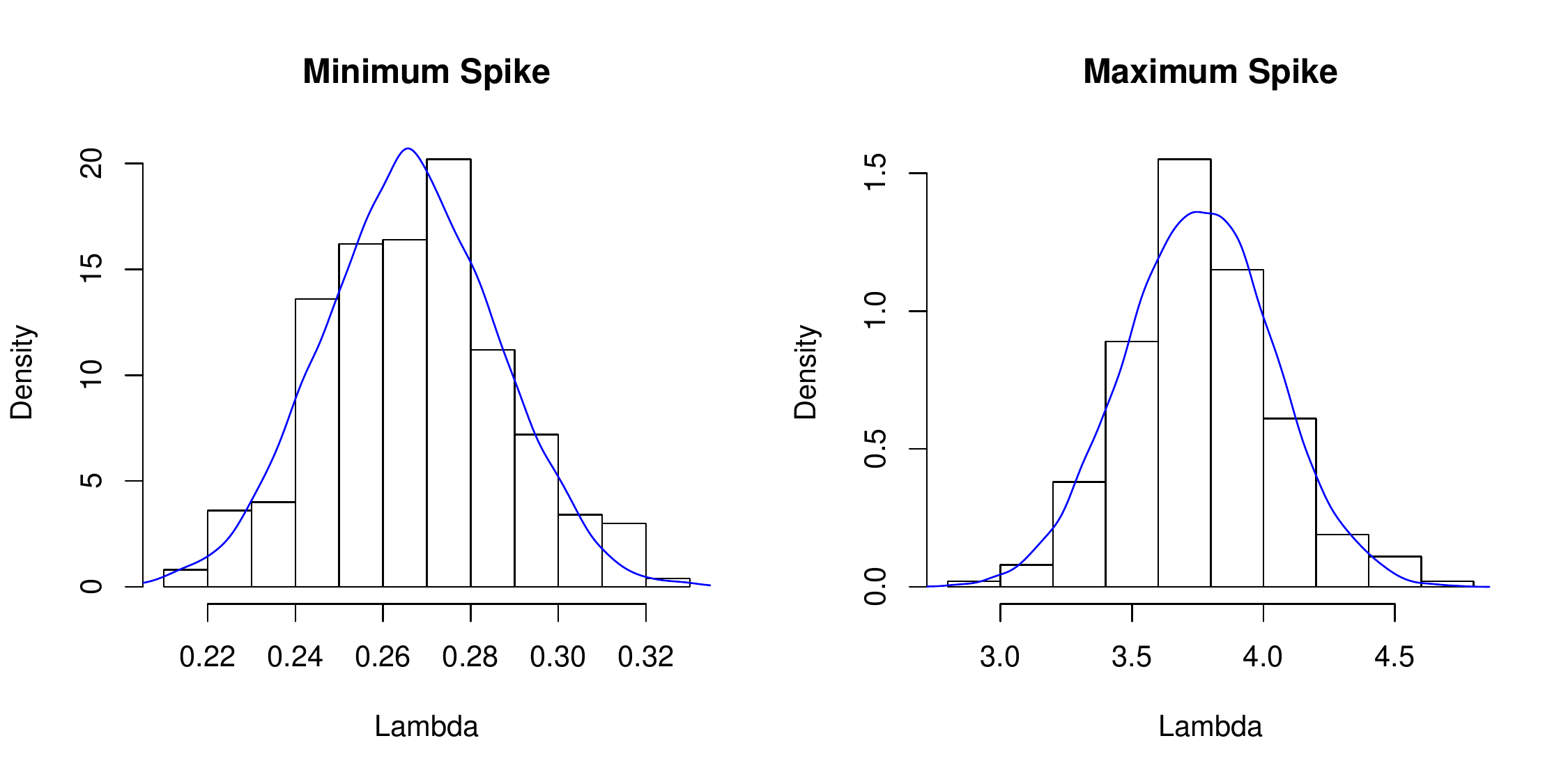}   
    \end{minipage}\\
\begin{minipage}{.45\textwidth} \centering
\scalebox{0.8}{
\begin{tabular}{c}
Scenario 1 \\
\begin{tabular}{c|c|c}
$\rho=0.5$ & $c_X=1/3$ & $c_Y=1/2$ \\ 
\hline 
$m=300$ & $n_X=900$ & $n_Y=600$ \\ 

\end{tabular} \\
 $\theta=5'000$.
\end{tabular}}
\end{minipage}
 & 
\begin{minipage}{.45\textwidth} \centering
\scalebox{0.8}{
\begin{tabular}{c}
Scenario 2 \\
\begin{tabular}{c|c|c}
$\rho=0.$ & $c_X=1$ & $c_Y=1$ \\ 
\hline 
$m=300$ & $n_X=300$ & $n_Y=300$ \\ 

\end{tabular} \\
 $\theta=5'000$.
\end{tabular}}
\end{minipage} \\
\begin{minipage}{.45\textwidth} \centering
      \includegraphics[width=1\textwidth]{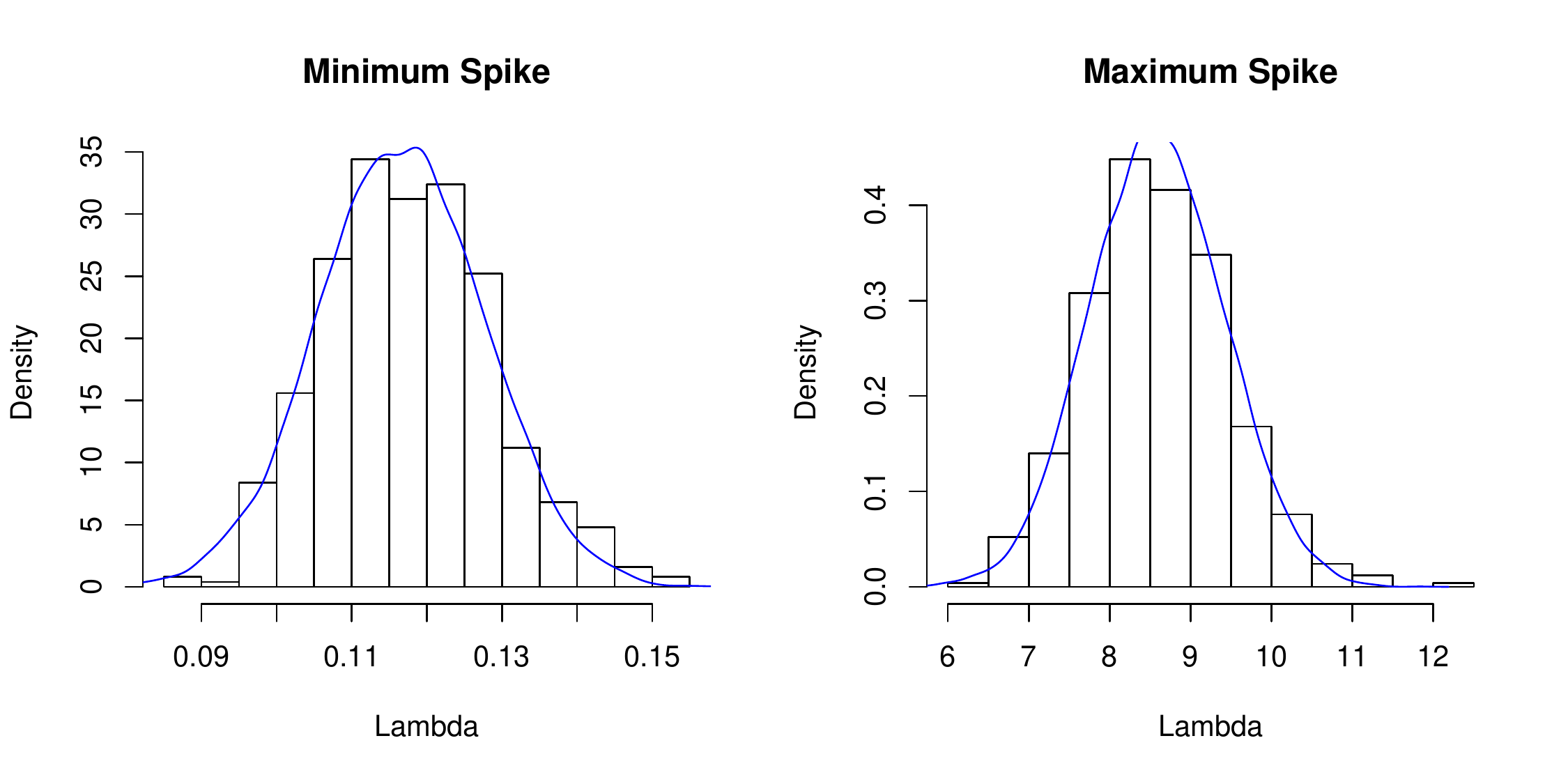}   
    \end{minipage} 
    &
    \begin{minipage}{.45\textwidth} \centering
      \includegraphics[width=1\textwidth]{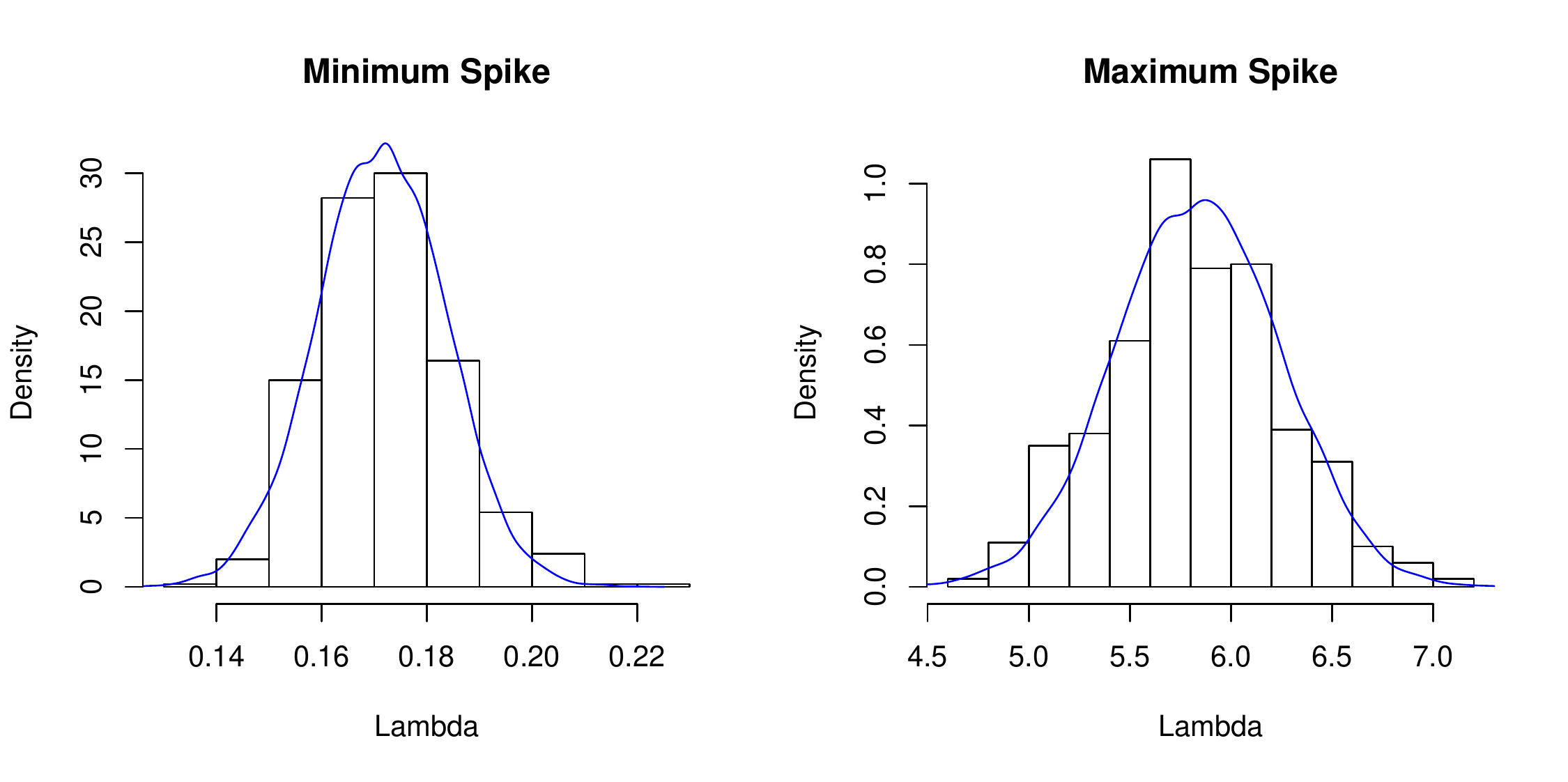}   
    \end{minipage}\\
\begin{minipage}{.45\textwidth} \centering
\scalebox{0.8}{
\begin{tabular}{c}
Scenario 3 \\
\begin{tabular}{c|c|c}
$\rho=0.5$ & $c_X=2$ & $c_Y=2$ \\ 
\hline 
$m=600$ & $n_X=300$ & $n_Y=300$ \\ 

\end{tabular} \\
 $\theta=5'000$.
\end{tabular}}
\end{minipage}
 & 
\begin{minipage}{.45\textwidth} \centering
\scalebox{0.8}{
\begin{tabular}{c}
Scenario 4 \\
\begin{tabular}{c|c|c}
$\rho=0.$ & $c_X=2$ & $c_Y=2$ \\ 
\hline 
$m=600$ & $n_X=300$ & $n_Y=300$ \\ 

\end{tabular} \\
 $\theta=5'000$.
\end{tabular}}
\end{minipage}
\end{tabular}
\caption{Empirical distributions of the residual spikes and the Gaussian densities from the theorem \ref{TH=Main} (in blue).  
} \label{tab:residual}
\end{table}

Appendix \ref{sec:Tablecomplete}, Table \ref{tab:Tablecomplete} contains a comparison of the estimates of the mean and the variance with empirical ones in diverse situations. It is noticeable that in situations where $m$ is large compared to $n_X$ or $n_Y$, the asymptotic results are less accurate.

\subsubsection{Increasing residual spike} \label{sec:Criter} 

\begin{figure}[hbt]
\center
\includegraphics[width=1\textwidth]{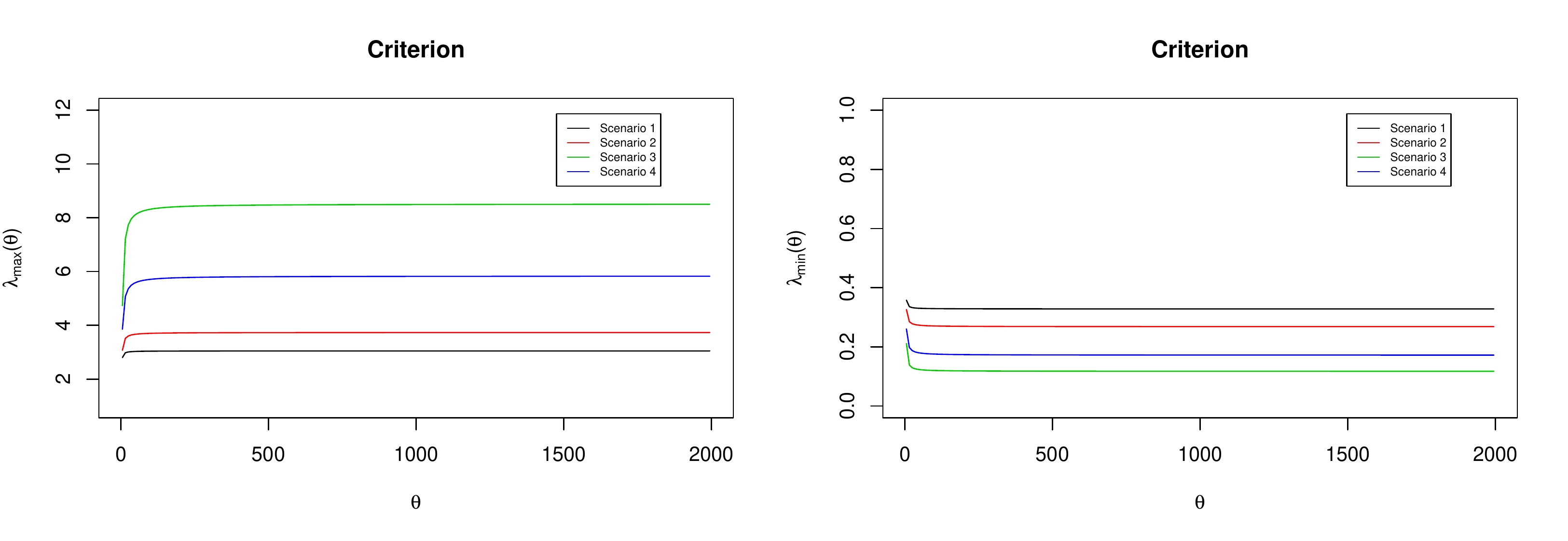} 
\caption{Plot of the criterion in the four scenarios used in the simulation.} \label{fig:increasingtheta}
\end{figure}

In the four scenarios used in the simulations, we can estimate the expectation of the residual spike. Figure \ref{fig:increasingtheta} shows that  the expectations of the largest residual spikes are always strictly increasing as a function of $\theta$ and the expectations of the smallest residual spikes are always strictly decreasing. This is, however, not universally true. To address this issue, the following criterion can be used.
\begin{Def}\label{CriteriaTsmall}
Suppose $\hat{\Sigma}_X$ and $\hat{\Sigma}_Y$ are two independent random estimated covariance matrices of the form described at the start of Section 2. Let
\begin{equation} \label{eq=CriteriaTsmall}
\resizebox{.85\hsize}{!}{$\mu_{\lambda}(\theta,S_{X},S_{Y})=\frac{1}{2} \left(\theta+\alpha^2-\theta \alpha^2+\frac{1+(\theta-1)\alpha^2 + \sqrt{-4\theta^2 + \left( 1+\theta^2 - (\theta-1)^2 \alpha^2 \right)^2}}{\theta} \right),$}
\end{equation} 
where $S_{X}=\left\lbrace \hat{\lambda}_{\hat{\Sigma}_X,2},\hat{\lambda}_{\hat{\Sigma}_X,3},...,\hat{\lambda}_{\hat{\Sigma}_X,m} \right\rbrace$ and $S_{Y}=\left\lbrace \hat{\lambda}_{\hat{\Sigma}_Y,2},\hat{\lambda}_{\hat{\Sigma}_Y,3},...,\hat{\lambda}_{\hat{\Sigma}_Y,m} \right\rbrace$ and 
\begin{equation*}
\resizebox{.95\hsize}{!}{$\begin{array}{lll}
\alpha=\alpha_{X} \alpha_{Y}, & \alpha_{X}^2=\frac{(m-1)\theta}{(\theta-1)^2 \hat{\theta}_{X}  \sum_{i=2}^m \frac{\hat{\lambda}_{\hat{\Sigma}_X,i}}{(\hat{\theta}_{X}-\hat{\lambda}_{\hat{\Sigma}_X,i})^2}}, & \alpha_{Y}^2=\frac{(m-1)\theta}{(\theta-1)^2 \hat{\theta}_{Y}  \sum_{i=2}^m \frac{\hat{\lambda}_{\hat{\Sigma}_Y,i}}{(\hat{\theta}_{Y}-\hat{\lambda}_{\hat{\Sigma}_Y,i})^2}},\\
& \hat{\theta}_{X} \ \bigg| \  \frac{1}{\theta-1} = \frac{1}{m-1} \sum_{i=2}^m \frac{\hat{\lambda}_{\hat{\Sigma}_X,i}}{(\hat{\theta}_{X}-\hat{\lambda}_{\hat{\Sigma}_X,i})}, &
\hat{\theta}_{Y} \ \bigg| \ \frac{1}{\theta-1} = \frac{1}{m-1} \sum_{i=2}^m \frac{\hat{\lambda}_{\hat{\Sigma}_Y,i}}{(\hat{\theta}_{Y}-\hat{\lambda}_{\hat{\Sigma}_Y,i})}.
\end{array}$}
\end{equation*} 
We say that the criterion is satisfied, if this estimate of the expectation of the residual spike is a monotone increasing function of $\theta$. 
\end{Def}

\begin{Rem}
The above estimate of the expectation of a residual spike fails when $\theta$ is large compared to $m$ and we should then use an asymptotic estimator of $\alpha$ based on: 

\scalebox{0.9}{
\begin{minipage}{1\textwidth}
 \begin{eqnarray*}
 \alpha_{X}^2&=&1+\frac{1}{\theta} \left(1- \frac{1}{m} 
\sum_{i=2}^m \hat{\lambda}_{\hat{\Sigma}_X,i}^2 \right)\\
&&\hspace{0.3cm}
 + \frac{1}{\theta^2} 
 \left(1+\frac{2}{m} \sum_{i=2}^m \hat{\lambda}_{\hat{\Sigma}_X,i}^2 + \frac{3}{m^2} \left( \sum_{i=2}^m \hat{\lambda}_{\hat{\Sigma}_X,i}^2 \right)^2- \frac{2}{m} \sum_{i=2}^m \hat{\lambda}_{\hat{\Sigma}_X,i}^3\right)
 + O_p\left( \frac{1}{\theta^3}\right),\\
 \theta_{X}&=&(\theta-1)+ \frac{1}{m}\sum_{i=2}^m \hat{\lambda}_{\hat{\Sigma}_X,i}^2+O_p\left( \frac{1}{\theta^2} \right).
 \end{eqnarray*}
 \end{minipage}}\\
 
Using this approximation, the estimated curve of the equation (\ref{eq=CriteriaTsmall}) makes an error of $O_p\left(1/\theta^2 \right)$.
\end{Rem}

\section{Further Results} \label{sec:Theorems}
The proof of the main distributional result (Theorem \ref{TH=Main}) is based on three results and a lemma that are worthwhile on their own right and will be presented and discussed in this section. Future papers will also use this result for extensions.

\subsection{Unit invariant vector statistic}
\begin{Th}
 \label{Thunitstatistic} \ \\
Let $W$ be a random matrix with spectrum $S_W=\left\lbrace\hat{\lambda}_{W,1},\hat{\lambda}_{W,2},...,\hat{\lambda}_{W,m} \right\rbrace$ and trace equal to $m$. We denote by $u_{p_1}$ and $u_{p_2}$, two orthonormal invariant random vectors of size $m$ and independent of the eigenvalues of $W$. We set 
\begin{eqnarray*}
\vec{B}_m\left(\rho,\vec{s},\vec{r},\vec{p} \right)=\sqrt{m} \left( \begin{pmatrix}
\sum_{i=1}^m \frac{\hat{\lambda}_{W,i}^{s_1}}{\left(\rho-\hat{\lambda}_{W,i}\right)^{s_2}} {u}_{p_1,i} {u}_{p_2,i} \\ 
\sum_{i=1}^m \frac{\hat{\lambda}_{W,i}^{r_1}}{\left(\rho-\hat{\lambda}_{W,i}\right)^{r_2}} {u}_{p_1,i} {u}_{p_2,i} 
\end{pmatrix} - \begin{pmatrix}
M_{s_1,s_2} \\ 
M_{r_1,r_2}
\end{pmatrix} \mathbf{1}_{p_1=p_2}\right),
\end{eqnarray*}
where $\vec{s}=\left(s_1,s_2 \right)$, $\vec{r}=\left(r_1,r_2 \right)$ and $\vec{p}=\left(p_1,p_2 \right)$ with indices $1 \leqslant p_1 \leqslant p_2 \leqslant m$ and $s_1,s_2,r_1,r_2 \in \mathbf{N}$.\\
If $p=p_1=p_2$, 
\begin{equation*}
\resizebox{.95\hsize}{!}{$\left. \vec{B}_m\left(\rho,\vec{s},\vec{r},\vec{p} \right) \right| S_W \sim {\Normal}\left(\vec{0},  \begin{pmatrix}
2 \left(M_{2 s_1, 2 s_2}-M_{s_1,s_1}^2 \right) & 2 \left( M_{s_1+r_1,s_2+r_2}-M_{s_1,s_2} M_{r_1,r_2} \right) \\ 
2 \left( M_{s_1+r_1,s_2+r_2}-M_{s_1,s_2} M_{r_1,r_2} \right) & 2 \left(M_{2 r_1, 2 r_2}-M_{r_1,r_1}^2 \right)
\end{pmatrix} \right)+ o_{p}(1),$}
\end{equation*}
where $M_{s,r}=M_{s,r}(\rho)=\frac{1}{m}\sum_{i=1}^m \frac{\hat{\lambda}_{W,i}^s}{\left(\rho-\hat{\lambda}_{W,i}\right)^r} $.\\
Moreover, for $p_1 \neq p_2$,
\begin{equation*}
\resizebox{.95\hsize}{!}{$\left.\vec{B}_m\left(\rho,\vec{s},\vec{r},\vec{p} \right) \right| S_W \sim {\Normal}\left(\vec{0},  \begin{pmatrix}
 M_{2 s_1, 2 s_2}-M_{s_1,s_1}^2  &   M_{s_1+r_1,s_2+r_2}-M_{s_1,s_2} M_{r_1,r_2}  \\ 
  M_{s_1+r_1,s_2+r_2}-M_{s_1,s_2} M_{r_1,r_2}  & M_{2 r_1, 2 r_2}-M_{r_1,r_1}^2 
\end{pmatrix} \right)+ o_{p}(1),$}
\end{equation*}
In particular, with the notation $M_{s,0}=M_s=\frac{1}{m}\sum_{i=1}^m \hat{\lambda}_{W,i}^s $,
\begin{equation*}
\resizebox{.95\hsize}{!}{$\left. \sqrt{m} \left( \begin{pmatrix}
\sum_{i=1}^m \hat{\lambda}_{W,i} {u}_{p,i}^2 \\ 
\sum_{i=1}^m \hat{\lambda}_{W,i}^2 {u}_{p,i}^2
\end{pmatrix} - \begin{pmatrix}
1 \\ 
M_2
\end{pmatrix}\right) \right| S_W \sim \Normal\left(\vec{0},  \begin{pmatrix}
2 \left(M_2-1 \right) & 2 \left( M_3-M_2\right) \\ 
2 \left( M_3-M_2 \right) & 2 \left(M_4-M_2^2 \right)
\end{pmatrix} \right)+ o_{p}(1),$}
\end{equation*}
and 
\begin{equation*}
\resizebox{.95\hsize}{!}{$\left.\sqrt{m} \left( \begin{pmatrix}
\sum_{i=1}^m \hat{\lambda}_{W,i} {u}_{p_1,i} {u}_{p_2,i} \\ 
\sum_{i=1}^m \hat{\lambda}_{W,i}^2 {u}_{p_1,i} {u}_{p_2,i}
\end{pmatrix} - \begin{pmatrix}
0 \\ 
0
\end{pmatrix}\right) \right| S_W \sim \Normal \left(\vec{0},  \begin{pmatrix}
M_2-1  &   M_3-M_2 \\ 
  M_3-M_2  &  M_4-M_2^2 
\end{pmatrix} \right) + o_{p}(1).$}
\end{equation*} 
Finally if we look at $K$ bivariate normal random variables :
\begin{equation*}
\resizebox{.95\hsize}{!}{$\mathbf{B}_{m}\left(\vec{\rho},\mathbf{{s}},\mathbf{{r}},\mathbf{{p}} \right)=\left(\vec{B}_m\left(\rho_1,\vec{s}_1,\vec{r}_1,\vec{p}_1 \right), 
\vec{B}_m\left(\rho_2,\vec{s}_2,\vec{r}_2,\vec{p}_2 \right),
...,
\vec{B}_m\left(\rho_K,\vec{s}_K,\vec{r}_K,\vec{p}_K \right) \right),$}
\end{equation*}
where $\vec{p}_i \neq\vec{p}_j$ if $i\neq j$. Then, conditioning on the spectrum $S_W$, \linebreak $\mathbf{{B}}_{m}\left(\vec{\rho},\mathbf{{s}},\mathbf{{r}},\mathbf{{p}} \right)$ tends to a multivariate Normal. Moreover, all the bivariate elements $\vec{B}_m\left(\rho_i,\vec{s}_i,\vec{r}_i,\vec{p}_i \right)$ are asymptotically independent.
\begin{Remth}{\ }
\begin{enumerate}
\item The trace of $W$ equal to $m$ can easily be obtained by rescaling the matrix.
\item Although the condition of independence between eigenvectors and eigenvalues of $W$ appears to be restrictive, it is an automatic consequence if the eigenvectors are Haar distributed. 
\item If $W$ is a rescaled standard Wishart, then
\begin{eqnarray*}
\resizebox{.9\hsize}{!}{$\sqrt{m} \left( \begin{pmatrix}
\sum_{i=1}^m \hat{\lambda}_{W,i} \hat{u}_i^2 \\ 
\sum_{i=1}^m \hat{\lambda}_{W,i}^2 \hat{u}_i^2
\end{pmatrix} - \begin{pmatrix}
1 \\ 
1+c
\end{pmatrix}\right) \underset{m\rightarrow \infty }{\rightarrow } \Normal \left(\vec{0}, \begin{pmatrix}
2 c & 2c(2+c) \\ 
2c(2+c) & 2c(c+1)(c+4)
\end{pmatrix} \right),$}
\end{eqnarray*}
where $\hat{u}_1,\hat{u}_2,...,\hat{u}_m$ are the eigenvectors of $W$.
\end{enumerate}
\end{Remth}
\end{Th}
(Proof in supplement material \cite{Suppmaterial}.)

\subsection{Characterisation and convergence of eigenvalues and angles}

In this section, we study the convergence of the random variable $\hat{\hat{\theta}}_X$ and the angle between the eigenvectors. The proof for the parts  $2.a$ and $2.b$ are given in \cite{deformedRMT}, which also provides the main idea for the proof.
We only show convergence results for perturbations of order $k=1$, although we express eigenvectors and eigenvalues of a matrix $P^{1/2} W P^{1/2}$ as a function of the eigenstructure of $W$ in general and $W$ can already be perturbed.

\begin{Th} 
\label{convergence}{\ \\}
In this theorem, $P=\I_m+ (\theta-1) u u^t$ is a finite perturbation of order $1$.
\begin{enumerate}
\item Suppose $W$ is a symmetric matrix with eigenvalues $\hat{\lambda}_{W,i}\geqslant 0$ and eigenvectors $\hat{u}_{W,i}$ for $i=1,2,...,m$. The perturbation of $W$ by $P$ leads to $\hat{\Sigma}=P^{1/2}WP^{1/2}$. \\
For $i=1,2,...,m$, we define $\tilde{u}_{\hat{\Sigma},i}$ and $\hat{\lambda}_{\hat{\Sigma},i}$ such that 
$$ W P \tilde{u}_{\hat{\Sigma},i} = \hat{\lambda}_{\hat{\Sigma},i} \tilde{u}_{\hat{\Sigma},i},$$ 
and the usual $\hat{u}_{\hat{\Sigma},i}$ such that if $\hat{\Sigma}=P^{1/2}WP^{1/2}$, then
$$\hat{\Sigma} \hat{u}_{\hat{\Sigma},i} = P^{1/2} W P^{1/2} \hat{u}_{\hat{\Sigma},i} = \hat{\lambda}_{\hat{\Sigma},i} \hat{u}_{\hat{\Sigma},i}.$$
Under these conditions, the following results hold:
\begin{itemize}
\item The eigenvalues $\hat{\lambda}_{\hat{\Sigma},s}$ are such that for $s=1,2,...,m$,
\begin{eqnarray*}
\sum_{i=1}^m \frac{\hat{\lambda}_{W,i}}{\hat{\lambda}_{\hat{\Sigma},s}-\hat{\lambda}_{W,i}}  \left\langle \hat{u}_{W,i},u \right\rangle^2 = \frac{1}{\theta-1}.
\end{eqnarray*}
\item The eigenvectors $\tilde{u}_{\hat{\Sigma},s}$ are such that 
\begin{eqnarray*} 
&&  \left\langle \tilde{u}_{\hat{\Sigma},s},v \right\rangle^2  =\frac{\left( 
 \sum_{i=1}^m \frac{\hat{\lambda}_{W,i}}{\hat{\lambda}_{\hat{\Sigma},s}-\hat{\lambda}_{W,i}} \left\langle \hat{u}_{W,i},v \right\rangle \left\langle \hat{u}_{W,i},u \right\rangle  \right)^2}{
\sum_{i=1}^m \frac{\hat{\lambda}_{W,i}^2}{(\hat{\lambda}_{\hat{\Sigma},s}-\hat{\lambda}_{W,i})^2} \left\langle \hat{u}_{W,i},u \right\rangle^2 }. 
\end{eqnarray*}
In particular if $v=u$,
\begin{eqnarray*}
 \left\langle \tilde{u}_{\hat{\Sigma},s},u \right\rangle^2 &=& \frac{1}{\left(\theta-1 \right)^2 \left(\sum_{i=1}^m \frac{\hat{\lambda}_{W,i}^2}{(\hat{\lambda}_{\hat{\Sigma},s}-\hat{\lambda}_{W,i})^2} \left\langle \hat{u}_{W,i},u \right\rangle ^2 \right)}.
\end{eqnarray*}
Moreover, 
\begin{eqnarray*}
\hat{u}_{\hat{\Sigma},s}= \frac{P^{1/2} \tilde{u}_{\hat{\Sigma},s}} {\sqrt{1+ \left(\theta-1 \right) \left\langle \tilde{u}_{\hat{\Sigma},s},u \right\rangle^2}}.
\end{eqnarray*}
Therefore,  for $u$ and $v$ such that $\left\langle v,u \right\rangle=0$,
\begin{eqnarray*}
  \left\langle \hat{u}_{\hat{\Sigma},s},u \right\rangle^2 &=& \frac{ \theta \left\langle \tilde{u}_{\hat{\Sigma},s},u \right\rangle^2}{1+ \left(\theta-1 \right) \left\langle \tilde{u}_{\hat{\Sigma},s},u \right\rangle^2}
  =-\frac{\theta}{(\theta-1)^2 \hat{\lambda}_{\hat{\Sigma},s} T_{W,u}'(\hat{\lambda}_{\hat{\Sigma},s})},\\
  \left\langle \hat{u}_{\hat{\Sigma},s},v \right\rangle^2 &=&\frac{\left\langle \tilde{u}_{\hat{\Sigma},s},u \right\rangle^2}{1+ \left(\theta-1 \right) \left\langle \tilde{u}_{\hat{\Sigma},s},u \right\rangle^2},
\end{eqnarray*}
where $T_{W,u}(z)= \sum_{i=1}^m \frac{\hat{\lambda}_{W,i}}{z-\hat{\lambda}_{W,i}} \left\langle \hat{u}_{W,i},u \right\rangle^2$ is a weighted empirical T-transform.
\end{itemize}

\item Suppose that $W_X$, $W_Y$ and $P=P_X=P_Y$ satisfy the conditions described at the start of Section 2. Moreover, suppose that $\theta$ is large enough to create detectable spikes, $\left(\hat{\theta}_X,\hat{u}_X\right)$ and $\left(\hat{\theta}_Y,\hat{u}_Y\right)$, in the matrices $\hat{\Sigma}_X=P^{1/2}W_XP^{1/2}$ and $\hat{\Sigma}_Y=P^{1/2}W_Y P^{1/2}$. Then,
\begin{eqnarray*}
&a) &\hat{\hat{\theta}}_X,\hat{\hat{\theta}}_Y \underset{n,m \rightarrow \infty}{\overset{P}{\longrightarrow}} \theta,\\
&b) &\left\langle \hat{u}_X,u \right\rangle-\alpha_X,\left\langle \hat{u}_Y,u \right\rangle-\alpha_Y \underset{n,m \rightarrow \infty}{\overset{P}{\longrightarrow}} 0, \\
&c) &\left\langle \hat{u}_X,\hat{u}_Y \right\rangle-\alpha_X \alpha_Y \underset{n,m \rightarrow \infty}{\overset{P}{\longrightarrow}} 0,
\end{eqnarray*}
where
\begin{eqnarray*}
&&\hat{\theta}_X \underset{n,m \rightarrow \infty}{\overset{P}{\longrightarrow}} \rho_X,\\
&&\hat{\hat{\theta}}_{X} = 1+\frac{1}{\hat{T}_{\hat{\Sigma}_X}\left( \hat{\theta}_{X} \right)}
= 1+ \frac{m}{\sum_{i=k+1}^m \frac{\hat{\lambda}_{\hat{\Sigma}_X,i}}{\hat{\theta}_{X}-\hat{\lambda}_{\hat{\Sigma}_X,i}}},\\
&&\alpha_X^2 = -\frac{\theta}{(\theta-1)^2 \rho_X \hat{T}_{W_X}'(\rho_X)},\\ 
&&\alpha_Y^2 = -\frac{\theta}{(\theta-1)^2 \rho_Y \hat{T}_{W_Y}'(\rho_Y)},\\
&&\hat{\lambda}_{\hat{\Sigma}_X,i} \text{ and } \hat{\lambda}_{\hat{\Sigma}_Y,i} \text{ are the eigenvalues of respectively } \hat{\Sigma}_{X} \text{ and }\hat{\Sigma}_{Y}.
\end{eqnarray*}
Note that $\hat{T}_{W_X}(z)=\frac{1}{m} \sum_{i=1}^m \frac{\hat{\lambda}_{W_X,i}}{z-\hat{\lambda}_{W_X,i}}$ and $\hat{T}_{\hat{\Sigma}_X}(z)=\frac{1}{m} \sum_{i=k+1}^m \frac{\hat{\lambda}_{\hat{\Sigma}_X,i}}{z-\hat{\lambda}_{\hat{\Sigma}_X,i}}$ are empirical T-transform and its estimation.

\end{enumerate}

\begin{Remth}\ \\
If the spectra of $W_X$ and $W_Y$ are Wishart random matrices of size $m$ and degree of freedom $n_X$, $n_Y$ respectively, then by setting $c_X=\frac{m}{n_X}$ and $c_Y=\frac{m}{n_Y}$
\begin{eqnarray*}
&&\alpha_X^2 =\frac{1-\frac{c_X}{(\theta-1)^2}}{1+\frac{c_X}{\theta-1}},\\
&&\hat{\hat{\theta}}_{X} \text{ is such that } \hat{\theta}_X=\hat{\hat{\theta}}_{X}\left( 1+\frac{c_X}{\hat{\hat{\theta}}_{X}-1} \right), \text{ and} \\
&& \underset{m \rightarrow \infty}{\lim} \hat{\theta}_X = \theta\left( 1+\frac{c_X}{\theta-1} \right).
\end{eqnarray*}
\end{Remth}

\end{Th}
(Proof in supplement material \cite{Suppmaterial}.)

The second part of Theorem \ref{convergence} is very surprising! We already knew that the eigenvectors are not consistent. We show in the proof that the dot product of $\hat{u}_X$ and $\hat{u}_Y$ is smaller than that of $\hat{u}_X$ and $u$ and that of $\hat{u}_Y$ and $u$.
Among the consequences of this theorem is the fact that there is always an asymptotic bias between two eigenvectors, even if they are equal.

\subsection{Asymptotic distribution of the eigenvalues and the angle}

Suppose that you observe a perturbation of order $k=1$ applied to two random matrices $W_X \in \mathbb{R}^{m \times m }$ and $W_Y \in \mathbb{R}^{m \times m}$. We investigate the distribution of $\hat{\hat{\theta}}_X$, $\hat{\hat{\theta}}_Y$, $\left\langle \hat{u}_X,u \right\rangle^2$ and $\left\langle \hat{u}_X,\hat{u}_Y \right\rangle^2$.

\begin{Th} \label{jointdistribution}\ \\
Suppose $W_X$ and $W_Y$ satisfy the conditions described at the start of Section \ref{section:teststat} with $P=P_X=P_Y=\I_m+ (\theta-1) u u^t$, a detectable perturbation of order $k=1$. Moreover, we assume $S_{W_X}=\left\lbrace  \hat{\lambda}_{W_X,1},\hat{\lambda}_{W_X,2},...,\hat{\lambda}_{W_X,m} \right\rbrace$ and $S_{W_Y}=\left\lbrace \hat{\lambda}_{W_Y,1},\hat{\lambda}_{W_Y,2},...,\hat{\lambda}_{W_Y,m} \right\rbrace$, the eigenvalues of $W_X$ and $W_Y$ as known. We defined
\begin{eqnarray*}
\hat{\Sigma}_X&=&P^{1/2} W_X P^{1/2},\\
\hat{\Sigma}_Y&=&P^{1/2} W_Y P^{1/2}.
\end{eqnarray*}
We construct the unbiased estimators of $\theta$, $\hat{\hat{\theta}}_X$ and $\hat{\hat{\theta}}_Y$ via the relationship
\begin{equation*}
\frac{1}{\hat{\hat{\theta}}_X-1}= \frac{1}{m}\sum_{i=1}^m \frac{\hat{\lambda}_{W_X,i}}{\hat{\theta}_X-\hat{\lambda}_{W_X,i}} \  \text{ and } \frac{1}{\hat{\hat{\theta}}_Y-1}= \frac{1}{m}\sum_{i=1}^m \frac{\hat{\lambda}_{W_Y,i}}{\hat{\theta}_Y-\hat{\lambda}_{W_Y,i}}\,,
\end{equation*}
where $\hat{\theta}_X=\hat{\lambda}_{\hat{\Sigma}_X,1}$ and $\hat{\theta}_Y=\hat{\lambda}_{\hat{\Sigma}_Y,1}$ are the largest eigenvalues of $\hat{\Sigma}_X$ and $\hat{\Sigma}_Y$ corresponding to the eigenvectors $\hat{u}_X=\hat{u}_{\hat{\Sigma}_X,1}$ and $\hat{u}_Y=\hat{u}_{\hat{\Sigma}_Y,1}$.

\begin{enumerate}
\item If $\frac{\theta}{\sqrt{m}} \rightarrow 0$, we define
\begin{eqnarray*}
\resizebox{.9\hsize}{!}{$
 M_{s,r,X}\equiv M_{s,r,X}(\rho_X)= \frac{1}{m} \sum_{i=1}^m \frac{\hat{\lambda}_{W_X,i}^s}{(\rho_X-\hat{\lambda}_{W_X,i})^r}, \  M_{s,r,X}\equiv M_{s,r,Y}(\rho_Y)=\frac{1}{m} \sum_{i=1}^m \frac{\hat{\lambda}_{W_X,i}^s}{(\rho_Y-\hat{\lambda}_{W_X,i})^r}\,,  
$}
\end{eqnarray*}
where we assume
\begin{eqnarray*}
\rho_X =\E\left[ \hat{\theta}_X \right]+o\left(\frac{\theta}{\sqrt{m}}\right),
&\rho_Y =\E\left[ \hat{\theta}_Y \right]+o\left(\frac{\theta}{\sqrt{m}}\right)
\end{eqnarray*}
and a convergence rate of $(\hat{\theta}_X,\hat{\theta}_Y)$ to $\left(\rho_X,\rho_Y \right)$ in $O_p\left(1/\sqrt{m}\right)$. Then
\begin{equation*}
 \resizebox{.9\hsize}{!}{$\left. \begin{pmatrix}
\hat{\hat{\theta}}_X  \\ 
\left\langle \hat{u}_X,u \right\rangle^2 
\end{pmatrix} \right| S_{W_X}
\sim  \Normal \left(
\begin{pmatrix}
\theta  \\ 
\alpha_X^2
\end{pmatrix}
,
 \frac{1}{m}\begin{pmatrix}
\sigma_{\theta,X}^2 & \sigma_{\theta,\alpha^2,X} \\ 
\sigma_{\theta,\alpha^2,X} & \sigma_{\alpha^2,X}^2
\end{pmatrix} \right)+\begin{pmatrix}
o_p\left(\frac{\theta}{\sqrt{m}}\right) \\ 
o_p\left( \frac{1}{\theta \sqrt{m}}\right)
\end{pmatrix},$}
\end{equation*}
\begin{eqnarray*}
\resizebox{.9\hsize}{!}{$\left. \begin{pmatrix}
\hat{\hat{\theta}}_X  \\ 
\hat{\hat{\theta}}_Y  \\ 
\left\langle \hat{u}_X,\hat{u}_Y \right\rangle^2
\end{pmatrix} \right| S_{W_X}, S_{W_Y} \sim \Normal \left(
\begin{pmatrix}
\theta \\
\theta  \\ 
\alpha_{X,Y}^2
\end{pmatrix}
,\frac{1}{m}\begin{pmatrix}
\sigma_{\theta,X}^2 & 0 & \sigma_{\theta,\alpha^2,X} \\ 
0 & \sigma_{\theta,Y}^2  &  \sigma_{\theta,\alpha^2,Y} \\
\sigma_{\theta,\alpha^2,X} & \sigma_{\theta,\alpha^2,Y} & \sigma_{\alpha^2,X,Y}^2
\end{pmatrix}
 \right)+\begin{pmatrix}
o_p\left(\frac{\theta}{\sqrt{m}}\right) \\ 
o_p\left(\frac{\theta}{\sqrt{m}}\right) \\ 
o_p\left( \frac{1}{\theta \sqrt{m}}\right)
\end{pmatrix},$}
\end{eqnarray*}

where  \\
\scalebox{0.67}{
\begin{minipage}{1\textwidth}
\begin{eqnarray*}
&&\alpha_X^2=\frac{\theta}{(\theta-1)^2} \frac{1}{\rho_X M_{1,2,X}},\\
&&\alpha_{X,Y}^2=\frac{\theta^2}{(\theta-1)^4} \frac{1}{\rho_X \rho_Y M_{1,2,X}M_{1,2,Y}},\\
&&\sigma_{\theta,X}^2 =\frac{2 \left(M_{2,2,X}-M_{1,1,X}^2\right)}{M_{1,1,X}^4}, \hspace{30cm}
\end{eqnarray*} 
\begin{eqnarray*}
&&\sigma_{\alpha^2,X}^2=\frac{2 \theta^2}{\left((\theta-1) \rho_X M_{1,2,X}\right)^4} \left( \rho_X^2 \left(M_{2,4,X}-M_{1,2,X}^2 \right) + \left( 2 \rho_X \frac{M_{1,3,X}}{M_{1,2,X}}-1 \right)^2 \left( M_{2,2,X}-M_{1,1,X}^2 \right) \right.\\
&&\hspace{2cm} \left.- 2 \rho_X \left( 2 \rho_X \frac{M_{1,3,X}}{M_{1,2,X}}-1 \right) \left( M_{2,3,X}-\frac{M_{1,1,X}}{ M_{1,2,X}} \right) \right),\\
&&\sigma_{\theta,\alpha^2,X}= \frac{2 \theta}{M_{1,1,X}^2 M_{1,2,X}^3 \rho_X^2 (-1+\theta)^2}    \bigg(M_{1,1,X} M_{1,2,X}^2 \rho_X+2 M_{1,3,X} M_{2,2,X} \rho_X \\
&&\hspace{2cm} +M_{1,1,X}^2 (M_{1,2,X}-2 M_{1,3,X} \rho_X)-M_{1,2,X} (M_{2,2,X}+M_{2,3,X} \rho_X) \bigg) ,\\
&&\sigma_{\alpha^2,X,Y}^2=\sigma_{\alpha^2,X}^2 \alpha_Y^4 + \sigma_{\alpha^2,Y}^2 \alpha_X^4+4 \alpha_{X,Y}^2 (1-\alpha_X^2)(1-\alpha_X^2).\hspace{30cm}
\end{eqnarray*}
\end{minipage}}\\

\item  If $\frac{\theta}{\sqrt{m}} \rightarrow \infty$, then we can simplify the formulas.
We define
\begin{eqnarray*}
M_{r,X}=\frac{1}{m} \sum_{i=1}^m \hat{\lambda}_{W_X,i}^r  &\text{ and }&
M_{r,Y}=\frac{1}{m} \sum_{i=1}^m \hat{\lambda}_{W_Y,i}^r .
\end{eqnarray*}
Using this notation,\\
\scalebox{0.55}{
\begin{minipage}{1\textwidth}
\begin{eqnarray*}
 && \left. \begin{pmatrix}
\hat{\hat{\theta}}_X  \\ 
\left\langle \hat{u}_X,u \right\rangle^2 
\end{pmatrix} \right| S_{W_X}
\sim  \Normal \left(
\begin{pmatrix}
\theta + O_p\left(1\right) \\ 
1+\frac{1-M_{2,X}}{\theta} +O_p\left( \frac{1}{\theta^2} \right)
\end{pmatrix}
,\right.\\
&&\hspace{2.3cm} \left.
 \frac{1}{m}\begin{pmatrix}
2 \theta^2 \left( 1-M_{2,X} \right) & 2 \left( 2 M_{2,X}^2- M_{2,X}- M_{3,X}\right) \\ 
2 \left( 2 M_{2,X}^2- M_{2,X}- M_{3,X}\right) & \frac{2}{\theta^2}\left( 4 M_{2,X}^3- M_{2,X}^2-4 M_{2,X} M_{3,X}+ M_{4,X} \right) 
\end{pmatrix} \right)+\begin{pmatrix}
o_p\left(\frac{\theta}{\sqrt{m}}\right) \\ 
o_p\left( \frac{1}{\theta \sqrt{m}}\right)
\end{pmatrix},
\end{eqnarray*}
\begin{eqnarray*}
&& \hspace{-0.5cm} \left. \begin{pmatrix}
\hat{\hat{\theta}}_X  \\ 
\hat{\hat{\theta}}_Y  \\ 
\left\langle \hat{u}_X,\hat{u}_Y \right\rangle^2
\end{pmatrix} \right| S_{W_X}, S_{W_Y} \sim \Normal\left(
\begin{pmatrix}
\theta + O_p\left(1 \right)\\
\theta + O_p\left(1 \right) \\ 
1+\frac{2-M_{2,X}-M_{2,Y}}{\theta} +O_p\left( \frac{1}{\theta^2} \right)
\end{pmatrix}
,\right.\\
&&\hspace{1.8cm}\left.\frac{1}{m}\begin{pmatrix}
2 \theta^2 \left( 1-M_{2,X} \right) & 0 & 2 \left( 2 M_{2,X}^2- M_{2,X}- M_{3,X}\right) \\ 
0 & 2 \theta^2 \left( 1-M_{2,Y} \right) & 2 \left( 2 M_{2,Y}^2- M_{2,Y}- M_{3,Y}\right) \\
2 \left( 2 M_{2,X}^2- M_{2,X}- M_{3,X}\right)& 2 \left( 2 M_{2,Y}^2- M_{2,Y}- M_{3,Y}\right)  & \frac{S}{\theta^2}
\end{pmatrix}
 \right)+\begin{pmatrix}
o_p\left(\frac{\theta}{\sqrt{m}}\right) \\ 
o_p\left(\frac{\theta}{\sqrt{m}}\right) \\ 
o_p\left( \frac{1}{\theta \sqrt{m}}\right)
\end{pmatrix},\\
 &&S= 2\left( 4 M_{2,X}^3- M_{2,X}^2-4 M_{2,X} M_{3,X}+ M_{4,X} \right)+2\left( 4 M_{2,Y}^3- M_{2,Y}^2-4 M_{2,Y} M_{3,Y}+ M_{4,Y} \right)+ 4 (M_{2,Y}-1) (M_{2,X}-1).
\end{eqnarray*}
\end{minipage}}\\

Moreover, the asymptotic distributions of $\hat{\theta}_X$ and $\hat{\hat{\theta}}_X$ are the same.

\item If $\frac{\theta}{\sqrt{m}} \rightarrow d$, a finite constant, then a mixture of the two first scenarios describes the first two moments of the joint distribution.\\
The formula of the second moment is asymptotically the same as the variance formula when $\frac{\theta}{\sqrt{m}} \rightarrow \infty$.\\
The formula of the first moment is asymptotically the same as the expectation formula when $ \frac{\theta}{\sqrt{m}} \rightarrow 0$.\\

\item The random variables can be expressed as functions of invariant unit random statistics of the form:
\begin{eqnarray*}
\overset{u}{M}_{r,s,X}(\rho)&=& \sum_{i=1}^m \frac{\hat{\lambda}^r_{W_X,i}}{\left(\rho-\hat{\lambda}_{W_X,i}\right)^s} \left\langle \hat{u}_{W_X,i},u \right\rangle^2.
\end{eqnarray*}
(Assuming a canonical perturbation leads to a simpler formula)\\
Knowing $S_{W_X}$ and $S_{W_Y}$ we have
\begin{itemize}
\item Exact distributions:
\begin{eqnarray*}
&&\hat{\theta}_X   \left|  \frac{1}{\theta-1}=\overset{  u}{M}_{1,1,X}\left( \hat{\theta}_X \right) \right. ,\\
&&\left\langle \hat{u}_X,e_1 \right\rangle^2 = \frac{\theta}{(\theta-1)^2} \frac{1}{\hat{\theta}_X \overset{u}{M}_{1,2,X}\left( \hat{\theta}_X \right)}.
\end{eqnarray*}
Moreover, \small
\begin{eqnarray*}
\left\langle \hat{u}_X,\hat{u}_Y \right\rangle&=& \left\langle \hat{u}_X,e_1 \right\rangle \left\langle \hat{u}_Y,e_1 \right\rangle+ \sqrt{1-\left\langle \hat{u}_X,e_1 \right\rangle^2} \sqrt{1-\left\langle \hat{u}_Y,e_1 \right\rangle^2} Z,\\
\sum_{i=2}^m \hat{u}_{X,1} \hat{u}_{Y,1} &=& \sqrt{1-\left\langle \hat{u}_X,e_1 \right\rangle^2} \sqrt{1-\left\langle \hat{u}_Y,e_1\right\rangle^2} Z, \\
Z &\sim & \Normal\left(0,\frac{1}{m}\right)+ O_p\left(\frac{1}{m} \right),
\end{eqnarray*} \normalsize
where $Z$ is independent of $\left\langle \hat{u}_X,e_1 \right\rangle$, $\left\langle \hat{u}_Y,e_1 \right\rangle$, $\hat{\theta}_X$ and $\hat{\theta}_Y$. In order to get the exact distribution, we should replace $Z$ by $\sum_{i=1}^{m-1} v_i\tilde{v}_i$ where $v_i $ and $\tilde{v}_i$ are independent unit invariant random vectors.
\item Approximations:
\begin{eqnarray*}
\hat{\theta}_X &=&\rho +  \frac{\left(\overset{  u}{M}_{1,1,X}\left(\rho\right)- {M}_{1,1,X}(\rho)  \right)}{M_{1,2,X}\left(\rho\right)} + O_p\left( \frac{\theta}{m} \right)\\
&=& \theta \overset{u}{M}_{1,X} +O_p\left(1\right),\\
\hat{\hat{\theta}}_X &=& \theta + (\theta-1)^2 \left(\overset{ u}{M}_{1,1,X}\left(\rho\right)- {M}_{1,1,X}(\rho_X)  \right)+O_p\left(\frac{\theta}{ m} \right).
\end{eqnarray*}
We provide three methods of estimation of the angle in order to estimate it for all $\theta$.

\scalebox{0.68}{
\begin{minipage}{1\textwidth}
\begin{eqnarray*}
\left\langle \hat{u}_X,e_1 \right\rangle^2&=&\frac{\theta}{(\theta-1)^2} \left( \frac{1}{\rho_X M_{1,2,X}(\rho_X)} + \left( \frac{2 M_{1,3,X}(\rho_X)}{M_{1,2,X}(\rho_X)}-\frac{1}{\rho_X}\right) \frac{\overset{u}{M}_{1,1,X}(\rho_X)-\frac{1}{\theta-1}}{\left( M_{1,2,X}(\rho_X) \right)^2} \right.\\
&&\hspace{2cm} \left. -\frac{\overset{u}{M}_{1,2,X}(\rho_X)-M_{1,2,X}(\rho_X)}{\rho_X \left(M_{1,2,X}(\rho_X) \right)^2} \right) + O_p\left( \frac{1}{m} \right),\\
&=& 1+\frac{1}{\theta} \left(1- \overset{u}{M}_{2,X}+2 M_{2,X} \left( \overset{u}{M}_{1,X}-1 \right)\right)\\
&&\hspace{2cm}+ \frac{1}{\theta^2} \left(1- 2 \overset{u}{M}_{2,X}+3 \overset{u}{M}_{2,X}^2-2 \overset{u}{M}_{3,X} \right)+O_p\left(\frac{1}{\theta^3}\right)+O_p\left( \frac{1}{\theta m} \right),\\
&=& 1+\frac{1}{\theta}-\frac{1}{\theta} \overset{u}{M}_{2,X}+ \frac{2}{\theta} M_{2,X} \left( \overset{u}{M}_{1,X}-1 \right)+O_p\left( \frac{1}{\theta^2} \right)+O_p\left( \frac{1}{\theta m} \right)
\end{eqnarray*}
\end{minipage}}\\

Finally, the double angle is such that \\
\scalebox{0.85}{
\begin{minipage}{1\textwidth}
\begin{eqnarray*}
\left\langle \hat{u}_X,\hat{u}_Y \right\rangle= \left\langle \hat{u}_X,e_1 \right\rangle \left\langle \hat{u}_Y,e_1 \right\rangle+ \frac{\sqrt{M_{2,X}-1} \sqrt{M_{2,Y}-1}}{\theta} Z
+O_p\left( \frac{1}{\theta^2 \sqrt{m}} \right).
\end{eqnarray*}
\end{minipage}}
\end{itemize}

\end{enumerate}

\begin{Remth}
If the spectra of $W_X$ and $W_Y$ are rescaled Wishart matrices of size $m$ with $n$ degree of freedom. By setting $c=\frac{m}{n}$,\\
\scalebox{0.78}{
\begin{minipage}{1\textwidth}
$$ \begin{pmatrix}
\hat{\hat{\theta}}_X \\ 
\left\langle \hat{u}_X,u_0 \right\rangle^2 
\end{pmatrix} \overset{Asy}{\sim} \Normal \left(\begin{pmatrix}
\theta \\ 
\frac{1-\frac{c}{(\theta-1)^2}}{1+\frac{c}{\theta-1}} 
\end{pmatrix} , \frac{1}{m}\left(
\begin{array}{cc}
-\frac{2 c (\theta-1)^2 \theta^2}{c-(\theta-1)^2} & -\frac{2 c^2 (\theta-1) \theta^3}{\left(c-(\theta-1)^2\right) (c+\theta-1)^2} \\
-\frac{2 c^2 (\theta-1) \theta^3}{\left(c-(\theta-1)^2\right) (c+\theta-1)^2} & -\frac{2 c^2 \theta^2 \left(c^2+(\theta (\theta+2)-2) c+(\theta-1)^2\right)}{\left(c-(\theta-1)^2\right) (c+\theta-1)^4}  \\
\end{array}
\right)\right)$$ 
\end{minipage}}\\
 and \\
\scalebox{0.63}{
\begin{minipage}{1\textwidth}
\begin{eqnarray*}
&&\hspace{-0.5cm}\begin{pmatrix}
\hat{\hat{\theta}}_X \\ 
\hat{\hat{\theta}}_Y  \\ 
\left\langle \hat{u}_X,\hat{u}_Y \right\rangle^2
\end{pmatrix} \overset{Asy}{\sim} \Normal \left(\begin{pmatrix}
\theta \\ 
\theta\\
\left( \frac{1-\frac{c}{(\theta-1)^2}}{1+\frac{c}{\theta-1}}\right)^2 
\end{pmatrix} \right. , \\
&& \hspace{1.5cm}\left. \frac{1}{m} \begin{pmatrix}
-\frac{2 c (\theta-1)^2 \theta^2}{c-(\theta-1)^2}  & 0 & -\frac{2 c^2 (\theta-1) \theta^3}{\left(c-(\theta-1)^2\right) (c+\theta-1)^2}  \\ 
0 & -\frac{2 c (\theta-1)^2 \theta^2}{c-(\theta-1)^2} & -\frac{2 c^2 (\theta-1) \theta^3}{\left(c-(\theta-1)^2\right) (c+\theta-1)^2}\\
-\frac{2 c^2 (\theta-1) \theta^3}{\left(c-(\theta-1)^2\right) (c+\theta-1)^2} & -\frac{2 c^2 (\theta-1) \theta^3}{\left(c-(\theta-1)^2\right) (c+\theta-1)^2} & \frac{4 c^2 \theta^2 \left(c-(\theta-1)^2\right)^2 \left(c^3+4 c^2 (\theta-1)+c (\theta-1) (\theta (\theta+5)-5)+2 (\theta-1)^3\right)}{(\theta-1)^4 (c+\theta-1)^7} 
\end{pmatrix}  \right). 
\end{eqnarray*}  
 \end{minipage}}\\
 
 If $\theta$ tends to infinity, then
{\small
$$ \begin{pmatrix}
\hat{\hat{\theta}}_X \\ 
\left\langle \hat{u}_X,u \right\rangle^2 
\end{pmatrix} \overset{Asy}{\sim} \Normal \left(\begin{pmatrix}
\theta \\ 
\frac{1-\frac{c}{(\theta-1)^2}}{1+\frac{c}{\theta-1}} 
\end{pmatrix} , \frac{1}{m}\begin{pmatrix}
2 c \theta^2& 2c^2 \\ 
2c^2 & \frac{2 c^2 (c+1)}{\theta^2} 
\end{pmatrix}\right).$$}
Moreover,{\small
$$ \begin{pmatrix}
\hat{\hat{\theta}}_X \\ 
\hat{\hat{\theta}}_Y  \\ 
\left\langle \hat{u}_X,\hat{u}_X \right\rangle^2
\end{pmatrix} \overset{Asy}{\sim} \Normal \left(\begin{pmatrix}
\theta \\ 
\theta\\
\frac{1-\frac{c}{(\theta-1)^2}}{1+\frac{c}{\theta-1}} 
\end{pmatrix}, \frac{1}{m}\begin{pmatrix}
2c \theta^2 & 0 & 2c^2  \\ 
0 & 2c \theta^2 & 2c^2\\
2c^2 & 2c^2 & \frac{4 c^2 (c+2)}{\theta^2} 
\end{pmatrix} \right).$$}

\end{Remth}
\end{Th}
(Proof in supplement material \cite{Suppmaterial}.)
\subsection{Residual spike as a function of the statistics}
Finally we present a simple result of linear algebra that express the residual spike as a function of the statistics.
\begin{Lem} \label{Lemrespiketheta} \ \\
Suppose 
\begin{equation*}
\resizebox{.95\hsize}{!}{$D=\left(\I_m+(\theta-1) u_X u_X^t \right)^{-1/2} \left(\I_m+(\theta-1) u_Y u_Y^t \right) \left(\I_m+(\theta-1) u_X u_X^t \right)^{-1/2}.$}
\end{equation*}
The eigenvalues of $D$ are $1$ and 

\begin{eqnarray*}
\resizebox{.95\hsize}{!}{$\lambda\left( D \right)=- \frac{1}{2 \theta} \left(-1+\alpha^2-2 \alpha^2 \theta -\theta^2 (1-\alpha^2) \pm \sqrt{-4 \theta^2 + \left[1+ \theta^2-(-1+\theta)^2\alpha^2 \right]^2} \right),$}
\end{eqnarray*}
where $\alpha^2=\left\langle u_X,u_Y \right\rangle^2$.\\
Moreover, if 
\begin{equation*}
\resizebox{.95\hsize}{!}{$D_2=\left(\I_m+(\theta_X-1) u_X u_X^t \right)^{-1/2} \left(\I_m+(\theta_Y-1) u_Y u_Y^t \right) \left(\I_m+(\theta_X-1) u_X u_X^t \right)^{-1/2}.$}
\end{equation*}
The eigenvalues of $D_2$ are $1$ and 
\begin{eqnarray*}
\resizebox{.95\hsize}{!}{$\lambda\left( D_2 \right)=\frac{1}{2} \left(\theta_Y+\alpha^2-\theta_Y \alpha^2+\frac{1+(\theta_Y-1)\alpha^2 \pm \sqrt{-4 \theta_Y \theta_X+ \left( 1+\theta_Y \theta_X - (\theta_Y-1) (\theta_X-1) \alpha^2 \right)^2}}{\theta_X} \right),$}
\end{eqnarray*}

where $\alpha^2=\left\langle u_X,u_Y \right\rangle^2$.
\end{Lem}
(Proof in supplement material \cite{Suppmaterial}.)

\section{Comparison with existing tests}
In the classical multivariate theory, \cite{multi22} proposes a log-ratio test for the equality of two covariance matrices.\\
Suppose
\begin{eqnarray*}
X_1,X_2,...,X_{n_X} \overset{i.i.d.}\sim \Normal_m(0,\Sigma_X),\\
Y_1,Y_2,...,Y_{n_Y} \overset{i.i.d.}\sim \Normal_m(0,\Sigma_Y).
\end{eqnarray*}
We want to test 
\begin{eqnarray*}
{\rm H}_0 : \Sigma_X=\Sigma_Y,\\
{\rm H}_1 : \Sigma_X \neq \Sigma_Y,
\end{eqnarray*} 
The log-likelihood ratio test look at the statistic 
\begin{eqnarray*}
T_1=n_X \log\left( \left| \frac{n_X}{n_X+n_Y} \I_m+\frac{n_Y}{n_X+n_Y} \hat{\Sigma}_X^{-1/2} \hat{\Sigma}_Y \hat{\Sigma}_X^{-1/2} \right| \right).
\end{eqnarray*}
Under ${\rm H}_0$ and if $m$ is finite, $T_1 \overset{D}{\rightarrow} \chi^2_{\frac{p(p-1)}{2}}$. Some other interesting tests propose to observe the determinant and the trace of $\hat{\Sigma}_X^{-1/2} \hat{\Sigma}_Y \hat{\Sigma}_X^{-1/2}$.\\
In this section we show that any test statistics using $T_2=\log\left| \hat{\Sigma}_X^{-1/2} \hat{\Sigma}_Y \hat{\Sigma}_X^{-1/2} \right|$ or $T_3=\Tr \left( \hat{\Sigma}_X^{-1/2} \hat{\Sigma}_Y \hat{\Sigma}_X^{-1/2} \right)$ have difficulties to test 
\begin{eqnarray*}
{\rm H}_0 : P_X=P_Y,\\
{\rm H}_1 : P_X \neq P_Y,
\end{eqnarray*}
when $P_X$ and $P_Y$ are finite perturbations.\\
We compare the performance of these tests with our procedure $T$ defined in Section \ref{sec=test} in the table \ref{Tabletest}.
\begin{enumerate}
\item When $T=\left(\lambda_{\min} \left(\hat{\hat{\Sigma}}_X^{-1/2}  \hat{\hat{\Sigma}}_Y \hat{\hat{\Sigma}}_X^{-1/2} \right),\lambda_{\max} \left( \hat{\hat{\Sigma}}_X^{-1/2}  \hat{\hat{\Sigma}}_Y \hat{\hat{\Sigma}}_X^{-1/2} \right) \right)$, the table shows 
\begin{equation*}
\resizebox{.9\hsize}{!}{$P_{{\rm H}_1}\left( \lambda_{\min} \left(\hat{\hat{\Sigma}}_X^{-1/2}  \hat{\hat{\Sigma}}_Y \hat{\hat{\Sigma}}_X^{-1/2} \right) < q_{\lambda_{\min},H_0} (0.025) \text{ or } \lambda_{\max} \left(\hat{\hat{\Sigma}}_X^{-1/2}  \hat{\hat{\Sigma}}_Y \hat{\hat{\Sigma}}_X^{-1/2} \right) < q_{\lambda_{\max},H_0} (0.975) \right),$}
\end{equation*}
where $q_{\lambda_{\min},H_0}$ and $q_{\lambda_{\max},H_0}$ give the quantiles of $T$ under ${\rm H}_0$ and are given in Theorem \ref{TH=Main}.
\item When $T_2=\log\left| \hat{\Sigma}_X^{-1/2} \hat{\Sigma}_Y \hat{\Sigma}_X^{-1/2} \right|$, the table show 
\begin{equation*}
P_{{\rm H}_1}\left( T_2 < q_{T_2,H_0}(0.025) \text{ or } T_2 > q_{T_2,H_0}(0.975) \right),
\end{equation*}
where $q_{T_2,H_0}$ gives the quantiles of $T_2$ under ${\rm H}_0$ and is found empirically.
\item When $T_3=\Tr \left( \hat{\Sigma}_X^{-1/2} \hat{\Sigma}_Y \hat{\Sigma}_X^{-1/2} \right)$, the table show 
\begin{equation*}
P_{{\rm H}_1}\left( T_3 < q_{T_3,H_0}(0.025) \text{ or } T_3 > q_{T_3,H_0}(0.975) \right),
\end{equation*}
where $q_{T_3,H_0}$ gives the quantiles of $T_3$ under ${\rm H}_0$ is are found empirically.
\end{enumerate}
\begin{Rem}
In order to generalise the test to degenerated matrices, the determinant is defined as the product of the non-null eigenvalues of the matrix and the inverse is the generalised inverse. 
\end{Rem}

\begin{table}[hbt] \scalebox{0.7}{
\begin{tabular}{cccccc}
\begin{minipage}{1cm}
\begin{eqnarray*}
&&m=500,\\
&&n_X=n_Y=250
\end{eqnarray*}
\end{minipage}
 & \\
&\begin{minipage}{3cm}
\begin{eqnarray*}
&&\theta_X=7, u_X=e_1,\\
&&\theta_Y=7, u_Y=e_2,
\end{eqnarray*}
\end{minipage} &
\begin{minipage}{3cm}
\begin{eqnarray*}
&&\theta_X=50, u_X=e_1,\\
&&\theta_Y=50, u_Y=e_2,
\end{eqnarray*}
\end{minipage}  &
\begin{minipage}{3cm}
\begin{eqnarray*}
&&\theta_X =7, u_X=e_1,\\
&&\theta_Y =17, u_Y=e_1,
\end{eqnarray*}
\end{minipage} &
\begin{minipage}{3cm}
\begin{eqnarray*}
&&\theta_X=300, u_X=e_1,\\
&&\theta_Y=600, u_Y=e_1,
\end{eqnarray*} 
\end{minipage} &\\
\hline
$T$  & 0.81 & 1 & 0.98 & 0.99\\
$T_2$ & 0.15 & 0.05 & 0.11 & 0.05 \\
$T_3$ & 0.11 & 1 & 0.13 & 0.1\\
\hline 
\hline
\begin{minipage}{1cm}
\begin{eqnarray*}
&&m=500,\\
&&n_X=1000,\\
&&n_Y=250
\end{eqnarray*}
\end{minipage}
 & \\
&\begin{minipage}{3cm}
\begin{eqnarray*}
&&\theta_X=5, u_X=e_1,\\
&&\theta_Y=5, u_Y=e_2,
\end{eqnarray*}
\end{minipage} &
\begin{minipage}{3cm}
\begin{eqnarray*}
&&\theta_X=50, u_X=e_1,\\
&&\theta_Y=50, u_Y=e_2,
\end{eqnarray*}
\end{minipage}  &
\begin{minipage}{3cm}
\begin{eqnarray*}
&&\theta_X =5, u_X=e_1,\\
&&\theta_Y =15, u_Y=e_1,
\end{eqnarray*}
\end{minipage} &
\begin{minipage}{3cm}
\begin{eqnarray*}
&&\theta_X=300, u_X=e_1,\\
&&\theta_Y=600, u_Y=e_1,
\end{eqnarray*} 
\end{minipage} &\\
\hline
$T$ & 1 & 1 & 0.99 & 1 \\
$T_2$ & 0.09 & 0.31 & 0.04 & 0.11 \\
$T_3$ & 0.47 & 1 & 0.15 & 0.12 \\
\hline 
\hline
\begin{minipage}{1cm}
\begin{eqnarray*}
&&m=500,\\
&&n_X=1000,\\
&&n_Y=1000
\end{eqnarray*}
\end{minipage}
 & \\
&\begin{minipage}{3cm}
\begin{eqnarray*}
&&\theta_X=5, u_X=e_1,\\
&&\theta_Y=5, u_Y=e_2,
\end{eqnarray*}
\end{minipage} &
\begin{minipage}{3cm}
\begin{eqnarray*}
&&\theta_X=50, u_X=e_1,\\
&&\theta_Y=50, u_Y=e_2,
\end{eqnarray*}
\end{minipage}  &
\begin{minipage}{3cm}
\begin{eqnarray*}
&&\theta_X =5, u_X=e_1,\\
&&\theta_Y =15, u_Y=e_1,
\end{eqnarray*}
\end{minipage} &
\begin{minipage}{3cm}
\begin{eqnarray*}
&&\theta_X=300, u_X=e_1,\\
&&\theta_Y=600, u_Y=e_1,
\end{eqnarray*} 
\end{minipage} &\\
\hline
$T$ & 1 & 1 & 1 & 1\\
$T_2$ & 0.07 & 0.12 & 0.07 & 0.02\\
$T_3$ & 0.41 & 1 & 0.08 & 0.05
\end{tabular}}
\caption{Probability to detect the alternative with a test at level $0.05$ when $P_X=\I_m+(\theta_X-1) u_X u_X^t$ and $P_Y=\I_m+(\theta_Y-1) u_Y u_Y^t$ for the different tests. The distribution of $T_2$ and $T_3$ is computed empirically by assuming the same perturbation $P_X$ for the two groups.} \label{Tabletest}
\end{table}
In the particular case of finite perturbation, the trace and the determinant have difficulties to catch the alternative. On the other hand, our procedure detects easily some small differences.\\
The statistic $T_2$ and $T_3$ would be interesting to detect perturbation of large order such as a global change of the variance.

\begin{Rem}
Assuming $\hat{\Sigma}_X=P^{1/2}_X W_X P^{1/2}_X$ and $\hat{\Sigma}_Y=P^{1/2}_Y W_Y P^{1/2}_Y$ are as described at the start of Section  \ref{section:teststat}, the procedure proposed in this paper required the estimation of $M_{s,X}=\frac{1}{m}\sum_{i=1}^m \hat{\lambda}_{W_X,i}^s$ and $M_{s,Y}=\frac{1}{m} \sum_{i=1}^m\hat{\lambda}_{W_Y,i}^s$ for $s=1,2,3,4$ in order to compute the quantile under $ \rm{H}_0$ of the residual spikes, $q_{\lambda_{\min},H_0}$ and $q_{\lambda_{\max},H_0}$. By Cauchy-Interlacing law and bounded eigenvalues $\hat{\lambda}_{W_X,i}$ and $\hat{\lambda}_{W_Y,i}$ we can use the following estimator
\begin{eqnarray*}
\hat{M}_{s,X}=\frac{1}{m-1}\sum_{i=2}^m \hat{\lambda}_{\hat{\Sigma}_X,i}^s= M_{s,X} +O\left(\frac{1}{m}\right).
\end{eqnarray*}

\end{Rem}

\subsection{Conclusion}
By studying perturbation of order $1$, this work highlights the particular behaviour of residual spikes. A future work will present the behaviours of residual spikes when the perturbations are of order $k$. Nevertheless, this task requires many intermediary results. Therefore an other future work will present only the joint distribution of some statistics as the eigenvalues and the eigenvectors.

\appendix

\section{Table} \label{sec:Tablecomplete}
We extend the simulations of Section \ref{sec:smallsimulationmainth}. We test our Main Theorem \ref{TH=Main} under different hypotheses on $\X\in\mathbf{R}^{m \times n_X}$ and $\Y\in\mathbf{R}^{m \times n_Y}$ (recall that $W_X=\frac{1}{n_X}\X \X^t$ and $W_Y=\frac{1}{n_Y}\Y \Y^t$):
\begin{enumerate}
\item The matrices $\X$ and $\Y$ contain independent standard normal entries.
\item The columns of the matrices $\X$ and $\Y$ are i.i.d. Multivariate Student with $8$ degrees of freedom. For $i=1,2,...,n_X$ and $j=1,2...,n_Y$,
$$X_{\cdot,i} \overset{i.i.d.}{\sim} \frac{\Normal\left(\vec{0},\I_m\right)}{\sqrt{\frac{\chi^2_8}{8}}} \text{ and } Y_{\cdot,j} \overset{i.i.d.}{\sim} \frac{\Normal\left(\vec{0},\I_m\right)}{\sqrt{\frac{\chi^2_8}{8}}}  $$ 
\item The rows of $\X$ and $\Y$ are i.i.d. ARMA entries of parameters $\text{AR}=(0.6,0.2)$ and $\text{MA}=(0.5,0.2)$. Moreover, the traces of the matrices are standardised by the estimated variance.
\end{enumerate}

\begin{table}[hbt!]
\begin{center}
\begin{adjustbox}{addcode={\begin{minipage}{\width}}{\caption{Simulations of the extreme residual spikes. The values $(\hat{\mu},\hat{\sigma})$ and $(\mu,\sigma)$ are respectively the estimations of the mean and the standard error of the residual spikes obtained by respectively the Main Theorem \ref{TH=Main} and empirical methods using 500 replicates and $\theta=5000$. } \label{tab:Tablecomplete}
\end{minipage}},rotate=90,center}
\scalebox{0.55}{
\begin{tabular}{c c}
&$\lambda_{\max}\left(\hat{\hat{\Sigma}}_{X,P_k}^{-1/2} \hat{\hat{\Sigma}}_{Y,P_k} \hat{\hat{\Sigma}}_{X,P_k}^{-1/2}\right)$ \\ 
 \rotatebox{90}{\hspace{-1cm}\tiny 1. Normal entries. }&\tiny 
\begin{tabular}{c @{\,\vrule width 0.75mm\,}  c @{\,\vrule width 0.75mm\,} c c c c @{\,\vrule width 0.25mm\,}  c c c c @{\,\vrule width 0.25mm\,} c c c c @{\,\vrule width 0.25mm\,} c c c c }
\backslashbox{$n_Y$}{$n_X$} & & \multicolumn{4}{c}{100} & \multicolumn{4}{c}{500} & \multicolumn{4}{c}{1000} & \multicolumn{4}{c}{2000}\\
\specialrule{0.75mm}{0pt}{0pt}
 & \backslashbox{$k$}{$m$} & \multicolumn{2}{c}{100}  & \multicolumn{2}{c}{1000}  &  \multicolumn{2}{c}{100}  & \multicolumn{2}{c}{1000} & \multicolumn{2}{c}{100}  & \multicolumn{2}{c}{1000} & \multicolumn{2}{c}{100}  & \multicolumn{2}{c}{1000} \\
 \specialrule{0.75mm}{0pt}{0pt}
 & &  $(\hat{\mu},\hat{\sigma})$  & $(\mu,\sigma)$ &  $(\hat{\mu},\hat{\sigma})$  & $(\mu,\sigma)$ & $(\hat{\mu},\hat{\sigma})$  & $(\mu,\sigma)$ & $(\hat{\mu},\hat{\sigma})$  & $(\mu,\sigma)$ & $(\hat{\mu},\hat{\sigma})$  & $(\mu,\sigma)$ & $(\hat{\mu},\hat{\sigma})$  & $(\mu,\sigma)$ & $(\hat{\mu},\hat{\sigma})$  & $(\mu,\sigma)$ & $(\hat{\mu},\hat{\sigma})$  & $(\mu,\sigma)$ \\
 \specialrule{0.75mm}{0pt}{0pt}
100 
&    & $(3.69,0.47 )$ & $(3.77,0.48 )$ & $(21.91,2.38 )$ & $(22.32,2.36 )$ & $(2.88,0.25 )$ & $(2.94,0.25 )$ & $(13.93,0.71 )$ & $(14.15,0.69 )$ & $(2.76,0.21 )$ & $(2.78,0.21 )$ & $(13.01,0.60 )$ & $(13.14,0.57 )$ & $(2.73,0.22 )$ & $(2.75,0.23 )$ & $(12.39,0.52 )$ & $(12.59,0.50)$ \\ 
500
&  &  &  &  &  & $( 1.87,0.12 )$ & $(1.87,0.12 )$ & $(5.85,0.33 )$ & $(5.86,0.34 )$ & $(1.72,0.09 )$ & $(1.73,0.09 )$ & $(4.79,0.20 )$ & $(4.81,0.21 )$ & $(1.63,0.07 )$ & $(1.63,0.08 )$ & $(4.28,0.15 )$ & $(4.29,0.16)$ \\ 
1000
&  &  &  &  &  &  &  &  &  & $(1.57,0.07 )$ & $(1.57,0.07 )$ & $(3.73,0.15 )$ & $(3.74,0.15 )$ & $(1.48,0.06 )$ & $(1.47,0.06 )$ & $(3.19,0.10 )$ & $(3.19,0.11)$ \\ 
2000
& &  & &  &  &  &  &  &  &  &  &  &  & $(1.37,0.04 )$ & $(1.37,0.04 )$ & $(2.62,0.08 )$ & $(2.63,0.08)$ \\ 
 \hline  \hline  
\end{tabular} \\
& \\
\rotatebox{90}{\hspace{-1cm} \tiny 2. Multivariate Student. } & \tiny \begin{tabular}{c @{\,\vrule width 0.75mm\,}  c @{\,\vrule width 0.75mm\,} c c c c @{\,\vrule width 0.25mm\,}  c c c c @{\,\vrule width 0.25mm\,} c c c c @{\,\vrule width 0.25mm\,} c c c c }
\backslashbox{$n_Y$}{$n_X$} & & \multicolumn{4}{c}{100} & \multicolumn{4}{c}{500} & \multicolumn{4}{c}{1000} & \multicolumn{4}{c}{2000}\\
\specialrule{0.75mm}{0pt}{0pt}
 & \backslashbox{$k$}{$m$} & \multicolumn{2}{c}{100}  & \multicolumn{2}{c}{1000}  &  \multicolumn{2}{c}{100}  & \multicolumn{2}{c}{1000} & \multicolumn{2}{c}{100}  & \multicolumn{2}{c}{1000} & \multicolumn{2}{c}{100}  & \multicolumn{2}{c}{1000} \\
 \specialrule{0.75mm}{0pt}{0pt}
 & &  $(\hat{\mu},\hat{\sigma})$  & $(\mu,\sigma)$ &  $(\hat{\mu},\hat{\sigma})$  & $(\mu,\sigma)$ & $(\hat{\mu},\hat{\sigma})$  & $(\mu,\sigma)$ & $(\hat{\mu},\hat{\sigma})$  & $(\mu,\sigma)$ & $(\hat{\mu},\hat{\sigma})$  & $(\mu,\sigma)$ & $(\hat{\mu},\hat{\sigma})$  & $(\mu,\sigma)$ & $(\hat{\mu},\hat{\sigma})$  & $(\mu,\sigma)$ & $(\hat{\mu},\hat{\sigma})$  & $(\mu,\sigma)$ \\
 \specialrule{0.75mm}{0pt}{0pt}
100 
& & $(5.00,1.23)$ & $(5.37,1.24)$ & $(31.66,7.70)$ & $(33.68,7.67)$ & $(3.27,0.37)$ & $(3.35,0.36)$ & $(18.45,2.87)$ & $(19.15,2.96)$ & $(3.24,0.39)$ & $(3.33,0.40)$ & $(16.33,1.70)$ & $(16.69,1.52)$ & $(3.16,0.42)$ & $(3.27,0.40)$ & $(18.30,3.61)$ & $(18.72,3.24)$ \\ 
500
&  &  &  &  &  & $(2.16,0.17)$ & $(2.15,0.16)$ & $(7.56,0.62)$ & $(7.62,0.57)$ & $(1.92,0.13)$ & $(1.91,0.13)$ & $(6.12,0.42)$ & $(6.22,0.43)$ & $(1.81,0.11)$ & $(1.84,0.11)$ & $(6.15,0.90)$ & $(6.32,1.02)$ \\ 
1000 
&  &  &  &  &  &  &  &  &  & $(1.72,0.10)$ & $(1.71,0.09)$ & $(5.28,0.87)$ & $(5.34,0.84)$ & $(1.70,0.17)$ & $(1.70,0.16)$ & $(3.94,0.21)$ & $(3.92,0.21)$ \\ 
2000 
&  &  &  &  &  &  &  &  &  &  &  &  &  & $(1.46,0.06)$ & $(1.46,0.05)$ & $(3.17,0.13)$ & $(3.18,0.12)$ \\ 
  \hline \hline 
\end{tabular}\\
\rotatebox{90}{\hspace{-2.3cm} \tiny 3. ARMA \scriptsize $\big((0.6,0.2),(0.5,0.2)\big)$. }& \tiny

\begin{tabular}{c @{\,\vrule width 0.75mm\,}  c @{\,\vrule width 0.75mm\,} c c c c @{\,\vrule width 0.25mm\,}  c c c c @{\,\vrule width 0.25mm\,} c c c c @{\,\vrule width 0.25mm\,} c c c c }
\backslashbox{$n_Y$}{$n_X$} & & \multicolumn{4}{c}{100} & \multicolumn{4}{c}{500} & \multicolumn{4}{c}{1000} & \multicolumn{4}{c}{2000}\\
\specialrule{0.75mm}{0pt}{0pt}
 & \backslashbox{$k$}{$m$} & \multicolumn{2}{c}{100}  & \multicolumn{2}{c}{1000}  &  \multicolumn{2}{c}{100}  & \multicolumn{2}{c}{1000} & \multicolumn{2}{c}{100}  & \multicolumn{2}{c}{1000} & \multicolumn{2}{c}{100}  & \multicolumn{2}{c}{1000} \\
 \specialrule{0.75mm}{0pt}{0pt}
 & &  $(\hat{\mu},\hat{\sigma})$  & $(\mu,\sigma)$ &  $(\hat{\mu},\hat{\sigma})$  & $(\mu,\sigma)$ & $(\hat{\mu},\hat{\sigma})$  & $(\mu,\sigma)$ & $(\hat{\mu},\hat{\sigma})$  & $(\mu,\sigma)$ & $(\hat{\mu},\hat{\sigma})$  & $(\mu,\sigma)$ & $(\hat{\mu},\hat{\sigma})$  & $(\mu,\sigma)$ & $(\hat{\mu},\hat{\sigma})$  & $(\mu,\sigma)$ & $(\hat{\mu},\hat{\sigma})$  & $(\mu,\sigma)$ \\
 \specialrule{0.75mm}{0pt}{0pt}
100
&  & $(28.96,12.69)$ & $ (36.79,16.64)$ & $ (251.42,105.05)$ & $ (333.04,156.12)$ & $ (19.48,7.79)$ & $ (25.23,9.14)$ & $ (154.16,50.52)$ & $ (187.38,61.43)$ & $ (15.59,5.62)$ & $ (18.96,5.99)$ & $ (139.57,46.00)$ & $ (170.16,49.33)$ & $ (15.10,5.05)$ & $ (18.11,5.88)$ & $ (132.97,47.80)$ & $ (163.83,51.00)$ \\ 
500
&  &  &  &  &  & $(7.24,1.55)$ & $ (7.70,1.67)$ & $ (55.71,11.52)$ & $ (58.73,12.43)$ & $ (5.86,1.06)$ & $ (6.15,1.21)$ & $ (42.51,7.07)$ & $ (45.19,8.13)$ & $ (5.29,0.96)$ & $ (5.57,1.01)$ & $ (36.00,5.82)$ & $ (36.68,5.54)$ \\ 
1000 
&   &  &  &  &  &  &  &  &  & $(4.46,0.66)$ & $ (4.52,0.73)$ & $ (28.97,3.98)$ & $ (29.43,4.26)$ & $ (3.83,0.48)$ & $ (3.93,0.54)$ & $ (22.44,2.76)$ & $ (23.08,2.87)$ \\ 
2000 
&  &  &  &  &  &  &  &  &  &  &  &  &  &  $(3.09,0.37)$ &  $ (3.17,0.34)$ &  $ (15.35,1.61)$ &  $ (15.58,1.60)$ \\ 
\end{tabular}\\
& \\
& \\
&$\lambda_{\min}\left(\hat{\hat{\Sigma}}_{X,P_k}^{-1/2} \hat{\hat{\Sigma}}_{Y,P_k} \hat{\hat{\Sigma}}_{X,P_k}^{-1/2}\right)$ \\ 
 \rotatebox{90}{\hspace{-1cm} \tiny 1. Normal entries. }&\tiny 
\begin{tabular}{c @{\,\vrule width 0.75mm\,}  c @{\,\vrule width 0.75mm\,} c c c c @{\,\vrule width 0.25mm\,}  c c c c @{\,\vrule width 0.25mm\,} c c c c @{\,\vrule width 0.25mm\,} c c c c }
\backslashbox{$n_Y$}{$n_X$} & & \multicolumn{4}{c}{100} & \multicolumn{4}{c}{500} & \multicolumn{4}{c}{1000} & \multicolumn{4}{c}{2000}\\
\specialrule{0.75mm}{0pt}{0pt}
 & \backslashbox{$k$}{$m$} & \multicolumn{2}{c}{100}  & \multicolumn{2}{c}{1000}  &  \multicolumn{2}{c}{100}  & \multicolumn{2}{c}{1000} & \multicolumn{2}{c}{100}  & \multicolumn{2}{c}{1000} & \multicolumn{2}{c}{100}  & \multicolumn{2}{c}{1000} \\
 \specialrule{0.75mm}{0pt}{0pt}
 & &  $(\hat{\mu},\hat{\sigma})$  & $(\mu,\sigma)$ &  $(\hat{\mu},\hat{\sigma})$  & $(\mu,\sigma)$ & $(\hat{\mu},\hat{\sigma})$  & $(\mu,\sigma)$ & $(\hat{\mu},\hat{\sigma})$  & $(\mu,\sigma)$ & $(\hat{\mu},\hat{\sigma})$  & $(\mu,\sigma)$ & $(\hat{\mu},\hat{\sigma})$  & $(\mu,\sigma)$ & $(\hat{\mu},\hat{\sigma})$  & $(\mu,\sigma)$ & $(\hat{\mu},\hat{\sigma})$  & $(\mu,\sigma)$ \\
 \specialrule{0.75mm}{0pt}{0pt}
100&
 & $(0.269,0.035)$ & $(0.278,0.035)$ & $(0.045,0.005)$ & $(0.045,0.005)$ & $(0.347,0.031)$ & $(0.348,0.042)$ & $(0.072,0.003)$ & $(0.070,0.009)$ & $(0.363,0.028)$ & $(0.353,0.043)$ & $(0.077,0.003)$ & $(0.076,0.010)$ & $(0.368,0.031)$ & $(0.362,0.042)$ & $(0.080,0.003)$ & $(0.079,0.011)$ \\ 
500  &
   & & &  &  & $(0.535,0.034)$ & $(0.544,0.035)$ & $(0.171,0.009)$ & $(0.171,0.010)$ & $(0.579,0.030)$ & $(0.579,0.033)$ & $(0.208,0.009)$ & $(0.208,0.011)$ & $(0.612,0.028)$ & $(0.619,0.030)$ & $(0.234,0.008)$ & $(0.233,0.013)$ \\ 
1000  &
   &  &  &  &  &  &  &  &  & $(0.637,0.028)$ & $(0.641,0.028)$ & $(0.269,0.011)$ & $(0.268,0.010)$ & $(0.676,0.025)$ & $(0.682,0.026)$ & $(0.314,0.010)$ & $(0.313,0.012)$ \\ 
2000 
   &  &  &  &  &  &  &  &  &  &  &  &  & $(0.730,0.023)$ & $(0.729,0.021)$ & $(0.382,0.011)$ & $(0.383,0.011)$ \\ 
 \hline  \hline
\end{tabular}\\   
 
\rotatebox{90}{\hspace{-1cm} \tiny 2. Multivariate Student. } &\tiny \begin{tabular}{c @{\,\vrule width 0.75mm\,}  c @{\,\vrule width 0.75mm\,} c c c c @{\,\vrule width 0.25mm\,}  c c c c @{\,\vrule width 0.25mm\,} c c c c @{\,\vrule width 0.25mm\,} c c c c }
\backslashbox{$n_Y$}{$n_X$} & & \multicolumn{4}{c}{100} & \multicolumn{4}{c}{500} & \multicolumn{4}{c}{1000} & \multicolumn{4}{c}{2000}\\
\specialrule{0.75mm}{0pt}{0pt}
 & \backslashbox{$k$}{$m$} & \multicolumn{2}{c}{100}  & \multicolumn{2}{c}{1000}  &  \multicolumn{2}{c}{100}  & \multicolumn{2}{c}{1000} & \multicolumn{2}{c}{100}  & \multicolumn{2}{c}{1000} & \multicolumn{2}{c}{100}  & \multicolumn{2}{c}{1000} \\
 \specialrule{0.75mm}{0pt}{0pt}
 & &  $(\hat{\mu},\hat{\sigma})$  & $(\mu,\sigma)$ &  $(\hat{\mu},\hat{\sigma})$  & $(\mu,\sigma)$ & $(\hat{\mu},\hat{\sigma})$  & $(\mu,\sigma)$ & $(\hat{\mu},\hat{\sigma})$  & $(\mu,\sigma)$ & $(\hat{\mu},\hat{\sigma})$  & $(\mu,\sigma)$ & $(\hat{\mu},\hat{\sigma})$  & $(\mu,\sigma)$ & $(\hat{\mu},\hat{\sigma})$  & $(\mu,\sigma)$ & $(\hat{\mu},\hat{\sigma})$  & $(\mu,\sigma)$ \\
 \specialrule{0.75mm}{0pt}{0pt}
\multirow{3}{*}{100} 
&1  & $(0.203,0.047)$ & $(0.217,0.034)$ & $(0.031,0.007)$ & $(0.035,0.005)$ & $(0.309,0.035)$ & $(0.311,0.042)$ & $(0.054,0.008)$ & $(0.058,0.009)$ & $(0.307,0.038)$ & $(0.308,0.042)$ & $(0.062,0.006)$ & $(0.056,0.008)$ & $(0.318,0.040)$ & $(0.299,0.040)$ & $(0.055,0.011)$ & $(0.054,0.009)$ \\ 
500
&  &   &   &   &   & $(0.465,0.036)$ & $(0.440,0.039)$ & $(0.132,0.011)$ & $(0.132,0.011)$ & $(0.522,0.034)$ & $(0.535,0.039)$ & $(0.163,0.012)$ & $(0.147,0.011)$ & $(0.550,0.032)$ & $(0.515,0.032)$ & $(0.161,0.025)$ & $(0.160,0.017)$ \\ 
1000 
&  &   &   &   &   &   &   &   &   & $(0.583,0.033)$ & $(0.562,0.034)$ & $(0.190,0.031)$ & $(0.196,0.022)$ & $(0.587,0.059)$ & $(0.628,0.036)$ & $(0.255,0.014)$ & $(0.249,0.013)$\\ 

2000 
&  &   &   &   &   &   &   &   &   &   &   &   &   & $(0.684,0.026)$ & $(0.676,0.027)$ & $(0.315,0.013)$ & $(0.310,0.013)$ \\ 
 \hline \hline
\end{tabular}\\  
\rotatebox{90}{\hspace{-2.3cm} \tiny 3. ARMA \scriptsize$\big((0.6,0.2),(0.5,0.2)\big)$. }& \tiny

\begin{tabular}{c @{\,\vrule width 0.75mm\,}  c @{\,\vrule width 0.75mm\,} c c c c @{\,\vrule width 0.25mm\,}  c c c c @{\,\vrule width 0.25mm\,} c c c c @{\,\vrule width 0.25mm\,} c c c c }
\backslashbox{$n_Y$}{$n_X$} & & \multicolumn{4}{c}{100} & \multicolumn{4}{c}{500} & \multicolumn{4}{c}{1000} & \multicolumn{4}{c}{2000}\\
\specialrule{0.75mm}{0pt}{0pt}
 & \backslashbox{$k$}{$m$} & \multicolumn{2}{c}{100}  & \multicolumn{2}{c}{1000}  &  \multicolumn{2}{c}{100}  & \multicolumn{2}{c}{1000} & \multicolumn{2}{c}{100}  & \multicolumn{2}{c}{1000} & \multicolumn{2}{c}{100}  & \multicolumn{2}{c}{1000} \\
 \specialrule{0.75mm}{0pt}{0pt}
 & &  $(\hat{\mu},\hat{\sigma})$  & $(\mu,\sigma)$ &  $(\hat{\mu},\hat{\sigma})$  & $(\mu,\sigma)$ & $(\hat{\mu},\hat{\sigma})$  & $(\mu,\sigma)$ & $(\hat{\mu},\hat{\sigma})$  & $(\mu,\sigma)$ & $(\hat{\mu},\hat{\sigma})$  & $(\mu,\sigma)$ & $(\hat{\mu},\hat{\sigma})$  & $(\mu,\sigma)$ & $(\hat{\mu},\hat{\sigma})$  & $(\mu,\sigma)$ & $(\hat{\mu},\hat{\sigma})$  & $(\mu,\sigma)$ \\
 \specialrule{0.75mm}{0pt}{0pt}
100 
&  & $(0.035,0.015)$ & $(0.033,0.017)$ & $(0.004,0.002)$ & $(0.004,0.002)$ & $(0.052,0.020)$ & $(0.041,0.017)$ & $(0.007,0.002)$ & $(0.005,0.002)$ & $(0.065,0.025)$ & $(0.054,0.020)$ & $(0.007,0.002)$ & $(0.006,0.002)$ & $(0.067,0.022)$ & $(0.057,0.021)$ & $(0.007,0.003)$ & $(0.006,0.002)$ \\ 
500
&  & & & & & $(0.137,0.030)$ & $(0.136,0.027)$ & $(0.018,0.004)$ & $(0.018,0.004)$ & $(0.170,0.032)$ & $(0.164,0.032)$ & $(0.023,0.004)$ & $(0.023,0.004)$ & $(0.189,0.035)$ & $(0.184,0.034)$ & $(0.028,0.005)$ & $(0.028,0.005)$ \\ 
1000
&  & & & & & & & & & $(0.225,0.034)$ & $(0.223,0.032)$ & $(0.035,0.005)$ & $(0.035,0.005)$ & $(0.265,0.035)$ & $(0.258,0.033)$ & $(0.045,0.006)$ & $(0.044,0.006)$ \\ 
2000
&  & & & & & & & & & & & & & $(0.322,0.036)$ & $(0.325,0.037)$ & $(0.065,0.007)$ & $(0.065,0.007)$ 
\end{tabular}\\
\normalsize \\ \\ \\ 
\end{tabular}}

\end{adjustbox}
\end{center}
\end{table}

The Table \ref{tab:Tablecomplete} compares the estimations of the mean and the standard error of the residual spikes $(\hat{\mu},\hat{\sigma})$ to their empirical values $(\mu,\sigma)$.\\
The simulations are computed for the three scenarios described above. The perturbation $P=\I_m+ (\theta-1) u u_i$ is without loss of generality assumed canonical and the eigenvalue $\theta$ is fixed to 5000.

\pagebreak

\begin{supplement}
\sname{Supplement A}\label{suppA}
\stitle{Statistical applications of Random matrix theory: comparison of two populations I, 
Supplement}
\slink[url]{The supplement is in the second part of this paper.}
\sdescription{Proofs of Theorems}
\end{supplement}

\bibliographystyle{imsart-nameyear}
\bibliography{biblio4} 

\begin{thebibliography}{11}

\bibitem[\protect\citeauthoryear{Anderson}{1958}]{multi22}
\begin{bbook}[author]
\bauthor{\bsnm{Anderson},~\bfnm{T.~W.}\binits{T.~W.}}
(\byear{1958}).
\btitle{An introduction to Multivariate Statistical Analysis}.
\bseries{Wiley publications in statistics}.
\bpublisher{Wiley}.
\end{bbook}
\endbibitem

\bibitem[\protect\citeauthoryear{Anderson}{2003}]{multi2}
\begin{bbook}[author]
\bauthor{\bsnm{Anderson},~\bfnm{T.~W.}\binits{T.~W.}}
(\byear{2003}).
\btitle{An introduction to Multivariate Statistical Analysis}.
\bseries{Wiley Series in Probability and Statistics}.
\bpublisher{Wiley}.
\end{bbook}
\endbibitem

\bibitem[\protect\citeauthoryear{Anderson, Guionnet and Zeitouni}{2009}]{Alice}
\begin{bbook}[author]
\bauthor{\bsnm{Anderson},~\bfnm{Greg~W.}\binits{G.~W.}},
  \bauthor{\bsnm{Guionnet},~\bfnm{Alice}\binits{A.}} \AND
  \bauthor{\bsnm{Zeitouni},~\bfnm{Ofer}\binits{O.}}
(\byear{2009}).
\btitle{An Introduction to Random Matrices}.
\bseries{Cambridge Studies in Advanced Mathematics}.
\bpublisher{Cambridge University Press}.
\end{bbook}
\endbibitem

\bibitem[\protect\citeauthoryear{Bai and Silverstein}{2010}]{Appliedbook}
\begin{bbook}[author]
\bauthor{\bsnm{Bai},~\bfnm{Zhidong}\binits{Z.}} \AND
  \bauthor{\bsnm{Silverstein},~\bfnm{Jack~W.}\binits{J.~W.}}
(\byear{2010}).
\btitle{Spectral Analysis of Large Dimensional Random Matrices}.
\bpublisher{Springer}.
\end{bbook}
\endbibitem

\bibitem[\protect\citeauthoryear{Benaych-Georges and Rao}{2009}]{deformedRMT}
\begin{barticle}[author]
\bauthor{\bsnm{Benaych-Georges},~\bfnm{Florent}\binits{F.}} \AND
  \bauthor{\bsnm{Rao},~\bfnm{Nadakuditi~Raj}\binits{N.~R.}}
(\byear{2009}).
\btitle{The eigenvalues and eigenvectors of finite, low rank perturbations of
  large random matrices}.
\bjournal{Advances in Mathematics}
\bvolume{227}
\bpages{494-521}.
\end{barticle}
\endbibitem

\bibitem[\protect\citeauthoryear{Bose}{2018}]{bookrecent}
\begin{bbook}[author]
\bauthor{\bsnm{Bose},~\bfnm{Arup}\binits{A.}}
(\byear{2018}).
\btitle{Patterned Random matrices}.
\bpublisher{Chapman and Hall/CRC}.
\end{bbook}
\endbibitem

\bibitem[\protect\citeauthoryear{Marchenko and Pastur}{1967}]{deformed}
\begin{barticle}[author]
\bauthor{\bsnm{Marchenko},~\bfnm{V.~A.}\binits{V.~A.}} \AND
  \bauthor{\bsnm{Pastur},~\bfnm{L.~A.}\binits{L.~A.}}
(\byear{1967}).
\btitle{Distribution of eigenvalues for some sets of random matrices}.
\bjournal{Math. USSR}
\bvolume{1}
\bpages{457-483}.
\end{barticle}
\endbibitem

\bibitem[\protect\citeauthoryear{Mardia, Kent and Bibby}{1979}]{multi3}
\begin{bbook}[author]
\bauthor{\bsnm{Mardia},~\bfnm{K.~V.}\binits{K.~V.}},
  \bauthor{\bsnm{Kent},~\bfnm{J.~T.}\binits{J.~T.}} \AND
  \bauthor{\bsnm{Bibby},~\bfnm{J.~M.}\binits{J.~M.}}
(\byear{1979}).
\btitle{Multivariate Analysis}.
\bseries{Probability and mathematical statistics}.
\bpublisher{Academic press}.
\end{bbook}
\endbibitem

\bibitem[\protect\citeauthoryear{Mari\'etan and
  Morgenthaler}{2019}]{Suppmaterial}
\begin{barticle}[author]
\bauthor{\bsnm{Mari\'etan},~\bfnm{R\'emy}\binits{R.}} \AND
  \bauthor{\bsnm{Morgenthaler},~\bfnm{Stephan}\binits{S.}}
(\byear{2019}).
\btitle{Statistical applications of Random matrix theory: comparison of two
  populations I, Supplement}.
\end{barticle}
\endbibitem

\bibitem[\protect\citeauthoryear{Muirhead}{2005}]{multi}
\begin{bbook}[author]
\bauthor{\bsnm{Muirhead},~\bfnm{Robb~J.}\binits{R.~J.}}
(\byear{2005}).
\btitle{Aspect of Multivariate Statistical Theory}.
\bseries{Wiley Series in Probability and Statistics}.
\bpublisher{Wiley-Interscience}.
\end{bbook}
\endbibitem

\bibitem[\protect\citeauthoryear{Tao}{2012}]{Tao}
\begin{bmisc}[author]
\bauthor{\bsnm{Tao},~\bfnm{Terence}\binits{T.}}
(\byear{2012}).
\btitle{Topics in random matrix theory}.
\bhowpublished{\url{http://www.math.hkbu.edu.hk/~ttang/UsefulCollections/matrix-book-2011-08.pdf}}.
\end{bmisc}
\endbibitem

\end{thebibliography}
\addcontentsline{toc}{section}{Bibliography}

\pagebreak



\section*{Supplementary material}

# Statistical applications of Random matrix theory: comparison of two populations I, Supplement [*]

**Rémy Mariétan[†] and Stephan Morgenthaler[‡]**

*Department of Mathematics*
*École Polytechnique Fédéral de Lausanne*
*1015 Lausanne*
*e-mail:* remy.marietan@alumni.epfl.ch, stephan.morgenthaler@epfl.ch

**Abstract:** This supplement presents the proofs of the theorems and lemmas of Mariétan and Morgenthaler (2020).

**Keywords and phrases:** High dimension, equality test of two covariance matrices, Random matrix theory, residual spike, spike model, dependent data, eigenvector, eigenvalue.

## 1. Notation, definitions and assumptions

As presented in Mariétan and Morgenthaler (2020) we use the following notation.

**Notation 1.1.**

- If $W$ is a symmetric random matrix, we denote by $\left(\hat{\lambda}_{W,i}, \hat{u}_{W,i}\right)$ its $i^{\text{th}}$ eigenvalues and eigenvectors, ordered by decreasing eigenvalue.
- A finite perturbation of order $k$ is denoted by $P_k = I_m + \sum_{i=1}^{k}(\theta_i - 1)u_i u_i^t \in \mathbb{R}^{m \times m}$ with $u_1, u_2, ..., u_k \in \mathbb{R}^{m \times m}$ orthonormal vectors.
- We denote by $W \in \mathbb{R}^{m \times m}$ an random matrix invariant by rotation as defined in Assumption 1.1. Moreover, the estimated covariance matrix is decomposed as $\hat{\Sigma} = P^{1/2} W P^{1/2}$. When comparing two groups, we use $W_X$, $W_Y$ and $\hat{\Sigma}_X$, $\hat{\Sigma}_Y$.
- Suppose two groups $X$ and $Y$. The perturbation of order 1 of the matrices $W_X$ and $W_Y$ are $\hat{\Sigma}_X = P^{1/2} W_X P^{1/2}$ and $\hat{\Sigma}_Y = P^{1/2} W_Y P^{1/2}$ respectively. Then, we define for $\hat{\Sigma}_X$ (and similarly for $\hat{\Sigma}_Y$):
    - $\hat{u}_X = \hat{u}_{\hat{\Sigma}_X}$, its first eigenvector.
    - $\hat{u}_{X,j}$, the $j^{\text{th}}$ component of the first eigenvector.
    - $\hat{\theta}_X = \hat{\lambda}_{\hat{\Sigma}_X, 1}$, the largest eigenvalue of $\hat{\Sigma}_X$.

---

- $\hat{\lambda}_{\hat{\Sigma}_X,i}$, its $i^{\text{th}}$ eigenvalue, keeping in mind that $\hat{\lambda}_{\hat{\Sigma}_X,1} = \hat{\theta}_X$.

When the results concern only one group, we use simpler notation, $\hat{u}$, $\hat{\theta}$ and $\hat{\lambda}_{\hat{\Sigma},i}$.

- We define the function $M_{s_1,s_2,X}(\rho_X)$, $M_{s_1,s_2,Y}(\rho_Y)$ and $M_{s_1,s_2}(\rho_X,\rho_Y)$ as

$$M_{s_1,s_2,X}(\rho_X) = \frac{1}{m}\sum_{i=1}^{m}\frac{\hat{\lambda}_{W_X,i}^{s_1}}{\left(\rho_X - \hat{\lambda}_{W_X,i}\right)^{s_2}}, \text{ and analogously for } Y$$

$$M_{s_1,s_2}(\rho_X,\rho_Y) = \frac{M_{s_1,s_2,X}(\rho_X) + M_{s_1,s_2,Y}(\rho_Y)}{2}.$$

When $s_2 = 0$, we use $M_{s_1,X} = M_{s_1,0,X}$. When we only study one group, we use the simpler notation $M_{s_1,s_2}(\rho)$ when no confusion is possible.

- We use two transforms inspired by the T-transform:
  - $T_{W,u}(z) = \sum_{i=1}^{m}\frac{\hat{\lambda}_{W,i}}{z-\hat{\lambda}_{W,i}}\langle \hat{u}_{W,i}, u\rangle^2$ is the T-transform in direction $u$ using the random matrix $W$.
  - $\hat{T}_{\hat{\Sigma}_X}(z) = \frac{1}{m}\sum_{i=k+1}^{m}\frac{\hat{\lambda}_{\hat{\Sigma}_X,i}}{z-\hat{\lambda}_{\hat{\Sigma}_X,i}}$, and $\hat{T}_{W_X}(z) = \frac{1}{m}\sum_{i=1}^{m}\frac{\hat{\lambda}_{W_X,i}}{z-\hat{\lambda}_{W_X,i}}$ are the estimated T-transforms using $\hat{\Sigma}_X$ and $W$ respectively.

The main assumption of Mariétan and Morgenthaler (2020) concerns the estimated covariance matrices.

**Assumption 1.1.** Let $W_X$ and $W_Y$ be such that

$$W_X = O_X \Lambda_X O_X \text{ and } W_Y = O_Y \Lambda_Y O_Y,$$

where

$O_X, O_Y$ are unit orthonormal invariant and independent random matrices,

$\Lambda_X, \Lambda_Y$ are diagonal bounded matrices and independent of $O_X, O_Y$,

Trace $(W_X) = m$ and Trace $(W_Y) = m$.

Assume $P_X = I_m + (\theta_X - 1)u_X u_X^t$ and $P_Y = I_m + (\theta_Y - 1)u_Y u_Y^t$. Then

$$\hat{\Sigma}_X = P_X^{1/2} W_X P_X^{1/2} \text{ and } \hat{\Sigma}_Y = P_Y^{1/2} W_Y P_Y^{1/2}.$$

We next define the filtered covariance matrix and the residual spike.

**Definition 1.1.** Suppose $\hat{\Sigma}$ satisfies Assumption 1.1.
The **unbiased estimator of** $\theta$ is

$$\hat{\hat{\theta}} = 1 + \frac{1}{\frac{1}{m-1}\sum_{i=2}^{m}\frac{\hat{\lambda}_{\hat{\Sigma},i}}{\hat{\theta}-\hat{\lambda}_{\hat{\Sigma},i}}}.$$

The **filtered estimated covariance matrix** is defined as

$$\hat{\hat{\Sigma}} = I_m + (\hat{\hat{\theta}} - 1)\hat{u}\hat{u}^t.$$

**Definition 1.2.** The **residual spikes** are the largest and the smallest eigenvalues of
$$\hat{\hat{\Sigma}}_X^{-1/2} \hat{\hat{\Sigma}}_Y \hat{\hat{\Sigma}}_X^{-1/2}.$$

Some results assume a large perturbation, where large has one of the following meanings:

**Assumption 1.2.**
(A1) $\frac{\theta}{\sqrt{n}} \to \infty$.
(A2) $\theta \to \infty$.

Finally we define the **detectability** of perturbations.

**Definition 1.3.** We assume that a perturbation $P = I_m + (\theta - 1)uu^t$ is **detectable** via $\hat{\Sigma} = P^{1/2}WP^{1/2}$ if the perturbation creates a largest isolated eigenvalue, $\hat{\theta}$.

## 2. Theorems and proofs

Mariétan and Morgenthaler (2020) contains one main theorem, three theorems and one lemma.

### 2.1. Main Theorem

**Theorem 2.1.** *Suppose $W_X, W_Y \in \mathbb{R}^{m \times m}$ satisfy Assumption 1.1 and $\theta$ satisfies Assumption 1.2 (A1).*
*Let $P = I_m + (\theta - 1)e_1 e_1^t \in \mathbb{R}^{m \times m}$ and define*
$$\hat{\Sigma}_X = P^{1/2} W_X P^{1/2} \text{ and } \hat{\Sigma}_Y = P^{1/2} W_Y P^{1/2}.$$

*The induced filtered estimators become*
$$\hat{\hat{\Sigma}}_X = I_m + (\hat{\hat{\theta}}_{\hat{\Sigma}_X} - 1)\hat{u}_{\hat{\Sigma}_X} \hat{u}_{\hat{\Sigma}_X}^t \text{ and analogously for } Y$$

*where $\hat{\hat{\theta}}_{\hat{\Sigma}_X}$ and $\hat{\hat{\theta}}_{\hat{\Sigma}_Y}$ are the unbiased estimators of the largest eigenvalues of $\hat{\Sigma}_X$ and $\hat{\Sigma}_Y$ respectively, as defined in Definition 1.1.*

*Then, conditional on the spectra $S_{W_X} = \left\{\hat{\lambda}_{W_X,1}, \hat{\lambda}_{W_X,2}, ..., \hat{\lambda}_{W_X,m}\right\}$ and $S_{W_Y} = \left\{\hat{\lambda}_{W_Y,1}, \hat{\lambda}_{W_Y,2}, ..., \hat{\lambda}_{W_Y,m}\right\}$ of $W_X$ and $W_Y$, we have*

$$\sqrt{m} \frac{\left(\lambda_{\max}\left(\hat{\hat{\Sigma}}_X^{-1/2} \hat{\hat{\Sigma}}_Y \hat{\hat{\Sigma}}_X^{-1/2}\right) - \lambda^+\right)}{\sigma^+} \Bigg| S_{W_X}, S_{W_Y} \sim \mathbf{N}(0,1) + o_p(1),$$

*where $o_p(1)$ is with regard to large values of $m$,*
$$\lambda^+ = \sqrt{M_2^2 - 1} + M_2,$$

$$\sigma^{+2} = \frac{1}{(M_{2,X} + M_{2,Y} - 2)(M_{2,X} + M_{2,Y} + 2)}$$

$$\left( 9M_{2,X}^4 M_{2,Y} + 4M_{2,X}^3 M_{2,Y}^2 + 4M_{2,X}^3 M_{2,Y} + 2M_{2,X}^3 M_{3,Y} - 2M_{2,X}^2 M_{2,Y}^3 \right.$$

$$+4M_{2,X}^2 M_{2,Y}^2 - 11M_{2,X}^2 M_{2,Y} - 8M_{3,X} M_{2,X}^2 M_{2,Y} + 2M_{2,X}^2 M_{2,Y} M_{3,Y}$$

$$-2M_{2,X}^2 M_{3,Y} + M_{2,X}^2 M_{4,Y} + 4M_{2,X} M_{2,Y}^3 + M_{2,X} M_{2,Y}^2 + 4M_{2,X} M_{2,Y}$$

$$-4M_{3,X} M_{2,X} M_{2,Y}^2 - 4M_{3,X} M_{2,X} M_{2,Y} - 2M_{2,X} M_{2,Y}^2 M_{3,Y} - 4M_{2,X} M_{2,Y} M_{3,Y}$$

$$-6M_{2,X} M_{3,Y} + 2M_{4,X} M_{2,X} M_{2,Y} + 2M_{2,X} M_{2,Y} M_{4,Y} - 2M_{3,X} M_{2,Y}^2$$

$$+2M_{3,X} M_{2,Y} + M_{4,X} M_{2,Y}^2 + 4M_{2,X}^5 + 2M_{2,X}^4 - 4M_{3,X} M_{2,X}^3 - 13M_{2,X}^3$$

$$-2M_{3,X} M_{2,X}^2 + M_{4,X} M_{2,X}^2 - 2M_{2,X}^2 + 10M_{3,X} M_{2,X} + 4M_{2,X} + 4M_{3,X}$$

$$-2M_{4,X} + M_{2,Y}^5 + 2M_{2,Y}^4 - M_{2,Y}^3 - 2M_{2,Y}^2 + 4M_{2,Y} - 2M_{2,Y}^3 M_{3,Y}$$

$$\left. -2M_{2,Y}^2 M_{3,Y} + 2M_{2,Y} M_{3,Y} + 4M_{3,Y} + M_{2,Y}^2 M_{4,Y} - 2M_{4,Y} - 4 \right)$$

$$+ \frac{1}{\sqrt{(M_{2,X} + M_{2,Y} - 2)(M_{2,X} + M_{2,Y} + 2)}}$$

$$\left( 5M_{2,X}^3 M_{2,Y} - M_{2,X}^2 M_{2,Y}^2 + 2M_{2,X}^2 M_{2,Y} + 2M_{2,X}^2 M_{3,Y} - M_{2,X} M_{2,Y}^3 \right.$$

$$+2M_{2,X} M_{2,Y}^2 - 4M_{2,X} M_{2,Y} - 4M_{3,X} M_{2,X} M_{2,Y} - 2M_{2,X} M_{3,Y} + M_{2,X} M_{4,Y}$$

$$-2M_{3,X} M_{2,Y} + M_{4,X} M_{2,Y} + 4M_{2,X}^4 + 2M_{2,X}^3 - 4M_{3,X} M_{2,X}^2 - 5M_{2,X}^2$$

$$-2M_{3,X} M_{2,X} + M_{4,X} M_{2,X} + 2M_{2,X} + 2M_{3,X} + M_{2,Y}^4 + 2M_{2,Y}^3 + M_{2,Y}^2$$

$$\left. +2M_{2,Y} - 2M_{2,Y}^2 M_{3,Y} - 2M_{2,Y} M_{3,Y} - 2M_{3,Y} + M_{2,Y} M_{4,Y} \right),$$

Moreover,

$$\sqrt{m} \frac{\left( \lambda_{\min}\left( \hat{\hat{\Sigma}}_X^{-1/2} \hat{\hat{\Sigma}}_Y \hat{\hat{\Sigma}}_X^{-1/2} \right) - \lambda^- \right)}{\sigma^-} \bigg| S_{W_X}, S_{W_Y} \sim \mathbf{N}(0,1) + o_p(1),$$

where

$$\lambda^- = -\sqrt{M_2^2 - 1} + M_2,$$
$$\sigma^{-2} = \left( \lambda^- \right)^4 \sigma^{+2}.$$

*Remark* 2.1.

1. If the spectra have the Marcenko-Pastur distribution, $c_X = n_X/m$ and $c_Y = n_Y/m$, then:

$$c = \frac{c_X + c_Y}{2},$$
$$\lambda^+ = c + \sqrt{c(c+2)} + 1,$$
$$\sigma^{+2} = c_X^3 + c_X^2 c_Y + 3c_X^2 + 4c_X c_Y - c_X + c_Y^2 + c_Y$$
$$+ \frac{(8c_X + 2c_X^2 + (c_X^3 + 5c_X^2 + c_X^2 c_Y + 4c_X c_Y + 5c_X + 3c_Y + c_Y^2)\sqrt{c(c+2)}}{c+2}.$$

2. If $c_X$ tends to 0, then

$$\sigma^{+2} = \left( M_{2,Y}^5 + 2M_{2,Y}^4 - 2M_{3,Y}M_{2,Y}^3 + M_{2,Y}^3 - 4M_{3,Y}M_{2,Y}^2 + M_{4,Y}M_{2,Y}^2 + 2M_{2,Y}^2 \right.$$
$$\left. + 2M_{4,Y}M_{2,Y} + 2M_{2,Y} - 2M_{3,Y} - M_{4,Y} - 2 \right) \Big/ \left( (M_{2,Y}-1)(M_{2,Y}+3) \right)$$
$$+ \left( M_{2,Y}^4 + M_{2,Y}^3 - 2M_{3,Y}M_{2,Y}^2 + 2M_{2,Y}^2 - 2M_{3,Y}M_{2,Y} + M_{4,Y}M_{2,Y} \right.$$
$$\left. - 2M_{3,Y} + M_{4,Y} \right) \Big/ \sqrt{(M_{2,Y}-1)(M_{2,Y}+3)}.$$

*Proof.* Theorem 2.1
The Theorem is a direct consequence of Theorem 2.4 and Lemma 2.1 together with an application of Slutsky's Theorem. □

### 2.2. Unit invariant vector statistic

**Theorem 2.2.**
Let $W$ be a random matrix of trace $m$ with spectrum $S_W = \left\{ \hat{\lambda}_{W,1}, \hat{\lambda}_{W,2}, ..., \hat{\lambda}_{W,m} \right\}$. We denote by $u_{p_1}$ and $u_{p_2}$, two orthonormal invariant random vectors of size $m$, independent of the eigenvalues of $W$. Let

$$\vec{B}_m(\rho, \vec{s}, \vec{r}, \vec{p}) = \sqrt{m} \left( \begin{pmatrix} \sum_{i=1}^m \frac{\hat{\lambda}_{W,i}^{s_1}}{(\rho-\hat{\lambda}_{W,i})^{s_2}} u_{p_1,i} u_{p_2,i} \\ \sum_{i=1}^m \frac{\hat{\lambda}_{W,i}^{r_1}}{(\rho-\hat{\lambda}_{W,i})^{r_2}} u_{p_1,i} u_{p_2,i} \end{pmatrix} - \begin{pmatrix} M_{s_1,s_2} \\ M_{r_1,r_2} \end{pmatrix} \mathbf{1}_{p_1=p_2} \right),$$

where $\vec{s} = (s_1, s_2)$, $\vec{r} = (r_1, r_2)$ and $\vec{p} = (p_1, p_2)$ with indices $1 \leqslant p_1 \leqslant p_2 \leqslant m$ and $s_1, s_2, r_1, r_2 \in \mathbf{N}$.
If $p = p_1 = p_2$,

$$\vec{B}_m(\rho, \vec{s}, \vec{r}, \vec{p}) \Big| S_W \sim \mathbf{N}\left( \vec{0}, \begin{pmatrix} 2\left(M_{2s_1,2s_2} - M_{s_1,s_1}^2\right) & 2\left(M_{s_1+r_1,s_2+r_2} - M_{s_1,s_2}M_{r_1,r_2}\right) \\ 2\left(M_{s_1+r_1,s_2+r_2} - M_{s_1,s_2}M_{r_1,r_2}\right) & 2\left(M_{2r_1,2r_2} - M_{r_1,r_1}^2\right) \end{pmatrix} \right) + o_p(1),$$

where $M_{s,r} = M_{s,r}(\rho) = \frac{1}{m} \sum_{i=1}^m \frac{\hat{\lambda}_{W,i}^s}{(\rho-\hat{\lambda}_{W,i})^r}$.
If $p_1 \neq p_2$,

$$\vec{B}_m(\rho, \vec{s}, \vec{r}, \vec{p}) \Big| S_W \sim \mathbf{N}\left( \vec{0}, \begin{pmatrix} M_{2s_1,2s_2} - M_{s_1,s_1}^2 & M_{s_1+r_1,s_2+r_2} - M_{s_1,s_2}M_{r_1,r_2} \\ M_{s_1+r_1,s_2+r_2} - M_{s_1,s_2}M_{r_1,r_2} & M_{2r_1,2r_2} - M_{r_1,r_1}^2 \end{pmatrix} \right) + o_p(1),$$

In particular, with the notation $M_{s,0} = M_s = \frac{1}{m} \sum_{i=1}^m \hat{\lambda}_{W,i}^s$,

$$\sqrt{m} \left( \begin{pmatrix} \sum_{i=1}^m \hat{\lambda}_{W,i} u_{p,i}^2 \\ \sum_{i=1}^m \hat{\lambda}_{W,i}^2 u_{p,i}^2 \end{pmatrix} - \begin{pmatrix} 1 \\ M_2 \end{pmatrix} \right) \Big| S_W \sim \mathbf{N}\left( \vec{0}, \begin{pmatrix} 2(M_2-1) & 2(M_3-M_2) \\ 2(M_3-M_2) & 2(M_4-M_2^2) \end{pmatrix} \right) + o_p(1),$$

and

$$\sqrt{m} \left( \begin{pmatrix} \sum_{i=1}^m \hat{\lambda}_{W,i} u_{p_1,i} u_{p_2,i} \\ \sum_{i=1}^m \hat{\lambda}_{W,i}^2 u_{p_1,i} u_{p_2,i} \end{pmatrix} - \begin{pmatrix} 0 \\ 0 \end{pmatrix} \right) \Big| S_W \sim \mathbf{N}\left( \vec{0}, \begin{pmatrix} M_2-1 & M_3-M_2 \\ M_3-M_2 & M_4-M_2^2 \end{pmatrix} \right) + o_p(1).$$

*Finally if we look at $K$ bivariate normal random variables:*

$$\mathbf{B}_m\left(\vec{\rho},\mathbf{s},\mathbf{r},\mathbf{p}\right) = \left(\vec{B}_m\left(\rho_1,\vec{s}_1,\vec{r}_1,\vec{p}_1\right), \vec{B}_m\left(\rho_2,\vec{s}_2,\vec{r}_2,\vec{p}_2\right),...,\vec{B}_m\left(\rho_K,\vec{s}_K,\vec{r}_K,\vec{p}_K\right)\right),$$

*where $\vec{p}_i \neq \vec{p}_j$ if $i \neq j$. Then, conditional on the spectrum $S_W$, $\tilde{\mathbf{B}}_m(\vec{\rho},\mathbf{s},\mathbf{r},\mathbf{p})$ tends to a multivariate Normal with block-diagonal covariance matrix of blocks of size $2 \times 2$.*

*Remark* 2.2.1.

1. We assume the trace of $W$ is $m$, which can be assured by rescaling the matrix.
2. The condition of independence between eigenvectors and eigenvalues of $W$ is strong, but follows automatically if the eigenvectors are Haar distributed.
3. If $W$ is a rescaled standard Wishart, then

$$\sqrt{m}\left(\begin{pmatrix}\sum_{i=1}^m \hat{\lambda}_{W,i}\hat{u}_i^2 \\ \sum_{i=1}^m \hat{\lambda}_{W,i}^2\hat{u}_i^2\end{pmatrix} - \begin{pmatrix}1 \\ 1+c\end{pmatrix}\right) \xrightarrow[m\to\infty]{} \mathbf{N}\left(\vec{0}, \begin{pmatrix}2c & 2c(2+c) \\ 2c(2+c) & 2c(c+1)(c+4)\end{pmatrix}\right),$$

where $\hat{u}_1, \hat{u}_2, ..., \hat{u}_m$ are the eigenvectors of $W$.

*Proof.* **Theorem 2.2**

The proof is divided into three steps. First, we compute the first two moments of the statistics. Then, we show the asymptotic joint normality. Finally, we show the asymptotic independence of the statistics. In each step of the proof, we condition on the spectrum of $W$.

For any non-random real-valued function $g$,

$$\mathrm{E}\left[\sum_{i=1}^m g(\hat{\lambda}_{W,i})u_{1,i}^2\right] = \frac{1}{m}\sum_{i=1}^m g(\hat{\lambda}_{W,i}),$$

$$\mathrm{E}\left[\sum_{i=1}^m g(\hat{\lambda}_{W,i})u_{1,i}u_{2,i}\right] = 0$$

This proves the formulas for the first moments.

Since for unit invariant eigenvectors $\mathrm{Cov}\left(u_{1,i}^2, u_{1,j}^2\right) = -2/(m^2(m-1))$, we have

$$\mathrm{Var}\left(\sum_{i=1}^m \frac{\hat{\lambda}_{W,i}^{s_1}}{\left(\rho - \hat{\lambda}_{W,i}\right)^{s_2}} u_{1,i}^2\right)$$

$$= \sum_{i=1}^m \frac{\hat{\lambda}_{W,i}^{2s_1}}{\left(\rho - \hat{\lambda}_{W,i}\right)^{2s_2}} \frac{2}{m^2} - \sum_{\substack{i,j=1 \\ i\neq j}}^m \left(\frac{\hat{\lambda}_{W,i}^{s_1}}{\left(\rho - \hat{\lambda}_{W,i}\right)^{s_2}} \frac{\hat{\lambda}_{W,j}^{s_1}}{\left(\rho - \hat{\lambda}_{W,j}\right)^{s_2}}\right) \frac{2}{m^2(m-1)}$$

$$= \frac{2M_{2s_1,2s_2}(\rho)}{m} - \frac{2M_{s_1,s_2}(\rho)^2}{m-1} + \frac{2M_{2s_1,2s_2}(\rho)}{m(m-1)}$$

$$= \frac{2}{m}\left(M_{2s_1,2s_2}(\rho) - M_{s_1,s_2}(\rho)^2\right) + O_p\left(\frac{1}{m^2}\right).$$

Similarly,

$$\text{Cov}\left(\sum_{i=1}^{m}\frac{\hat{\lambda}_{W,i}^{s_1}}{\left(\rho-\hat{\lambda}_{W,i}\right)^{s_2}}u_{1,i}^2, \sum_{j=1}^{m}\frac{\hat{\lambda}_{W,i}^{r_1}}{\left(\rho-\hat{\lambda}_{W,i}\right)^{r_2}}u_{1,j}^2\right)$$
$$=\frac{2}{m}\left(M_{s_1+r_1,s_2+r_2}(\rho)-M_{s_1,s_2}(\rho)M_{r_1,r_2}(\rho)\right)+O_p\left(\frac{1}{m^2}\right),$$

$$\text{Var}\left(\sum_{i=1}^{m}\frac{\hat{\lambda}_{W,i}^{s_1}}{\left(\rho-\hat{\lambda}_{W,i}\right)^{s_2}}u_{1,i}u_{2,i}\right)=\frac{M_{2s_1,2s_2}(\rho)-M_{s_1,s_2}(\rho)^2}{m}+O_p\left(\frac{1}{m^2}\right),$$

$$\text{Cov}\left(\sum_{i=1}^{m}\frac{\hat{\lambda}_{W,i}^{s_1}}{\left(\rho-\hat{\lambda}_{W,i}\right)^{s_2}}u_{1,i}u_{2,i}, \sum_{j=1}^{m}\frac{\hat{\lambda}_{W,i}^{r_1}}{\left(\rho-\hat{\lambda}_{W,i}\right)^{r_2}}u_{1,j}u_{2,i}\right)$$
$$=\frac{M_{s_1+r_1,s_2+r_2}(\rho)-M_{s_1,s_2}(\rho)M_{r_1,r_2}(\rho)}{m}+O_p\left(\frac{1}{m^2}\right).$$

In order to prove normality, we study projections. Let

$$\left(\sqrt{m}\sum_{i=1}^{m}a_{r,s}\left(\hat{\lambda}_{W,i}\right)\left(u_{s,i}u_{r,i}-\mathbf{1}_{r=s}\frac{1}{m}\right)\right),$$

where $a_{r,s}$ are real-valued functions and $r,s$ take values between 1 and $m$. We need to show that sums on $r$ and $s$ of these statistics are still normal. Furthermore, let $v_1, v_2,..., v_{p_2} \in \mathbb{R}^m$ ($p_2 \geqslant p_1$) be independent normal random vectors with $\text{Var}(v_i)=\frac{1}{m}\mathbf{I}_m$. Applying Gram-Schmidt we obtain

$$\begin{aligned}u_s &= \frac{v_s-\sum_{i=1}^{s-1}\left(\langle v_s,v_i\rangle+O_p\left(\frac{1}{m}\right)\right)\frac{v_i}{||v_i||^2+O_p\left(\frac{1}{m}\right)}}{||v_s||+O_p\left(\frac{1}{m}\right)}\\ &= \sum_{i=1}^{s}b_{s,i}v_i.\end{aligned}$$

The central limit theorem of Bentkus (2005) shows that for any $t_1,t_2=1,2,...k$,

$$\sqrt{m}\sum_{i=1}^{m}a_{r,s}\left(\hat{\lambda}_{W,i}\right)\left(v_{t_1,i}v_{t_2,i}-\mathbf{1}_{t_1=t_2}\frac{1}{m}\right) \tag{2.1}$$

can for finite $k$ be jointly approximated by a centred multivariate normal. In particular, using Slutsky's theorem we show that

$$\forall s,t=1,2,...,k,\ s\neq t,\ \sqrt{m}\,b_{s,t} \tag{2.2}$$

tends to centred joint Normal random variable with variance of size $O_p(1)$. Moreover,

$$\forall s=1,2,...,k,\ \sqrt{m}\,(b_{s,s}-1) \tag{2.3}$$

tends to a centred joint normal random variable with variance of size $O_p(1)$.

Finally, we consider the elements in the statistics $B_m$.

$$\sum_{s,r=1}^{k} \sqrt{m} \sum_{i=1}^{m} a_{r,s}\left(\hat{\lambda}_{W,i}\right) \left(u_{s,i} u_{r,i} - \mathbf{1}_{s=r} \frac{1}{m}\right)$$

$$= \sum_{s,r=1}^{k} \sum_{t_1=1}^{s} \sum_{t_2=1}^{r} b_{s,t_1} b_{r,t_2} \sqrt{m} \sum_{i=1}^{m} a_{r,s}\left(\hat{\lambda}_{W,i}\right) v_{t_1,i} v_{t_2,i} - \frac{1}{m} \sum_{s=1}^{k} \sqrt{m} \sum_{i=1}^{m} a_{s,s}\left(\hat{\lambda}_{W,i}\right)$$

$$= \sum_{s,r=1}^{k} \sum_{t_1=1}^{s} \sum_{t_2=1}^{r} b_{s,t_1} b_{r,t_2} \sqrt{m} \sum_{i=1}^{m} a_{r,s}\left(\hat{\lambda}_{W,i}\right) \left(v_{t_1,i} v_{t_2,i} - \mathbf{1}_{t_1=t_2} \frac{1}{m}\right)$$

$$+ \sum_{s,r=1}^{k} \sum_{t=1}^{\min(r,s)} b_{s,t} b_{r,t} \sqrt{m} \frac{1}{m} \sum_{i=1}^{m} a_{r,s}\left(\hat{\lambda}_{W,i}\right) - \sum_{s=1}^{k} \sqrt{m} \frac{1}{m} \sum_{i=1}^{m} a_{s,s}\left(\hat{\lambda}_{W,i}\right)$$

$$= \underbrace{\sum_{s,r=1}^{k} \sum_{t_1=1}^{s} \sum_{t_2=1}^{r} b_{s,t_1} b_{r,t_2} \sqrt{m} \sum_{i=1}^{m} a_{r,s}\left(\hat{\lambda}_{W,i}\right) \left(v_{t_1,i} v_{t_2,i} - \mathbf{1}_{t_1=t_2} \frac{1}{m}\right)}_{A_1}$$

$$+ \underbrace{\sum_{\substack{s,r=1 \\ s \neq r}}^{k} \sum_{t=1}^{\min(r,s)} b_{s,t} b_{r,t} \sqrt{m} \frac{1}{m} \sum_{i=1}^{m} a_{r,s}\left(\hat{\lambda}_{W,i}\right)}_{A_2}$$

$$- \underbrace{\sum_{s=1}^{k} \sqrt{m} \left(1 - b_{s,s}^2\right) \frac{1}{m} \sum_{i=1}^{m} a_{s,s}\left(\hat{\lambda}_{W,i}\right)}_{A_3}$$

We show asymptotic normality of the three elements using Slutsky's Theorem:

$$\begin{aligned} A_1 &= b_{s,t_1} b_{r,t_2} \sqrt{m} \sum_{i=1}^{m} a_{r,s}\left(\hat{\lambda}_{W,i}\right) \left(v_{t_1,i} v_{t_2,i} - \mathbf{1}_{t_1=t_2} \frac{1}{m}\right) \\ &= \mathrm{E}\left[b_{s,t_1} b_{r,t_2}\right] \sqrt{m} \sum_{i=1}^{m} a_{r,s}\left(\hat{\lambda}_{W,i}\right) \left(v_{t_1,i} v_{t_2,i} - \mathbf{1}_{t_1=t_2} \frac{1}{m}\right) + O_p\left(\frac{1}{\sqrt{m}}\right). \end{aligned}$$

We recognise the equation 2.1 on the right and we know that we can neglect all the terms when $s \neq t_1$ or $r \neq t_2$. Indeed in these cases by equation 2.2, $b_{s,t_1} b_{r,t_2} = O_p(1/\sqrt{m})$ and the terms become of order $O_p(1/\sqrt{m})$.

$$\begin{aligned} A_2 &= \sum_{\substack{s,r=1 \\ s \neq r}}^{k} \sum_{t=1}^{\min(r,s)} b_{s,t} b_{r,t} \sqrt{m} \frac{1}{m} \sum_{i=1}^{m} a_{r,s}\left(\hat{\lambda}_{W,i}\right) \\ &= \sum_{\substack{s,r=1 \\ s \neq r}}^{k} b_{s,\min(r,s)} b_{r,\min(r,s)} \sqrt{m} \frac{1}{m} \sum_{i=1}^{m} a_{r,s}\left(\hat{\lambda}_{W,i}\right) + \sum_{\substack{s,r=1 \\ s \neq r}}^{k} \sum_{t=1}^{\min(r,s)-1} b_{s,t} b_{r,t} \sqrt{m} \frac{1}{m} \sum_{i=1}^{m} a_{r,s}\left(\hat{\lambda}_{W,i}\right) \\ &= \sum_{\substack{s,r=1 \\ s \neq r}}^{k} \mathrm{E}\left[b_{\min(r,s),\min(r,s)}\right] \sqrt{m} b_{\max(r,s),\min(r,s)} \frac{1}{m} \sum_{i=1}^{m} a_{r,s}\left(\hat{\lambda}_{W,i}\right) + O_p\left(\frac{1}{\sqrt{m}}\right) \\ &= \sum_{\substack{s,r=1 \\ s \neq r}}^{k} \sqrt{m} b_{\max(r,s),\min(r,s)} \frac{1}{m} \sum_{i=1}^{m} a_{r,s}\left(\hat{\lambda}_{W,i}\right) + O_p\left(\frac{1}{\sqrt{m}}\right). \end{aligned}$$

We recognise the equation 2.2 multiplied by a constant.

$$A_3 = \sum_{s=1}^{k} \sqrt{m} \left(1 - b_{s,s}^2\right) \frac{1}{m} \sum_{i=1}^{m} a_{s,s}\left(\hat{\lambda}_{W,i}\right).$$

We recognise up to a constant Equation 2.3. Therefore, we obtain:

$$\sum_{s,r=1}^{k} \sqrt{m} \sum_{i=1}^{m} a_{r,s}\left(\hat{\lambda}_{W,i}\right)\left(u_{s,i}u_{r,i} - \mathbf{1}_{s=r}\frac{1}{m}\right)$$

$$= \sum_{s,r=1}^{k} \sqrt{m} \sum_{i=1}^{m} a_{r,s}\left(\hat{\lambda}_{W,i}\right)\left(v_{s,i}v_{r,i} - \mathbf{1}_{s=r}\frac{1}{m}\right) + \sum_{\substack{s,r=1 \\ s\neq r}}^{k} \sqrt{m} b_{\min(s,t),\max(s,t)} \frac{1}{m} \sum_{i=1}^{m} a_{r,s}\left(\hat{\lambda}_{W,i}\right)$$

$$+ \sum_{s=1}^{k} 2\sqrt{m}(b_{s,s} - 1)) \frac{1}{m} \sum_{i=1}^{m} a_{s,s}\left(\hat{\lambda}_{W,i}\right).$$

Because all the quantities in the three sums are approximately jointly normal, the quantity is approximately normal with an error $o_p(1)$ for large $m$. Assuming $s_i \leqslant r_i$, if $(s_1, r_1) \neq (s_2, r_2)$, then, when either $s_1 \neq r_1$ or $s_2 \neq r_2$, by invariance under rotation,

$$\operatorname{Cov}\left(\sum_{i=1}^{m} a_{r_1,s_1}\left(\hat{\lambda}_{W,i}\right)\left(u_{s_1,i}u_{r_1,i} - \mathbf{1}_{r_1=s_1}\frac{1}{m}\right), \sum_{i=1}^{m} a_{r_2,s_2}\left(\hat{\lambda}_{W,i}\right)\left(u_{s_2,i}u_{r_2,i} - \mathbf{1}_{r_2=s_2}\frac{1}{m}\right)\right) = 0.$$

Therefore, when $(s_1, r_1) \neq (s_2, r_2)$, the resulting joint Normal statistics are asymptotically independent. However, when $(s_1, r_1) = (s_2, r_2)$, then the statistics are jointly Normal and correlated.
Note that when $s_1 = r_1$ and $s_2 = r_2$ we have to use the linearity of the covariance and the fact that

$$\operatorname{Corr}\left(\sum_{i=1}^{m} a_i u_{1,i}^2, \sum_{i=1}^{m} \tilde{a}_i u_{2,i}^2\right) = \operatorname{Corr}\left(\sum_{i=1}^{m} a_i u_{1,i}^2, \sum_{i=1}^{m} \tilde{a}_i \frac{1 - u_{1,i}^2}{m - 1}\right)$$

that tends to 0. □

### 2.3. Characterisation and convergence of eigenvalues and angles

Parts 2.a and 2.b of the following theorem are due to Benaych-Georges and Rao (2009). This paper provides also the main idea of the remaining parts of the theorem.

**Theorem 2.3.**
Let $P = I_m + (\theta - 1)uu^t$ be a finite perturbation of order 1.

1. Suppose $W$ is a symmetric random matrix with eigenvalues $\hat{\lambda}_{W,i} \geqslant 0$ and eigenvectors $\hat{u}_{W,i}$ for $i = 1, 2, ..., m$. The perturbation of $W$ by $P$ leads to $\hat{\Sigma} = P^{1/2}WP^{1/2}$.
   For $i = 1, 2, ..., m$, we define $\tilde{u}_{\hat{\Sigma},i}$ and $\hat{\lambda}_{\hat{\Sigma},i}$ such that

   $$WP\tilde{u}_{\hat{\Sigma},i} = \hat{\lambda}_{\hat{\Sigma},i}\tilde{u}_{\hat{\Sigma},i},$$

   and the usual $\hat{u}_{\hat{\Sigma},i}$ such that if $\hat{\Sigma} = P^{1/2}WP^{1/2}$, then

   $$\hat{\Sigma}\hat{u}_{\hat{\Sigma},i} = P^{1/2}WP^{1/2}\hat{u}_{\hat{\Sigma},i} = \hat{\lambda}_{\hat{\Sigma},i}\hat{u}_{\hat{\Sigma},i}.$$

- The eigenvalues $\hat{\lambda}_{\hat{\Sigma},s}$ are such that for $s = 1, 2, ..., m$,

$$\sum_{i=1}^{m} \frac{\hat{\lambda}_{W,i}}{\hat{\lambda}_{\hat{\Sigma},s} - \hat{\lambda}_{W,i}} \langle \hat{u}_{W,i}, u \rangle^2 = \frac{1}{\theta_k - 1}.$$

- The eigenvectors $\tilde{u}_{\hat{\Sigma},s}$ are such that

$$\left\langle \tilde{u}_{\hat{\Sigma},s}, v \right\rangle^2 = \frac{\left( \sum_{i=1}^{m} \frac{\hat{\lambda}_{W,i}}{\hat{\lambda}_{\hat{\Sigma},s} - \hat{\lambda}_{W,i}} \langle \hat{u}_{W,i}, v \rangle \langle \hat{u}_{W,i}, u \rangle \right)^2}{\sum_{i=1}^{m} \frac{\hat{\lambda}_{W,i}^2}{(\hat{\lambda}_{\hat{\Sigma},s} - \hat{\lambda}_{W,i})^2} \langle \hat{u}_{W,i}, u \rangle^2}.$$

In particular if $v = u$,

$$\left\langle \tilde{u}_{\hat{\Sigma},s}, u \right\rangle^2 = \frac{1}{(\theta_k - 1)^2 \left( \sum_{i=1}^{m} \frac{\hat{\lambda}_{W,i}^2}{(\hat{\lambda}_{\hat{\Sigma},s} - \hat{\lambda}_{W,i})^2} \langle \hat{u}_{W,i}, u \rangle^2 \right)}.$$

Moreover,

$$\hat{u}_{\hat{\Sigma},s} = \frac{P^{1/2} \tilde{u}_{\hat{\Sigma},s}}{\sqrt{1 + (\theta - 1) \left\langle \tilde{u}_{\hat{\Sigma},s}, u \right\rangle^2}}.$$

Therefore, for $u$ and $v$ such that $\langle v, u \rangle = 0$,

$$\left\langle \hat{u}_{\hat{\Sigma},s}, u \right\rangle^2 = \frac{\theta \left\langle \tilde{u}_{\hat{\Sigma},s}, u \right\rangle^2}{1 + (\theta - 1) \left\langle \tilde{u}_{\hat{\Sigma},s}, u \right\rangle^2} = -\frac{\theta}{(\theta - 1)^2 \hat{\lambda}_{\hat{\Sigma},s} T'_{W,u}(\hat{\lambda}_{\hat{\Sigma},s})},$$

$$\left\langle \hat{u}_{\hat{\Sigma},s}, v \right\rangle^2 = \frac{\left\langle \tilde{u}_{\hat{\Sigma},s}, u \right\rangle^2}{1 + (\theta - 1) \left\langle \tilde{u}_{\hat{\Sigma},s}, u \right\rangle^2},$$

where $T_{W,u}(z) = \sum_{i=1}^{m} \frac{\hat{\lambda}_{W,i}}{z - \hat{\lambda}_{W,i}} \langle \hat{u}_{W,i}, u \rangle^2$.

2. Suppose that $W_X$, $W_Y$ and $P = P_X = P_Y$ satisfy Assumption 1.1 and that $\theta$ is large enough to create detectable spikes, $\left( \hat{\theta}_X, \hat{u}_X \right)$ and $\left( \hat{\theta}_Y, \hat{u}_Y \right)$, in the matrices $\hat{\Sigma}_X = P^{1/2} W_X P^{1/2}$ and $\hat{\Sigma}_Y = P^{1/2} W_Y P^{1/2}$. Then,

$$a) \quad \hat{\hat{\theta}}_X, \hat{\hat{\theta}}_Y \xrightarrow[n,m \to \infty]{P} \theta,$$

$$b) \quad \langle \hat{u}_X, u \rangle - \alpha_X, \langle \hat{u}_Y, u \rangle - \alpha_Y \xrightarrow[n,m \to \infty]{P} 0,$$

$$c) \quad \langle \hat{u}_X, \hat{u}_Y \rangle - \alpha_X \alpha_Y \xrightarrow[n,m \to \infty]{P} 0,$$

where

$$\hat{\theta}_X \xrightarrow[n,m\to\infty]{P} \rho_X,$$

$$\hat{\hat{\theta}}_X = 1 + \frac{1}{\hat{T}_{\hat{\Sigma}_X}(\hat{\theta}_X)} = 1 + \frac{m}{\sum_{i=k+1}^{m} \frac{\hat{\lambda}_{\hat{\Sigma}_X,i}}{\hat{\theta}_X - \hat{\lambda}_{\hat{\Sigma}_X,i}}},$$

$$\alpha_X^2 = -\frac{\theta}{(\theta-1)^2 \rho_X \hat{T}'_{W_X}(\rho_X)},$$

$$\alpha_Y^2 = -\frac{\theta}{(\theta-1)^2 \rho_Y \hat{T}'_{W_Y}(\rho_Y)},$$

$\hat{\lambda}_{\hat{\Sigma}_X,i}$ and $\hat{\lambda}_{\hat{\Sigma}_Y,i}$ are the eigenvalues of respectively $\hat{\Sigma}_X$ and $\hat{\Sigma}_Y$.

*Remark* 2.3.1.
If $W_X$ and $W_Y$ are Wishart random matrices of size $m$ and degree of freedom $n_X$, $n_Y$, respectively, then by setting $c_X = \frac{m}{n_X}$ and $c_Y = \frac{m}{n_Y}$

$$\alpha_X^2 = \frac{1 - \frac{c_X}{(\theta-1)^2}}{1 + \frac{c_X}{\theta-1}},$$

$\hat{\hat{\theta}}_X$ is such that $\hat{\theta}_X = \hat{\hat{\theta}}_X \left(1 + \frac{c_X}{\hat{\hat{\theta}}_X - 1}\right)$, and

$$\lim_{m\to\infty} \hat{\theta}_X = \theta \left(1 + \frac{c_X}{\theta - 1}\right).$$

*Proof.* **Theorem 2.3**

In the proof of this theorem, we use the two transforms:

- $T_{W_X,u}(z) = \sum_{i=1}^{m} \frac{\hat{\lambda}_{W_X,i}}{z - \hat{\lambda}_{W_X,i}} \langle \hat{u}_{W_X,i}, u \rangle^2$,
- $\hat{T}_{W_X}(z) = \frac{1}{m} \sum_{i=1}^{m} \frac{\hat{\lambda}_{W_X,i}}{z - \hat{\lambda}_{W_X,i}}$.

1. To study $P^{1/2} W P^{1/2}$, note the spectral decomposition, $W = \hat{U}_W^t \Lambda_W \hat{U}_W$.

   **Eigenvalues**   The above matrix has the same eigenvalues as $WP$. Thus, for $z$ an eigenvalue of $P^{1/2} W P^{1/2}$,

$$\begin{aligned}
0 &= \det\left(z I_m - P^{1/2} W P^{1/2}\right) \\
&= \det\left(z I_m - W \left(I_m + (\theta-1) u u^t\right)\right) \\
&= \det\left(z I_m - \hat{U}_W^t \Lambda_W \hat{U}_W \left(I_m + (\theta-1) u u^t\right)\right) \\
&= \det\left(z I_m - \Lambda_W \left(I_m + (\theta-1) \left(\hat{U}_W^t u\right)\left(\hat{U}_W^t u\right)^t\right)\right) \\
&= \det(z I_m - \Lambda_W) \det\left(I_m - (z I_m - \Lambda_W)^{-1} \Lambda_W \left((\theta-1) \left(\hat{U}_W^t u\right)\left(\hat{U}_W^t u\right)^t\right)\right).
\end{aligned}$$

If $z$ is not an eigenvalue of $W$ but an eigenvalue of $WP$, then

$$\det(z I_m - \Lambda_W) \neq 0.$$

Thus,

$$\det\left(I_m - (zI - \Lambda_W)^{-1}\Lambda_W\left((\theta-1)\left(\hat{U}_W^t u\right)\left(\hat{U}_W^t u\right)^t\right)\right) = 0$$

$$\Rightarrow \text{Trace}\left((zI_m - \Lambda_W)^{-1}\Lambda_W\left((\theta-1)\left(\hat{U}_W^t u\right)\left(\hat{U}_W^t u\right)^t\right)\right) = 1$$

$$\Rightarrow \text{Trace}\left((zI_m - \Lambda_W)^{-1}\Lambda_W\left(\left(\hat{U}_W^t u\right)\left(\hat{U}_W^t u\right)^t\right)\right) = \frac{1}{\theta-1}$$

$$\Rightarrow \sum_{i=1}^{m} \frac{\hat{\lambda}_{W,i}}{z - \hat{\lambda}_{W,i}} \langle \hat{u}_{W,i}, u \rangle^2 = \frac{1}{\theta-1}.$$

In our notation, we can replace $z$ by the eigenvalues of $\hat{\Sigma} = P^{1/2}WP^{1/2}$, $\hat{\lambda}_{\hat{\Sigma},s}$ for $s = 1, 2, ..., m$.

**Eigenvectors**  We first study $\tilde{u}_{\hat{\Sigma},s}$ for $s = 1, 2, ..., m$, the eigenvectors of $WP$. Because $\hat{\lambda}_{\hat{\Sigma},s}$ is the eigenvalue of $WP$ corresponding to $\tilde{u}_{\hat{\Sigma},s}$,

$$\hat{\lambda}_{\hat{\Sigma},s}\tilde{u}_{\hat{\Sigma},s} = WP\tilde{u}_{\hat{\Sigma},s} = W(I_m + \bar{P})\tilde{u}_{\hat{\Sigma},s} = (W + W\bar{P})\tilde{u}_{\hat{\Sigma},s},$$

where $\bar{P} = (\theta-1)uu^t$. Therefore, we have

$$\tilde{u}_{\hat{\Sigma},s} = ((\theta-1)u^t\tilde{u}_{\hat{\Sigma},s})\left(\hat{\lambda}_{\hat{\Sigma},s}I_m - W\right)^{-1}Wu.$$

Using the fact that $\tilde{u}_{\hat{\Sigma},s}$ is proportional to $\left(\hat{\lambda}_{\hat{\Sigma},s}I_m - W\right)^{-1}Wu$ and its norm is 1, we get

$$\tilde{u}_{\hat{\Sigma},s} = \frac{\left(\hat{\lambda}_{\hat{\Sigma},s}I_m - W\right)^{-1}W}{\sqrt{u^tW\left(\hat{\lambda}_{\hat{\Sigma},s}I_m - W\right)^{-2}Wu}}u,$$

which leads to most of the results about eigenvectors we need.
Assuming $v \in \mathbf{R}^m$, we have

$$\left\langle \tilde{u}_{\hat{\Sigma},s}, v \right\rangle^2 = \frac{\left(v^tW\left(\hat{\lambda}_{\hat{\Sigma},s}I_m - W\right)^{-1}u\right)^2}{u^tW\left(\hat{\lambda}_{\hat{\Sigma},s}I_m - W\right)^{-2}W_X u}$$

$$= \frac{\left(\sum_{i=1}^{m} \frac{\hat{\lambda}_{W,i}}{\hat{\lambda}_{\hat{\Sigma},s}-\hat{\lambda}_{W,i}} \langle \hat{u}_{W,i}, v \rangle \langle \hat{u}_{W,i}, u \rangle\right)^2}{\sum_{i=1}^{m} \frac{\hat{\lambda}_{W,i}^2}{(\hat{\lambda}_{P,s}-\hat{\lambda}_{W,i})^2} \langle \hat{u}_{W,i}, u \rangle^2}.$$

In particular, when $v = u$, $u^tW_X(\hat{\lambda}_{\hat{\Sigma},s}I - W_X)^{-1}u = 1/(\theta-1)$ by the previous equation on eigenvalues.

In order to obtain a more elegant formula, we need to study the denominator and its transform:

$$\int \frac{x^2}{(z-x)^2} f(x)dx = -\int \frac{x}{(z-x)} f(x)dx + z \int \frac{x}{(z-x)^2} f(x)dx = -T(z) - zT'(z),$$

where empirical version of the T-transform $T$ is defined above. It follows that

$$\left\langle \tilde{u}_{\hat{\Sigma},s}, u \right\rangle^2 = \frac{T_{W,u}(\hat{\lambda}_{\hat{\Sigma},s})^2}{-T_{W,u}(\hat{\lambda}_{\hat{\Sigma},s}) - \hat{\lambda}_{\hat{\Sigma},s} T'_u(\hat{\lambda}_{\hat{\Sigma},s})}$$

$$= -\frac{1}{(\theta-1)^2 \hat{\lambda}_{\hat{\Sigma},s} T'_{W,u}(\hat{\lambda}_{\hat{\Sigma},s}) + (\theta-1)}.$$

It is easy to see that

$$\hat{u}_{\hat{\Sigma},s} = \frac{P^{1/2} \tilde{u}_{\hat{\Sigma},s}}{\sqrt{\tilde{u}^t_{\hat{\Sigma},s} P \tilde{u}_{\hat{\Sigma},s}}}.$$

Using the fact that $\tilde{u}^t_{\hat{\Sigma},s} P \tilde{u}_{\hat{\Sigma},s} = 1 + (\theta-1)\left\langle \tilde{u}_{\hat{\Sigma},s}, u \right\rangle^2$, the angle can be written,

$$\left\langle \hat{u}_{\hat{\Sigma},s}, v \right\rangle^2 = \frac{\left\langle P^{1/2} \tilde{u}_{\hat{\Sigma},s}, v \right\rangle^2}{\tilde{u}^t_{\hat{\Sigma},s} P \tilde{u}_{\hat{\Sigma},s}}$$

$$= \frac{\left\langle \left(I_m + (\sqrt{\theta}-1)uu^t\right) \tilde{u}_{\hat{\Sigma},s}, v \right\rangle^2}{1 + (\theta-1)\left\langle \tilde{u}_{\hat{\Sigma},s}, u \right\rangle^2}$$

$$= \frac{\left(\left\langle \tilde{u}_{\hat{\Sigma},s}, v \right\rangle + \left(\sqrt{\theta}-1\right)\left\langle \tilde{u}_{\hat{\Sigma},s}, u \right\rangle \langle u, v \rangle\right)^2}{1 + (\theta-1)\left\langle \tilde{u}_{\hat{\Sigma},s}, u \right\rangle^2}.$$

For $u$ and $v$ such that $\langle v, u \rangle = 0$, we easily obtain,

$$\left\langle \hat{u}_{\hat{\Sigma},s}, u \right\rangle^2 = \frac{\theta \left\langle \tilde{u}_{\hat{\Sigma},s}, u \right\rangle^2}{1 + (\theta-1)\left\langle \tilde{u}_{\hat{\Sigma},s}, u \right\rangle^2}$$

$$= -\frac{\theta}{(\theta-1)^2 \hat{\lambda}_{\hat{\Sigma},s} T'_{W,u}(\hat{\lambda}_{\hat{\Sigma},s})},$$

$$\left\langle \hat{u}_{\hat{\Sigma},s}, v \right\rangle^2 = \frac{\left\langle \tilde{u}_{\hat{\Sigma},s}, u \right\rangle^2}{1 + (\theta-1)\left\langle \tilde{u}_{\hat{\Sigma},s}, u \right\rangle^2}.$$

2. The second part of the theorem is a direct consequence of the first part. Because $W_X$ and $W_Y$ satisfy Assumption 1.1 and $P$ is assumed detectable, the leading eigenvalues $\hat{\lambda}_{\hat{\Sigma}_X,1} = \hat{\theta}_X$ and $\hat{\lambda}_{\hat{\Sigma}_Y,1} = \hat{\theta}_Y$ are the spikes.

(a) We have

$$\frac{1}{\theta - 1} = T_u(\hat{\theta}_X) = \sum_{i=1}^{m} \frac{\hat{\lambda}_{W_X,i}}{\hat{\theta}_X - \hat{\lambda}_{W_X,i}} \langle \hat{u}_{W,i}, u \rangle^2,$$

where $\langle \hat{u}_{W_X,i}, u \rangle^2 = w_i$ creates an unit uniform vector $w$ independent of the spectrum of $W_X$.
Because

$$\frac{1}{\theta - 1} - \frac{1}{\hat{\theta}_X - 1} = T_u(\hat{\theta}_X) - \frac{1}{m}\sum_{i=1}^{m} \frac{\hat{\lambda}_{W_X,i}}{\hat{\theta}_X - \hat{\lambda}_{W_X,i}} \xrightarrow[n,m\to\infty]{P} 0,$$

we obtain $\hat{\theta}_X \xrightarrow[n,m\to\infty]{P} \theta$. In order to replace $\hat{\lambda}_{W_X,i}$ by $\hat{\lambda}_{\hat{\Sigma}_X,i}$, we make use of Cauchy's interlacing law and the boundedness of the eigenvalues.

(b) Moreover, assuming that $\hat{\theta}_X \xrightarrow[n,m\to\infty]{P} \rho_X$, we find

$$\langle \hat{u}_X, u \rangle^2 - \alpha_X^2 = -\frac{\theta}{(\theta-1)^2 \hat{\theta}_X T_u'(\hat{\theta}_X)} + \frac{\theta}{(\theta-1)^2 \rho_X T_{W_X}'(\rho_X)} \xrightarrow[n,m\to\infty]{P} 0.$$

(c) Finally, an elementary result from linear algebra proves the limit of $\langle \hat{u}_X, \hat{u}_Y \rangle$, which we call the double angle. If $W = O^t \Lambda O$, with invariant uniform $O$ independent of the bounded spectrum $\Lambda$, it follows that for any two fixed vectors $v, u \in \mathbf{R}^m$,

$$v^t W u = \tilde{v}^t \Lambda \tilde{u} = \sum_{i=1}^{m} \hat{\lambda}_{W,i} \tilde{v}_i \tilde{u}_i = \langle u, v \rangle \operatorname{Trace}(W)/m + O_p(1/m).$$

This result is easily demonstrated, because $\tilde{v}, \tilde{u}$ are uniform, $\langle \tilde{v}, \tilde{u} \rangle = \langle v, u \rangle$ and it is easy to show that

$$\mathrm{E}\left[\tilde{v}_i \tilde{u}_i\right] = \mathrm{E}\left[\sum_{s,r=1}^{m} O_{i,s} O_{i,r} v_s v_r\right] = \langle v, u \rangle /m,$$

$$\operatorname{Var}\left(\sum_{i=1}^{m} \tilde{v}_i \tilde{u}_i\right) = 0 \text{ and } \operatorname{Var}(\tilde{v}_i \tilde{u}_i) = O_p(1/m^2),$$

$$\operatorname{Cov}(\tilde{v}_i \tilde{u}_i, \tilde{v}_j \tilde{u}_j) = O_p(1/m^3).$$

Using the notation $\hat{\theta}_X = \hat{\lambda}_{\hat{\Sigma}_X,1}$ and $\hat{\theta}_Y = \hat{\lambda}_{\hat{\Sigma}_Y,1}$ for the estimated eigenvalues, we can combine the first part and the above result to obtain

$$\langle \tilde{u}_X, \tilde{u}_Y \rangle^2 = \frac{\left(u^t W_X(\hat{\theta}_X I_m - W_X)^{-1}(\hat{\theta}_Y I_m - W_Y)^{-1} W_Y u\right)^2}{u^t W_X(\hat{\theta}_X I - W_X)^{-2} W_X W_Y(\hat{\theta}_Y I - W_Y)^{-2} W_Y u}$$

$$= \frac{\left(u^t W_X(\hat{\theta}_X I_m - W_X)^{-1}(\hat{\theta}_Y I_m - W_Y)^{-1} W_Y u\right)^2}{\left(-T_u\left(\hat{\theta}_X\right) - \hat{\theta}_X T'_u\left(\hat{\theta}_X\right)\right)\left(-T_u\left(\hat{\theta}_Y\right) - \hat{\theta}_Y T'_u\left(\hat{\theta}_Y\right)\right)}$$

$$= \frac{\left(\text{Trace}\left(W_X(\hat{\theta}_X I_m - W_X)^{-1}\right) u^t(\hat{\theta}_Y I_m - W_Y)^{-1} W_Y u + O_p\left(\frac{1}{m\theta}\right)\right)^2}{\left(-T_u\left(\hat{\theta}_X\right) - \hat{\theta}_X T'_u\left(\hat{\theta}_X\right)\right)\left(-T_u\left(\hat{\theta}_Y\right) - \hat{\theta}_Y T'_u\left(\hat{\theta}_Y\right)\right)}$$

$$= \langle u, \tilde{u}_X \rangle^2 \langle u, \tilde{u}_Y \rangle^2 + O_p\left(\frac{1}{m\theta}\right),$$

where $\tilde{u}_X$ and $\tilde{u}_Y$ are the eigenvectors of $W_X P$ and $W_Y P$, respectively. As in the first part, we use the relation between $\hat{u}_X$ and $\tilde{u}_X$ to get

$$\langle \hat{u}_X, \hat{u}_Y \rangle^2 = \frac{\left\langle \left(I_m + \bar{P}\right)^{1/2} \tilde{u}_X, \left(I_m + \bar{P}\right)^{1/2} \tilde{u}_Y \right\rangle^2}{\left(\tilde{u}_X^t (I_m + \bar{P})\tilde{u}_X\right)\left(\tilde{u}_Y^t (I_m + \bar{P})\tilde{u}_Y\right)}$$

$$= \frac{\left(\langle \tilde{u}_X, u \rangle \langle u, \tilde{u}_Y \rangle \theta + \overbrace{(\langle \tilde{u}_X, \tilde{u}_Y \rangle - \langle \tilde{u}_X, u \rangle \langle u, \tilde{u}_Y \rangle)}^{\to 0}\right)^2}{\left(1 + (\theta - 1)\langle \tilde{u}_X, u \rangle^2\right)\left(1 + (\theta - 1)\langle \tilde{u}_Y, u \rangle^2\right)}.$$

This leads to the result,

$$\langle \hat{u}_X, \hat{u}_Y \rangle^2 - \alpha_X^2 \alpha_Y^2 \xrightarrow[n,m\to\infty]{P} 0.$$

$\square$

### 2.4. Asymptotic distribution of the eigenvalues and the angle

**Theorem 2.4.**
*Suppose $W_X$ and $W_Y$ are independent random matrices satisfying 1.1 with $P = P_X = P_Y$, a detectable perturbation of order $k = 1$. Moreover, we assume $S_{W_X} = \left\{\hat{\lambda}_{W_X,1}, \hat{\lambda}_{W_X,2}, ..., \hat{\lambda}_{W_X,m}\right\}$ and $S_{W_Y} = \left\{\hat{\lambda}_{W_Y,1}, \hat{\lambda}_{W_Y,2}, ..., \hat{\lambda}_{W_Y,m}\right\}$, the eigenvalues of $W_X$ and $W_Y$ as known. Let*

$$\hat{\Sigma}_X = P^{1/2} W_X P^{1/2},$$
$$\hat{\Sigma}_Y = P^{1/2} W_Y P^{1/2},$$
$$P = I_m + (\theta - 1)uu^t,$$

*where $u$ is fixed. We construct the unbiased estimators of $\theta$,*

$$\hat{\hat{\theta}}_X \Bigg|_{\frac{1}{\hat{\theta}_X - 1} = \frac{1}{m}\sum_{i=1}^m \frac{\hat{\lambda}_{W_X,i}}{\hat{\theta}_X - \hat{\lambda}_{W_X,i}}} \quad \text{and} \quad \hat{\hat{\theta}}_Y \Bigg|_{\frac{1}{\hat{\theta}_Y - 1} = \frac{1}{m}\sum_{i=1}^m \frac{\hat{\lambda}_{W_Y,i}}{\hat{\theta}_Y - \hat{\lambda}_{W_Y,i}}}$$

*where $\hat{\theta}_X = \hat{\lambda}_{\hat{\Sigma}_X,1}$ and $\hat{\theta}_Y = \hat{\lambda}_{\hat{\Sigma}_Y,1}$ are the largest eigenvalues of $\hat{\Sigma}_X$ and $\hat{\Sigma}_Y$ corresponding to the eigenvectors $\hat{u}_X = \hat{u}_{\hat{\Sigma}_X,1}$ and $\hat{u}_Y = \hat{u}_{\hat{\Sigma}_Y,1}$.*

1. If $\frac{\theta}{\sqrt{m}} \to 0$, we define

$$M_{s,r,X}(\rho_X) = \frac{1}{m}\sum_{i=1}^m \frac{\hat{\lambda}_{W_X,i}^s}{(\rho_X - \hat{\lambda}_{W_X,i})^r}, M_{s,r,Y}(\rho_Y) = \frac{1}{m}\sum_{i=1}^m \frac{\hat{\lambda}_{W_X,i}^s}{(\rho_Y - \hat{\lambda}_{W_X,i})^r}$$

where we assume

$$\rho_X = \mathrm{E}\left[\hat{\theta}_X\right] + o\left(\frac{\theta}{\sqrt{m}}\right), \quad \rho_Y = \mathrm{E}\left[\hat{\theta}_Y\right] + o\left(\frac{\theta}{\sqrt{m}}\right)$$

and a convergence rate in probability of $(\hat{\theta}_X, \hat{\theta}_Y)$ to $(\rho_X, \rho_Y)$ in $1/\sqrt{m}$. Then

$$\begin{pmatrix} \hat{\hat{\theta}}_X \\ \langle \hat{u}_X, u\rangle^2 \end{pmatrix}\bigg| S_{W_X} \sim \mathbf{N}\left(\begin{pmatrix} \theta \\ \alpha_X^2 \end{pmatrix}, \frac{1}{m}\begin{pmatrix} \sigma_{\theta,X}^2 & \sigma_{\theta,\alpha^2,X} \\ \sigma_{\theta,\alpha^2,X} & \sigma_{\alpha^2,X}^2 \end{pmatrix}\right) + \begin{pmatrix} o_p\left(\frac{\theta}{\sqrt{m}}\right) \\ o_p\left(\frac{1}{\theta\sqrt{m}}\right) \end{pmatrix},$$

$$\begin{pmatrix} \hat{\hat{\theta}}_X \\ \hat{\hat{\theta}}_Y \\ \langle \hat{u}_X, \hat{u}_Y\rangle^2 \end{pmatrix}\bigg| S_{W_X}, S_{W_Y} \sim \mathbf{N}\left(\begin{pmatrix} \theta \\ \theta \\ \alpha_{X,Y}^2 \end{pmatrix}, \frac{1}{m}\begin{pmatrix} \sigma_{\theta,X}^2 & 0 & \sigma_{\theta,\alpha^2,X} \\ 0 & \sigma_{\theta,Y}^2 & \sigma_{\theta,\alpha^2,Y} \\ \sigma_{\theta,\alpha^2,X} & \sigma_{\theta,\alpha^2,Y} & \sigma_{\alpha^2,X,Y}^2 \end{pmatrix}\right) + \begin{pmatrix} o_p\left(\frac{\theta}{\sqrt{m}}\right) \\ o_p\left(\frac{\theta}{\sqrt{m}}\right) \\ o_p\left(\frac{1}{\theta\sqrt{m}}\right) \end{pmatrix},$$

where

$$\alpha_X^2 = \frac{\theta}{(\theta-1)^2}\frac{1}{\rho_X M_{1,2,X}},$$

$$\alpha_{X,Y}^2 = \frac{\theta^2}{(\theta-1)^4}\frac{1}{\rho_X \rho_Y M_{1,2,X} M_{1,2,Y}},$$

$$\sigma_{\theta,X}^2 = \frac{2\left(M_{2,2,X} - M_{1,1,X}^2\right)}{M_{1,1,X}^4},$$

$$\sigma_{\alpha^2,X}^2 = \frac{2\theta^2}{((\theta-1)\rho_X M_{1,2,X})^4}\left(\rho_X^2\left(M_{2,4,X} - M_{1,2,X}^2\right) + \left(2\rho_X \frac{M_{1,3,X}}{M_{1,2,X}} - 1\right)^2\left(M_{2,2,X} - M_{1,1,X}^2\right)\right.$$
$$\left. -2\rho_X\left(2\rho_X \frac{M_{1,3,X}}{M_{1,2,X}} - 1\right)\left(M_{2,3,X} - \frac{M_{1,1,X}}{M_{1,2,X}}\right)\right),$$

$$\sigma_{\theta,\alpha^2,X} = \frac{2\theta}{M_{1,1,X}^2 M_{1,2,X}^3 \rho_X^2(-1+\theta)^2}\Big(M_{1,1,X}M_{1,2,X}^2\rho_X + 2M_{1,3,X}M_{2,2,X}\rho_X$$
$$+M_{1,1,X}^2(M_{1,2,X} - 2M_{1,3,X}\rho_X) - M_{1,2,X}(M_{2,2,X} + M_{2,3,X}\rho_X)\Big),$$

$$\sigma_{\alpha^2,X,Y}^2 = \sigma_{\alpha^2,X}^2\alpha_Y^4 + \sigma_{\alpha^2,Y}^2\alpha_X^4 + 4\alpha_{X,Y}^2(1-\alpha_X^2)(1-\alpha_X^2).$$

2. If $\frac{\theta}{\sqrt{m}} \to \infty$ (Assumption 1.2 (A1)), we can simplify the above as follows. Let

$$M_{r,X} = \frac{1}{m}\sum_{i=1}^m \hat{\lambda}_{W_X,i}^r \quad \text{and} \quad M_{r,Y} = \frac{1}{m}\sum_{i=1}^m \hat{\lambda}_{W_Y,i}^r.$$

We then have

$$\begin{pmatrix} \hat{\hat{\theta}}_X \\ \langle \hat{u}_X, u \rangle^2 \end{pmatrix} \bigg| S_{W_X} \sim \mathbf{N}\left( \begin{pmatrix} \theta + O_p(1) \\ 1 + \frac{1 - M_{2,X}}{\theta} + O_p\left(\frac{1}{\theta^2}\right) \end{pmatrix}, \right.$$

$$\left. \frac{1}{m}\begin{pmatrix} 2\theta^2(1 - M_{2,X}) & 2(2M_{2,X}^2 - M_{2,X} - M_{3,X}) \\ 2(2M_{2,X}^2 - M_{2,X} - M_{3,X}) & \frac{2}{\theta^2}(4M_{2,X}^3 - M_{2,X}^2 - 4M_{2,X}M_{3,X} + M_{4,X}) \end{pmatrix} \right) + \begin{pmatrix} O_p\left(\frac{\theta}{\sqrt{m}}\right) \\ O_p\left(\frac{1}{\theta\sqrt{m}}\right) \end{pmatrix},$$

$$\begin{pmatrix} \hat{\hat{\theta}}_X \\ \hat{\hat{\theta}}_Y \\ \langle \hat{u}_X, \hat{u}_Y \rangle^2 \end{pmatrix} \bigg| S_{W_X}, S_{W_Y} \sim \mathbf{N}\left( \begin{pmatrix} \theta + O_p(1) \\ \theta + O_p(1) \\ 1 + \frac{2 - M_{2,X} - M_{2,Y}}{\theta} + O_p\left(\frac{1}{\theta^2}\right) \end{pmatrix}, \right.$$

$$\left. \frac{1}{m}\begin{pmatrix} 2\theta^2(1 - M_{2,X}) & 0 & 2(2M_{2,X}^2 - M_{2,X} - M_{3,X}) \\ 0 & 2\theta^2(1 - M_{2,Y}) & 2(2M_{2,Y}^2 - M_{2,Y} - M_{3,Y}) \\ 2(2M_{2,X}^2 - M_{2,X} - M_{3,X}) & 2(2M_{2,Y}^2 - M_{2,Y} - M_{3,Y}) & \frac{S}{\theta^2} \end{pmatrix} \right) + \begin{pmatrix} O_p\left(\frac{\theta}{\sqrt{m}}\right) \\ O_p\left(\frac{\theta}{\sqrt{m}}\right) \\ O_p\left(\frac{1}{\theta\sqrt{m}}\right) \end{pmatrix},$$

$$S = 2\left(4M_{2,X}^3 - M_{2,X}^2 - 4M_{2,X}M_{3,X} + M_{4,X}\right) + 2\left(4M_{2,Y}^3 - M_{2,Y}^2 - 4M_{2,Y}M_{3,Y} + M_{4,Y}\right) + 4(M_{2,Y} - 1)(M_{2,X} - 1).$$

*Moreover, the asymptotic distributions of $\hat{\theta}_X$ and $\hat{\hat{\theta}}_X$ are the same.*

3. *If $\frac{\theta}{\sqrt{m}} \to d$, a finite positive constant, then a mixture of the two first scenarios holds.*

   *The formula of the second moment is asymptotically the same as the variance formula when $\frac{\theta}{\sqrt{m}} \to \infty$, whereas the formula of the first moment is asymptotically the same as the expectation formula when $\frac{\theta}{\sqrt{m}} \to 0$.*

4. *The random variables can be expressed as functions of invariant unit random statistics of the form:*

$$\overset{u}{M}_{r,s,X}(\rho) = \sum_{i=1}^{m} \frac{\hat{\lambda}_{W_X,i}^r}{\left(\rho - \hat{\lambda}_{W_X,i}\right)^s} \langle \hat{u}_{W_X,i}, u \rangle^2.$$

   *(Assuming a canonical perturbation ($u = e_1$) leads to a simpler formula.) Knowing $S_{W_X}$ and $S_{W_Y}$ we have*

   - *Exact distributions:*

$$\hat{\theta}_X \text{ satisfies } \frac{1}{\theta - 1} = \overset{u}{M}_{1,1,X}\left(\hat{\theta}_X\right),$$

$$\langle \hat{u}_X, e_1 \rangle^2 = \frac{\theta}{(\theta - 1)^2} \frac{1}{\hat{\theta}_X \overset{u}{M}_{1,2,X}\left(\hat{\theta}_X\right)}.$$

   *Moreover,*

$$\langle \hat{u}_X, \hat{u}_Y \rangle = \langle \hat{u}_X, e_1 \rangle \langle \hat{u}_Y, e_1 \rangle + \sqrt{1 - \langle \hat{u}_X, e_1 \rangle^2}\sqrt{1 - \langle \hat{u}_Y, e_1 \rangle^2} Z,$$

$$\sum_{i=2}^{m} \hat{u}_{X,1} \hat{u}_{Y,1} = \sqrt{1 - \langle \hat{u}_X, e_1 \rangle^2}\sqrt{1 - \langle \hat{u}_Y, e_1 \rangle^2} Z,$$

$$\text{where } Z \sim \mathbf{N}\left(0, \frac{1}{m}\right) + O_p\left(\frac{1}{m}\right).$$

   *The random variable $Z$ is independent of $\langle \hat{u}_X, e_1 \rangle$, $\langle \hat{u}_Y, e_1 \rangle$, $\hat{\theta}_X$ and $\hat{\theta}_Y$. In order to get the exact marginal distribution of $Z$, we could replace it by $\sum_{i=1}^{m-1} v_i \tilde{v}_i$ with $v_i \tilde{v}_i$ being independent unit invariant random vectors.*

- *Approximations:*

$$\hat{\theta}_X = \rho + \frac{\left(\overset{u}{M}_{1,1,X}(\rho) - M_{1,1,X}(\rho)\right)}{M_{1,2,X}(\rho)} + O_p\left(\frac{\theta}{m}\right)$$

$$= \theta \overset{u}{M}_{1,X} + O_p(1),$$

$$\hat{\hat{\theta}}_X = \theta + (\theta-1)^2 \left(\overset{u}{M}_{1,1,X}(\rho) - M_{1,1,X}(\rho_X)\right) + O_p\left(\frac{\theta}{m}\right).$$

We provide three estimators of the angle in order to cover all behaviors of $\theta$:

$$\langle \hat{u}_X, e_1 \rangle^2 = \frac{\theta}{(\theta-1)^2}\left(\frac{1}{\rho_X M_{1,2,X}(\rho_X)} + \left(\frac{2M_{1,3,X}(\rho_X)}{M_{1,2,X}(\rho_X)} - \frac{1}{\rho_X}\right)\frac{\overset{u}{M}_{1,1,X}(\rho_X) - \frac{1}{\theta-1}}{(M_{1,2,X}(\rho_X))^2}\right.$$

$$\left. -\frac{\overset{u}{M}_{1,2,X}(\rho_X) - M_{1,2,X}(\rho_X)}{\rho_X(M_{1,2,X}(\rho_X))^2}\right) + O_p\left(\frac{1}{m}\right),$$

$$= 1 + \frac{1}{\theta}\left(1 - \overset{u}{M}_{2,X} + 2M_{2,X}\left(\overset{u}{M}_{1,X} - 1\right)\right)$$

$$+ \frac{1}{\theta^2}\left(1 - 2\overset{u}{M}_{2,X} + 3\overset{u}{M}_{2,X}^2 - 2\overset{u}{M}_{3,X}\right) + O_p\left(\frac{1}{\theta^3}\right) + O_p\left(\frac{1}{\theta m}\right),$$

$$= 1 + \frac{1}{\theta} - \frac{1}{\theta}\overset{u}{M}_{2,X} + \frac{2}{\theta}M_{2,X}\left(\overset{u}{M}_{1,X} - 1\right) + O_p\left(\frac{1}{\theta^2}\right) + O_p\left(\frac{1}{\theta m}\right)$$

Finally, the double angle is such that

$$\langle \hat{u}_X, \hat{u}_Y \rangle = \langle \hat{u}_X, e_1 \rangle \langle \hat{u}_Y, e_1 \rangle + \frac{\sqrt{M_{2,X}-1}\sqrt{M_{2,Y}-1}}{\theta}Z + O_p\left(\frac{1}{\theta^2\sqrt{m}}\right).$$

*Remark* 2.4.1. If the spectra of $W_X$ and $W_Y$ are those of rescaled Wishart matrices of size $m$ with $n$ degree of freedom, we obtain the following results. Let $c = \frac{m}{n}$,

$$\begin{pmatrix} \hat{\hat{\theta}}_X \\ \langle \hat{u}_X, u_0 \rangle^2 \end{pmatrix} \overset{Asy}{\sim} \mathbf{N}\left(\begin{pmatrix} \theta \\ 1 - \frac{\frac{c}{(\theta-1)^2}}{1+\frac{c}{\theta-1}} \end{pmatrix}, \frac{1}{m}\begin{pmatrix} -\frac{2c(\theta-1)^2\theta^2}{c-(\theta-1)^2} & -\frac{2c^2(\theta-1)\theta^3}{(c-(\theta-1)^2)(c+\theta-1)^2} \\ -\frac{2c^2(\theta-1)\theta^3}{(c-(\theta-1)^2)(c+\theta-1)^2} & -\frac{2c^2\theta^2\left(c^2+(\theta(\theta+2)-2)c+(\theta-1)^2\right)}{(c-(\theta-1)^2)(c+\theta-1)^4} \end{pmatrix}\right)$$

and

$$\begin{pmatrix} \hat{\hat{\theta}}_X \\ \hat{\hat{\theta}}_Y \\ \langle \hat{u}_X, \hat{u}_Y \rangle^2 \end{pmatrix} \overset{Asy}{\sim} \mathbf{N}\left(\begin{pmatrix} \theta \\ \theta \\ \left(1-\frac{\frac{c}{(\theta-1)^2}}{1+\frac{c}{\theta-1}}\right)^2 \end{pmatrix}, \right.$$

$$\left. \frac{1}{m}\begin{pmatrix} -\frac{2c(\theta-1)^2\theta^2}{c-(\theta-1)^2} & 0 & -\frac{2c^2(\theta-1)\theta^3}{(c-(\theta-1)^2)(c+\theta-1)^2} \\ 0 & -\frac{2c(\theta-1)^2\theta^2}{c-(\theta-1)^2} & -\frac{2c^2(\theta-1)\theta^3}{(c-(\theta-1)^2)(c+\theta-1)^2} \\ -\frac{2c^2(\theta-1)\theta^3}{(c-(\theta-1)^2)(c+\theta-1)^2} & -\frac{2c^2(\theta-1)\theta^3}{(c-(\theta-1)^2)(c+\theta-1)^2} & \frac{4c^2\theta^2(c-(\theta-1)^2)^2\left(c^3+4c^2(\theta-1)+c(\theta-1)(\theta(\theta+5)-5)+2(\theta-1)^3\right)}{(\theta-1)^4(c+\theta-1)^7} \end{pmatrix}\right).$$

If $\theta$ tends to infinity, then

$$\begin{pmatrix} \hat{\hat{\theta}}_X \\ \langle \hat{u}_X, u \rangle^2 \end{pmatrix} \overset{Asy}{\sim} \mathbf{N}\left(\begin{pmatrix} \theta \\ 1-\frac{\frac{c}{(\theta-1)^2}}{1+\frac{c}{\theta-1}} \end{pmatrix}, \frac{1}{m}\begin{pmatrix} 2c\theta^2 & 2c^2 \\ 2c^2 & \frac{2c^2(c+1)}{\theta^2} \end{pmatrix}\right).$$

Moreover,

$$\begin{pmatrix} \hat{\hat{\theta}}_X \\ \hat{\hat{\theta}}_Y \\ \langle \hat{u}_X, \hat{u}_X \rangle^2 \end{pmatrix} \overset{Asy}{\sim} \mathbf{N}\left(\begin{pmatrix} \theta \\ \theta \\ 1-\frac{\frac{c}{(\theta-1)^2}}{1+\frac{c}{\theta-1}} \end{pmatrix}, \frac{1}{m}\begin{pmatrix} 2c\theta^2 & 0 & 2c^2 \\ 0 & 2c\theta^2 & 2c^2 \\ 2c^2 & 2c^2 & \frac{4c^2(c+2)}{\theta^2} \end{pmatrix}\right).$$

Rémy Mariétan and Stephan Morgenthaler/Comparison of two populations   19*Proof.* **Theorem 2.4** The first three parts of the theorem are consequences of the fourth part. The computation of the particular case assuming a fast convergence of the spectrum to the Marcenko-Pastur distribution is left to the reader.

**Prerequisite discussion**   Using the property of invariance by rotation applied to $P^{-1/2}\mathbf{X}$ and $P^{-1/2}\mathbf{Y}$ (see Assumption 1.1), we can assume without loss of generality that $P$ is canonical. Therefore, $u = \epsilon_1$.

In this proof, we use the following notation

$$\overset{u}{M}_{r,s,X}(\rho) = \sum_{i=1}^{m} \frac{\hat{\lambda}_{W_X,i}^r}{\left(\rho - \hat{\lambda}_{W_X,i}\right)^s} \hat{u}_{W_X,i,1}^2,$$

$$M_{r,s,X}(\rho) = \sum_{i=1}^{m} \frac{\hat{\lambda}_{W_X,i}^r}{\left(\rho - \hat{\lambda}_{W_X,i}\right)^s},$$

$$M_{r,s,\hat{\Sigma}_X}(\rho) = \sum_{i=2}^{m} \frac{\hat{\lambda}_{\hat{\Sigma}_X,i}^r}{\left(\rho - \hat{\lambda}_{\hat{\Sigma}_X,i}\right)^s}.$$

In particular, $\overset{u}{M}_{1,1,X}(\rho) = T_{W_X,\epsilon_1}(\rho)$ and $M_{1,1,X}(\rho) = T_{W_X}(\rho)$.

**Characterisation**   Using the characterisations presented in Benaych-Georges and Rao (2009) we have

$$\theta = 1 + \frac{1}{\overset{u}{M}_{1,1,X}(\hat{\theta}_X)}$$
$$= 1 + \frac{1}{\sum_{i=1}^{m} \frac{\hat{\lambda}_{W_X,i}}{\hat{\theta}_X - \hat{\lambda}_{W_X,i}} \hat{u}_{W_X,i,1}^2}$$

and

$$\langle \hat{u}_X, u \rangle^2 = \frac{\theta}{(\theta-1)^2 \hat{\theta}_X \overset{u}{M}_{1,2,X}(\hat{\theta}_X)}.$$

This result is generalised in Theorem 2.3 (see proof on page 11).

**Unit invariant vector statistics for eigenvalues**   Based on the condition that $\mathrm{E}\left[\hat{\theta}_X\right] = \rho_X + o\left(\frac{1}{\sqrt{m}}\right)$ and $\hat{\theta}_X = \rho_X + O_p\left(\frac{1}{\sqrt{m}}\right)$, there exists $\tilde{\theta}_X \in \left[\min\left(\hat{\theta}_X, \rho_X\right), \max\left(\hat{\theta}_X, \rho_X\right)\right]$ such that

$$\frac{1}{\theta-1} = \sum_{i=1}^{m} \left( \frac{\hat{\lambda}_{W_X,i}}{\rho_X - \hat{\lambda}_{W_X,i}} \hat{u}_{W_X,i,1}^2 - \frac{\hat{\lambda}_{W_X,i}}{\left(\rho_X - \hat{\lambda}_{W_X,i}\right)^2} \hat{u}_{W_X,i,1}^2 \left(\hat{\theta}_X - \rho_X\right) \right.$$
$$\left. + \frac{\hat{\lambda}_{W_X,i}}{\left(\tilde{\theta}_X - \hat{\lambda}_{W_X,i}\right)^3} \hat{u}_{W_X,i,1}^2 \left(\hat{\theta}_X - \rho_X\right)^2 \right).$$

Thus,

$$\hat{\theta}_X = \rho_X + \frac{\overset{u}{M}_{1,1,X}(\rho_X) - \frac{1}{\theta-1}}{\overset{u}{M}_{1,2,X}(\rho_X)} + \frac{\overset{u}{M}_{1,3,X}(\tilde{\theta}_X)}{\overset{u}{M}_{1,2,X}(\tilde{\theta}_X)} \left(\hat{\theta}_X - \rho_X\right)^2, \text{ where } \tilde{\theta}_X \sim \rho_X,$$

$$= \rho_X + \frac{\overset{u}{M}_{1,1,X}(\rho) - \frac{1}{\theta-1}}{M_{1,2,X}(\rho)} - \frac{\overset{u}{M}_{1,1,X}(\rho_X) - \frac{1}{\theta-1}}{M_{1,2,X}(\rho_X)^2} \left(\overset{u}{M}_{1,2,X}(\rho) - M_{1,2,X}(\rho_X)\right) + O_p\left(\frac{\theta}{m}\right).$$

Next, we show that $\overset{u}{M}_{1,1,X}(\rho_X) - \frac{1}{\theta-1} = O_p\left(\frac{1}{\sqrt{m\theta}}\right).$

$$\overset{u}{M}_{1,1,X}(\rho_X) - \frac{1}{\theta-1} = \overset{u}{M}_{1,1,X}(\rho_X) - \overset{u}{M}_{1,1,X}(\hat{\theta}_X)$$
$$= \overset{u}{M}_{1,2,X}(\rho_X) \left(\hat{\theta}_X - \rho_X\right) + O_p\left(\frac{1}{m\theta}\right)$$
$$= O_p\left(\frac{1}{\sqrt{m\theta}}\right).$$

A similar argument shows that

$$\left(\overset{u}{M}_{1,2,X}(\rho_X) - M_{1,2,X}(\rho_X)\right) = O_p\left(\frac{1}{\sqrt{m\theta^2}}\right).$$

This leads to the first notable equation:

$$\hat{\theta}_X = \rho_X + \frac{\overset{u}{M}_{1,1,X}(\rho_X) - \frac{1}{\theta-1}}{M_{1,2,X}(\rho_X)} + O_p\left(\frac{\theta}{m}\right) \qquad (2.4)$$

The unbiased eigenvalue satisfies

$$\hat{\hat{\theta}}_X - 1 = \frac{1}{M_{1,2,X}(\hat{\theta}_X)}$$
$$= \frac{m}{\sum_{i=1}^{m} \frac{\hat{\lambda}_{W_X,i}}{\hat{\theta}_X - \hat{\lambda}_{W_X,i}}}$$
$$= \frac{1}{M_{1,1,X}(\rho_X) - M_{1,2,X}(\rho_X)\left(\hat{\theta}_X - \rho_X\right) + O_p\left(\frac{1}{\theta m}\right)}$$
$$= \frac{1}{M_{1,1,X}(\rho_X) - \overset{u}{M}_{1,1,X}(\rho_X) + \frac{1}{\theta-1}} + O_p\left(\frac{\theta}{m}\right)$$
$$= \theta - 1 - (\theta-1)^2 \left(M_{1,1,X}(\rho_X) - \overset{u}{M}_{1,1,X}(\rho_X)\right) + O_p\left(\frac{\theta}{m}\right).$$

This leads to the second notable equation:

$$\hat{\theta}_X = \theta + (\theta-1)^2 \left( \overset{u}{M}_{1,1,X}(\rho_X) - M_{1,1,X}(\rho_X) \right) + O_p\left(\tfrac{\theta}{m}\right). \tag{2.5}$$

When $\theta$ is large we have

$$\frac{1}{\theta-1} = \frac{1}{\hat{\theta}_X} \sum_{i=1}^m \hat{\lambda}_{W_X,i} \hat{u}^2_{W_X,i,1} + \frac{1}{\hat{\theta}_X^2} \sum_{i=1}^m \hat{\lambda}^2_{W_X,i} \hat{u}^2_{W_X,i,1} + O_p\left(\tfrac{1}{\theta^3}\right),$$

which shows that

$$\hat{\theta}_X = (\theta-1) \overset{u}{M}_{1,X} + \frac{\overset{u}{M}_{2,X}}{M_{1,X}} + O_p\left(\tfrac{1}{m}\right) + O_p\left(\tfrac{1}{\theta}\right) \tag{2.6}$$

$$= \theta \overset{u}{M}_{1,X} + O_p(1). \tag{2.7}$$

**Unit invariant vector statistics for angles**   The following expansions are useful:

1: Applying Taylor expansions around $\rho_X$ as computed previously for $\hat{\theta}$ (2.4) and a similar one for $\overset{u}{M}_{1,2,X}(\hat{\theta}_X)$ lead to

$$\langle \hat{u}_X, u \rangle^2 = \frac{\theta}{(\theta-1)^2 \hat{\theta}_X \overset{u}{M}_{1,2,X}(\hat{\theta}_X)}$$

$$= \frac{\theta}{(\theta-1)^2} \left( \frac{1}{\rho_X M_{1,2,X}(\rho_X)} + \left( \frac{2M_{1,3,X}(\rho_X)}{M_{1,2,X}(\rho_X)} - \frac{1}{\rho_X} \right) \frac{\overset{u}{M}_{1,1,X}(\rho_X) - \frac{1}{\theta-1}}{(M_{1,2,X}(\rho_X))^2} \right.$$

$$\left. - \frac{\overset{u}{M}_{1,2,X}(\rho_X) - M_{1,2,X}(\rho_X)}{\rho_X (M_{1,2,X}(\rho_X))^2} \right) + O_p\left(\tfrac{1}{m}\right).$$

2: For large $\theta$ (2.6) and using the fact that $\overset{u}{M}_{1,1}(\hat{\theta}_X) = 1/(\theta-1)$, leads to

$$\langle \hat{u}_X, u \rangle^2 = 1 + \frac{1}{\theta} - \frac{\theta}{\hat{\theta}_X^2} \overset{u}{M}_{2,X} + O_p\left(\tfrac{1}{\theta^2}\right)$$

$$= 1 + \frac{1}{\theta} - \frac{1}{\theta} \frac{\overset{u}{M}_{2,X}}{\left(\overset{u}{M}_{1,X}\right)^2} + O_p\left(\tfrac{1}{\theta^2}\right)$$

$$= 1 + \frac{1}{\theta} - \frac{1}{\theta} \overset{u}{M}_{2,X} + \frac{2}{\theta} M_{2,X} \left( \overset{u}{M}_{1,X} - 1 \right) + O_p\left(\tfrac{1}{\theta^2}\right) + O_p\left(\tfrac{1}{\theta m}\right).$$

3: Finally, using a higher order approximation of $\hat{\theta}_X$ (2.7) and $\overset{u}{M}_{1,2,X}(\hat{\theta}_X)$ when $\theta$ is large shows that

$$\langle \hat{u}_X, u \rangle^2 = 1 + \frac{1}{\theta} \left( 1 - \overset{u}{M}_{2,X} + 2M_{2,X} \left( \overset{u}{M}_{1,X} - 1 \right) \right)$$

$$+ \frac{1}{\theta^2} \left( 1 - 2\overset{u}{M}_{2,X} + 3\overset{u}{M}_{2,X}^2 - 2\overset{u}{M}_{3,X} \right) + O_p\left(\tfrac{1}{\theta^3}\right) + O_p\left(\tfrac{1}{\theta m}\right).$$

The study of the angle $\langle \hat{u}_X, u \rangle^2$ is separated into several cases, depending on the behavior of $\theta$,

| $\theta$ finite | | $\frac{\theta^3}{\sqrt{m}} \to$ constant | | $\frac{\theta^2}{\sqrt{m}} \to$ constant | | $\frac{\theta^2}{m} \to$ constant | |
|---|---|---|---|---|---|---|---|
| A | B | C | D | E | F | G | H |

B, D, F and H are intermediate order compared to the neighbours.

**A**   An application of the delta method using the first expansion of the angle and Equation (2.4) leads to the result. In this case the standard error is of order $1/\sqrt{m}$; therefore, the error of order $1/m$ is negligible.

**B through F**   The first expansion shows that the error is not important. Indeed assuming that $\theta \to \infty$ makes the first expansion applicable which shows that
$$\mathrm{Var}\left(\langle \hat{u}_X, u \rangle^2\right) \sim \frac{1}{\theta^2 m}.$$

Moreover, $1 - \mathrm{E}\left[\langle \hat{u}_X, u \rangle^2\right]$ is of order $1/\theta$. If $\epsilon \sim 1/m$ is the error of the first expansion, then

$$\frac{\epsilon}{1/\theta} \sim \frac{\frac{1}{m}}{\frac{1}{\theta}} = \frac{\theta}{m} \to 0,$$

$$\frac{\epsilon}{\sqrt{\mathrm{Var}\left(\langle \hat{u}_X, u \rangle^2\right)}} \sim \frac{\frac{1}{m}}{\frac{1}{\theta\sqrt{m}}} = \frac{\theta}{\sqrt{m}} \to 0.$$

**D-E-F**   The previous arguments are still valid in these cases. Alternatively, we can prove the result using the third expansion of the angle and Equation (2.6). Thus, the estimations are asymptotically equivalent.

**G**   This scenario is the most difficult and only the third approximation of the angle and Equation (2.6) work in this case. However,

- The expectations obtained by expansions 1 and 3 of the angle lead to the same result.
- The variances obtained by expansions 2 and 3 of the angle lead to the same result.

The equivalence between the variance is direct, but the expectation requires some computations which we leave to the reader.

**H**   Using expansion 2 of the angle and Equation (2.7) leads to a negligible error compared to the standard error and the expectation.

This study provides the behaviour of the statistic for all $\theta$

- If $\frac{\theta}{\sqrt{m}} \to 0$, then we use the first expansion of the angle and Equation (2.5).
- If $\frac{\theta}{\sqrt{m}} \to \infty$, then we use the second expansion of the angle and Equation (2.7).

- If $\frac{\theta}{\sqrt{m}} \to d$, where $d$ is a fixed constant, then we use the first expansion with Equation (2.5) to estimate the expectation and the second expansion with Equation (2.7) to estimate the variance.

**Unit invariant vector statistics for double angle** Let $\hat{u}_X$ and $\hat{u}_Y$ be the first eigenvectors of $\hat{\Sigma}_X$ and $\hat{\Sigma}_Y$, respectively, then without loss of generality, if $u = e_1$,

$$\begin{aligned}
\langle \hat{u}_X, \hat{u}_Y \rangle &= \hat{u}_{X,1}\hat{u}_{Y,1} + \sum_{i=2}^{m} \hat{u}_{X,i}\hat{u}_{Y,i} \\
&= \hat{u}_{X,1}\hat{u}_{Y,1} + \sqrt{1 - \hat{u}_{X,1}^2}\sqrt{1 - \hat{u}_{Y,1}^2}Z \\
&= \hat{u}_{X,1}\hat{u}_{Y,1} + \sqrt{1 - \alpha_X^2}\sqrt{1 - \alpha_Y^2}\left(\tilde{Z}/\sqrt{m} + O_p\left(\frac{1}{m}\right)\right),
\end{aligned}$$

where $\alpha_X^2$ and $\alpha_Y^2$ are the limits of $\hat{u}_{X,1}^2$ and $\hat{u}_{Y,1}^2$. Because the eigenvectors are invariant by rotations that do not affect the first component, it follows that $Z$ is asymptotic independent of $\hat{u}_{X,1}$ and $\hat{u}_{Y,1}$. Moreover, $Z$ is the scalar product of two independent unit invariant vectors of dimension $m-1$ and can be estimated by a standard Normal $\tilde{Z}$ divided by $\sqrt{m}$.

**Joint distribution** Applying the delta method to the estimation of $\hat{\theta}_X$, $\hat{\theta}_Y$, $\langle \hat{u}_X, \hat{u}_Y \rangle$, $\langle \hat{u}_X, u \rangle$, $\langle \hat{u}_Y, u \rangle$ and $Z$ leads directly to the joint distribution.
The computations are lengthy but straightforward and based on the separate consideration of the three terms with rates of $\theta/\sqrt{m}$ as above. The details are left to the reader. □

### 2.5. Residual spikes

**Lemma 2.1.**
*Suppose*

$$D = (I_m + (\theta - 1)u_X u_X^t)^{-1/2} (I_m + (\theta - 1)u_Y u_Y^t) (I_m + (\theta - 1)u_X u_X^t)^{-1/2}.$$

*The eigenvalues of $D$ are equal to $1$ and*

$$\lambda(D) = -\frac{1}{2\theta}\left(-1 + \alpha^2 - 2\alpha^2\theta - \theta^2(1 - \alpha^2) \pm \sqrt{-4\theta^2 + [1 + \theta^2 - (-1 + \theta)^2\alpha^2]^2}\right),$$

*where $\alpha^2 = \langle u_X, u_Y \rangle^2$.*
*Moreover, for*

$$D_2 = (I_m + (\theta_X - 1)u_X u_X^t)^{-1/2} (I_m + (\theta_Y - 1)u_Y u_Y^t) (I_m + (\theta_X - 1)u_X u_X^t)^{-1/2}.$$

*the eigenvalues are $1$ and*

$$\lambda(D_2) = \frac{1}{2}\left(\theta_Y + \alpha^2 - \theta_Y\alpha^2 + \frac{1 + (\theta_Y - 1)\alpha^2 \pm \sqrt{-4\theta_Y\theta_X + (1 + \theta_Y\theta_X - (\theta_Y - 1)(\theta_X - 1)\alpha^2)^2}}{\theta_X}\right).$$

*Proof.* Lemma 2.1

For any set of $m-2$ orthogonal vectors of dimension $m$, orthogonal to $u_X$ and $u_Y$, the eigenvalues are 1. Consequently, we just need to compute the two remaining eigenvalues, $\lambda_1$ and $\lambda_2$. We will consider the Trace of $D_2$ and $D_2^2$.

$$
\begin{aligned}
\text{Trace}(D_2) &= \text{Trace}\left(\left(I_m + (\theta_Y - 1)u_Y u_Y^t\right)\left(I_m + \left(\frac{1}{\theta_X} - 1\right) u_X u_X^t\right)\right) \\
&= m + (\theta_Y - 1) + \left(\frac{1}{\theta_X} - 1\right) + (\theta_Y - 1)\left(\frac{1}{\theta_X} - 1\right)\langle u_X, u_Y\rangle^2 \\
&= m - 2 + \langle u_X, u_Y\rangle^2 + \theta_Y\left(1 - \langle u_X, u_Y\rangle^2\right) + \frac{1 + \langle u_Y, u_Y\rangle^2(\theta_Y - 1)}{\theta_X},
\end{aligned}
$$

$$
\begin{aligned}
\text{Trace}(D_2^2) &= \text{Trace}\left(\left(\left(I_m + (\theta_Y - 1)u_Y u_Y^t\right)\left(I_m + \left(\frac{1}{\theta_X} - 1\right) u_X u_X^t\right)\right)^2\right) \\
&= m - 2 + \langle u_X, u_Y\rangle^4 + \theta_Y^2\left(\langle u_X, u_Y\rangle^2 - 1\right)^2 - 2\theta_Y\langle u_X, u_Y\rangle^2\left(\langle u_X, u_Y\rangle^2 - 1\right) \\
&\quad + \frac{\left(1 + (\theta_Y - 1)\langle u_X, u_Y\rangle^2\right)^2}{\theta_X^2} - \frac{2(\theta_X - 1)^2\langle u_X, u_Y\rangle^2\left(\langle u_X, u_Y\rangle^2 - 1\right)}{\theta_X}.
\end{aligned}
$$

Since

$$
\begin{aligned}
\text{Trace}(D_2) &= m - 2 + \lambda_1 + \lambda_2, \\
\text{Trace}(D_2^2) &= m - 2 + \lambda_1^2 + \lambda_2^2,
\end{aligned}
$$

the result follows by solving the two equations. □

## References

BENAYCH-GEORGES, F. and RAO, N. R. (2009). The eigenvalues and eigenvectors of finite, low rank perturbations of large random matrices. *Advances in Mathematics* **227** 494-521.

BENTKUS, V. (2005). A Lyapunov type bound in $\mathbf{R}^d$. *Theory of Probability and Its Applications* **49** 311–323.

MARIÉTAN, R. and MORGENTHALER, S. (2020). Statistical applications of Random matrix theory: comparison of two populations I. *arXiv:2002.12741*.



\end{document}